\documentclass[a4paper,10pt,twoside]{article}

\usepackage[utf8]{inputenc}
\usepackage{amsmath, amssymb, amsthm}
\usepackage{amstext, amsfonts, a4}
\usepackage{dsfont}
\usepackage{latexsym}
\usepackage{mathrsfs}
\usepackage[all,cmtip]{xy}
\usepackage{graphicx}
\usepackage{enumerate}
\usepackage{hyperref,color}
\usepackage{stmaryrd}
\usepackage[shortlabels]{enumitem}
\usepackage{bbm}

\usepackage{url}
\usepackage[pdftex,a4paper,left=3cm,right=3cm,top=3cm,bottom=3cm]{geometry}

\def\ra{\rightarrow}

\def\lra{\longrightarrow}

\def\lmapsto{\longmapsto}

\def\sA{\mathscr{A}} \def\sB{\mathscr{B}}

 \def\sS{\mathscr{S}} 
  
 \def\sZ{\mathscr{Z}}

\def\bbA{\mathbb{A}}\def\bbC{\mathbb{C}}
\def\bbF{\mathbb{F}}\def\bbG{\mathbb{G}}

\def\bbN{\mathbb{N}}
\def\bbQ{\mathbb{Q}}\def\bbR{\mathbb{R}}\def\bbT{\mathbb{T}}
\def\bbX{\mathbb{X}}
\def\bbZ{\mathbb{Z}}

\def\cA{\mathcal{A}}\def\cC{\mathcal{C}}
\def\cH{\mathcal{H}}
\def\cI{\mathcal{I}}
\def\cM{\mathcal{M}}\def\cO{\mathcal{O}}
\def\cS{\mathcal{S}}
\def\cW{\mathcal{W}}
\def\cZ{\mathcal{Z}}

\def\bfB{\mathbf{B}}
\def\bfG{\mathbf{G}}

\def\bfT{\mathbf{T}}

\def\bfq{\mathbf{q}}

\def\whB{\widehat{B}}
\def\whG{\widehat{G}}

\def\whT{\widehat{T}}
\def\whU{\widehat{U}}

\def\ts{\tilde{s}}

\DeclareMathOperator{\aff}{aff}

\DeclareMathOperator{\Sph}{Sph}

\DeclareMathOperator{\BP}{BP}

\DeclareMathOperator{\Ch}{Ch}

\DeclareMathOperator{\cl}{cl}

\DeclareMathOperator{\diag}{diag}
\DeclareMathOperator{\Diag}{Diag}

\DeclareMathOperator{\End}{End}

\DeclareMathOperator{\ev}{ev}

\DeclareMathOperator{\Her}{Her}

\DeclareMathOperator{\Id}{Id}

\DeclareMathOperator{\Ind}{Ind}

\DeclareMathOperator{\Infl}{Infl}

\DeclareMathOperator{\Mat}{Mat}
\DeclareMathOperator{\Mod}{Mod}

\DeclareMathOperator{\nonreg}{non-reg}

\DeclareMathOperator{\pr}{pr}

\DeclareMathOperator{\QCoh}{QCoh}

\DeclareMathOperator{\rel}{rel}
\DeclareMathOperator{\reg}{reg}
\DeclareMathOperator{\Rep}{Rep}
\DeclareMathOperator{\Res}{Res}

\DeclareMathOperator{\sign}{sign}
\DeclareMathOperator{\Sing}{Sing}
\DeclareMathOperator{\SingDiag}{SingDiag}
\DeclareMathOperator{\SingUpTriang}{SingUpTriang}
\DeclareMathOperator{\Spec}{Spec}
\DeclareMathOperator{\sph}{sph}

\DeclareMathOperator{\SP}{SP}

\DeclareMathOperator{\sss}{ss}

\DeclareMathOperator{\Tr}{Tr}

\DeclareMathOperator{\UpTriang}{UpTriang}

\def\lan{\langle}

\def\ran{\rangle}

\newtheorem{counter}[subsection]{$\!\!$}

\newtheorem{counter*}[subsubsection]{$\!\!$}
\newenvironment{Def*}{\begin{counter*} {\bf Definition.}}{\end{counter*}}
\newenvironment{Not*}{\begin{counter*} \rm {\bf Notation.}}{\end{counter*}}
\newenvironment{Notss*}{\begin{counter*} \rm {\bf Notations.}}{\end{counter*}}
\newenvironment{DefNot*}{\begin{counter*} \rm {\bf Definition-Notation.}}{\end{counter*}}
\newenvironment{Nots*}{\begin{counter*} \rm {\bf Notations.}}{\end{counter*}}
\newenvironment{Prop*}{\begin{counter*} {\bf Proposition.}}{\end{counter*}}
\newenvironment{Lem*}{\begin{counter*} {\bf Lemma.}}{\end{counter*}}
\newenvironment{Cor*}{\begin{counter*} {\bf Corollary.}}{\end{counter*}}
\newenvironment{Th*}{\begin{counter*} {\bf Theorem.}}{\end{counter*}}
\newenvironment{Rem*}{\begin{counter*} \rm {\bf Remark.}}{\end{counter*}}
\newenvironment{Ex*}{\begin{counter*} \rm {\bf Example.}}{\end{counter*}}
\newenvironment{Exs*}{\begin{counter*} \rm {\bf Examples.}}{\end{counter*}}
\newenvironment{Pt*}{\begin{counter*} \rm}{\end{counter*}}
\newenvironment{Q*}{\begin{counter*} \rm {\bf Question.}}{\end{counter*}}

\newcommand{\iso}{\stackrel{\sim}{\longrightarrow}}

\title{\textbf{\huge{Generic and Mod $p$ \\ Kazhdan-Lusztig Theory for $GL_2$}}}

\author{Cédric PEPIN and Tobias SCHMIDT}
\date{September 20, 2021}

\begin{document}

\maketitle

\begin{abstract}
Let $F$ be a non-archimedean local field with residue field $\bbF_q$ and let $\mathbf{G}=GL_{2/F}$. Let $\bfq$ be an indeterminate and let $\cH^{(1)}(\bfq)$ be the generic pro-$p$ Iwahori-Hecke algebra of the $p$-adic group $\mathbf{G}(F)$. Let $V_{\mathbf{\whG}}$ be the Vinberg monoid of the dual group $\mathbf{\whG}$. We establish a generic version for $\cH^{(1)}(\bfq)$ of the Kazhdan-Lusztig-Ginzburg spherical representation, the Bernstein map and the Satake isomorphism. We define the flag variety for the monoid $V_{\mathbf{\whG}}$ and establish the characteristic map in its equivariant $K$-theory. These generic constructions recover the classical ones after the specialization 
$\bfq=q\in\bbC$. At $\bfq=q=0\in \overline{\bbF}_q$, the spherical map provides a dual parametrization of all the irreducible 
$\cH^{(1)}_{\overline{\bbF}_q}(0)$-modules. 

\end{abstract}

\tableofcontents

\section{Introduction}

Let $F$ be a non-archimedean local field with ring of integers $o_F$ and residue field $\bbF_q$. Let $\bfG$ be a connected split reductive group over $F$. Let $\cH_k=(k[I\setminus\bfG(F)/I],\star) $ be the Iwahori-Hecke algebra, i.e. the convolution algebra associated to an Iwahori subgroup $I\subset \bfG(F)$, with coefficients in an algebraically closed field $k$. On the other hand, let $\widehat{\bfG}$ be the Langlands dual group of $\bfG$ over $k$, with maximal torus and Borel subgroup 
$\widehat{\bfT}\subset \widehat{\bfB}$ respectively. Let $W_0$ be the finite Weyl group.

\vskip5pt

When $k=\bbC$, the irreducible $\cH_{\bbC}$-modules appear as subquotients of the Grothendieck group 
$K^{\widehat{\bfG}}( \widehat{\bfG}/ \widehat{\bfB})_{\bbC}$ of
$\widehat{\bfG}$-equivariant coherent sheaves on the dual flag variety $\widehat{\bfG}/ \widehat{\bfB}$. As such they can be parametrized by the isomorphism classes of irreducible tame $\widehat{\bfG}(\bbC)$-representations of the Weil group $\cW_F$ of $F$ with unipotent inertial type, thereby realizing a tame part of the local Langlands correspondence (in this setting also called the Deligne-Lusztig conjecture for Hecke modules): Kazhdan-Lusztig \cite{KL87}, Ginzburg \cite{CG97}. 
Their approach to the Deligne-Lusztig conjecture can be divided into two parts: the first part develops the theory of the so-called {\it spherical representation} leading to a certain dual parametrization of Hecke modules. The second part links these dual data to representations of the group $\cW_F$.

\vskip5pt 

The spherical representation is a distinguished faithful action of the Hecke algebra $\cH_{\bbC}$ on a maximal commutative subring $\cA_{\bbC}\subset\cH_{\bbC}$ via $\cA_{\bbC}^{W_0}$-linear operators: elements of the subring 
$\cA_{\bbC}$ act by multiplication, whereas the standard Hecke operators $T_s\in\cH_{\bbC}$, supported on double cosets indexed by simple reflections $s\in W_0$, act via the classical Demazure operators \cite{D73,D74}. The link with the geometry of the dual group comes then in two steps. First, the classical Bernstein map $\tilde{\theta}$ identifies the ring of functions $\bbC[\widehat{\bfT}]$ with $\cA_{\bbC}$, such that the invariants  $\bbC[\widehat{\bfT}]^{W_0}$ become the center $Z(\cH_{\bbC})=\cA_{\bbC}^{W_0}$. Second, the characteristic homomorphism $c_{\mathbf{\whG}}$ of equivariant $K$-theory identifies the rings $\bbC[\widehat{\bfT}]$ and $K^{\widehat{\bfG}}( \widehat{\bfG}/ \widehat{\bfB})_{\bbC}$ as 
algebras over the representation ring $\bbC[\widehat{\bfT}]^{W_0}=R(\widehat{\bfG})_{\bbC}$.

\vskip5pt

When $k=\overline{\bbF}_q$, any irreducible
$\widehat{\bfG}(\overline{\bbF}_q)$-representation of
$\cW_F$ is tame, with semisimple inertial type. Dually,  one replaces the Iwahori-Hecke algebra by the bigger pro-$p$-Iwahori-Hecke algebra $$\cH_{\overline{\bbF}_q}^{(1)}=(\overline{\bbF}_q[I^{(1)}\setminus \bfG(F)/I^{(1)}],\star),$$ 
where $I^{(1)}\subset I$ is the pro-$p$-radical of $I$.
The algebra $\cH_{\overline{\bbF}_q}^{(1)}$ was introduced by Vign\'eras and its structure theory developed in a series of papers
\cite{V04,V05,V06,V14,V15,V16,V17}. More generally, Vign\'eras introduces and studies a generic version $\cH^{(1)}(\bfq)$ of this algebra
 which is defined over a polynomial ring $\bbZ[\bfq]$ in an indeterminate $\bfq$. The mod $p$ ring $\cH_{\overline{\bbF}_q}^{(1)}$ is obtained by specialization $\bfq=q$ followed by extension of scalars from $\bbZ$ to $\overline{\bbF}_q$, in short $\bfq=q=0$.


\vskip5pt
The present paper is the first in a series of papers in which we will show that there is a {\it generic} version of 
Kazhdan-Lusztig theory, which applies to the generic pro-$p$ Iwahori-Hecke algebra $\cH^{(1)}(\bfq)$. On the one hand, it gives back (and actually improves) the classical theory after passing to the direct summand 
$\cH(\bfq)\subset \cH^{(1)}(\bfq)$ and then specializing $\bfq=q\in\bbC$. On the other hand, it gives a genuine mod $p$ theory after specializing $\bfq=q=0\in \overline{\bbF}_q$. Our key observation is that, in the generic setting, \emph{the Langlands dual group $\widehat{\bfG}$ needs to be enlarged to its Vinberg monoid} $V_{\widehat{\bfG}}$ \cite{V95}.

\vskip5pt

We will work in increasing generality, starting in the present paper with the theory of the spherical representation and the dual parametrization in the simplest case of the group $\bfG=\mathbf{GL_2}$. In this case, all geometry becomes explicit and the representation theory reduces to linear algebra computations. The monoid $V_{\widehat{\bfG}}$ comes with a fibration $\bfq : V_{\widehat{\bfG}}\rightarrow\bbA^1$ and the dual parametrization of $\cH_{\overline{\bbF}_q}^{(1)}$-modules is achieved by working over the $0$-fiber $V_{\widehat{\bfG},0}$.


\vskip5pt

So let $\bfG=\mathbf{GL_2}$ from now on. Let $k=\overline{\bbF}_q$ and $\bfq$ be an indeterminate. 
Let $\bfT\subset\bfG$ be the torus of diagonal matrices.
Let $\cA^{(1)}(\bfq) \subset \cH^{(1)}(\bfq)$ be the maximal commutative subring\footnote{for the choice of the antidominant spherical orientation}
and $\cA^{(1)}(\bfq)^{W_0} = Z(\cH^{(1)}(\bfq))$ be its ring of invariants. We let $\tilde{\bbZ}:=\bbZ[\frac{1}{q-1},\mu_{q-1}]$ and denote by $\tilde{\bullet}$ the base change from $\bbZ$ to $\tilde{\bbZ}$. The algebra $\tilde{\cH}^{(1)}(\bfq)$ splits as a direct product of subalgebras $\tilde{\cH}^{\gamma}(\bfq)$ indexed by the orbits $\gamma$ of $W_0$ in the set of characters of the finite torus $\bbT:=\bfT(\bbF_q)$. There are regular resp. non-regular components 
corresponding to $|\gamma|=2$ resp. $|\gamma|=1$ and the algebra structure of $\tilde{\cH}^{\gamma}(\bfq)$ in these two cases is fundamentally different. 
We define an analogue of the Demazure operator for the regular components and call it
the {\it Vign\'eras operator}. Passing to the product over all $\gamma$, this allows us to single out a distinguished $Z(\tilde{\cH}^{(1)}(\bfq))$-linear operator on 
$\tilde{\cA}^{(1)}(\bfq)$. Our first main result is the existence of the {\it generic pro-$p$ spherical representation}: 

\vskip5pt

{\bf Theorem A.} (cf. \ref{sA2q}, \ref{sA1q}) {\it  There is a (essentially unique) faithful representation 

$$
\xymatrix{
\tilde{\sA}^{(1)}(\bfq):\tilde{\cH}^{(1)}(\bfq) \ar[r] & \End_{Z(\tilde{\cH}^{(1)}(\bfq))}(\tilde{\cA}^{(1)}(\bfq))
}
$$
such that 
\begin{itemize}
\item[(i)]
$
\tilde{\sA}^{(1)}(\bfq)|_{\tilde{\cA}^{(1)}(\bfq)}=\textrm{ the natural inclusion $\tilde{\cA}^{(1)}(\bfq)\subset\End_{Z(\tilde{\cH}^{(1)}(\bfq))}(\tilde{\cA}^{(1)}(\bfq))$}
$
\item[(ii)]
$
\tilde{\sA}^{(1)}(\bfq)(T_s)=\textrm{ the Demazure-Vign\'eras operator on $\tilde{\cA}^{(1)}(\bfq)$}.\;\;\;\;\;\;\;\;\;\;
$
\end{itemize}

Restricting the representation $\tilde{\sA}^{(1)}(\bfq)$ to the Iwahori component, its base change 
$\bbZ[\bfq]\rightarrow \bbZ[\bfq^{\pm \frac{1}{2}} ]$ coincides with the classical spherical representation of Kazhdan-Lusztig and Ginzburg. }

\vskip5pt
We call the left $\tilde{\cH}^{(1)}(\bfq)$-module defined by $\tilde{\sA}^{(1)}(\bfq)$ the \emph{generic spherical module} $\tilde{\cM}^{(1)}$.

\vskip5pt 

Let $\Mat_{2\times 2}$ be the $\bbZ$-monoid scheme of $2\times 2$-matrices. The Vinberg monoid $V_{\widehat{\bfG}}$, as introduced in \cite{V95}, in the particular case of $\mathbf{GL_2}$ is the  
$\bbZ$-monoid scheme
$$
V_{\mathbf{GL_2}}:=\Mat_{2\times 2}\times\bbG_m.
$$ 
It implies the striking interpretation of the formal indeterminate $\bfq$ as a regular function. Indeed,
denote by $z_2$ the canonical coordinate on $\bbG_m$. Let $\bfq$ be the homomorphism from $V_{\mathbf{GL_2}}$ to the multiplicative monoid $(\bbA^1,\cdot )$ defined by $(f,z_2)\mapsto \det(f)z_2^{-1}$:
$$
\xymatrix{
V_{\mathbf{GL_2}} \ar[d]_{\bfq} \\
\bbA^1.
}
$$
The fibration $\bfq$ is trivial over  $\bbA^{1}\setminus \{ 0 \}$ with fibre $\mathbf{GL_2}$. 
The special fiber at $\bfq=0$ is the $\bbZ$-semigroup scheme
$$
V_{\mathbf{GL_2},0} := \bfq^{-1}(0)= \Sing_{2\times 2}\times\bbG_m,
$$
where $\Sing_{2\times 2}$ represents the singular $2\times 2$-matrices. Let $\Diag_{2\times 2}\subset \Mat_{2\times 2}$ be the submonoid scheme of diagonal $2\times 2$-matrices, and set
$$
V_{\mathbf{\whT}}:=\Diag_{2\times 2}\times\bbG_m\subset V_{\mathbf{GL_2}}= \Mat_{2\times 2}\times\bbG_m.
$$
This is a diagonalizable $\bbZ$-monoid scheme. Restricting the above $\bbA^1$-fibration to $V_{\mathbf{\whT}}$ we obtain a fibration, trivial 
over $\bbA^{1}\setminus \{ 0 \}$ with fibre $\mathbf{\whT}$. Its special fibre at $\bfq=0$ is the $\bbZ$-semigroup scheme 
$$
V_{\mathbf{\whT},0} := \bfq|_{V_{\mathbf{\whT}}} ^{-1}(0)= \SingDiag_{2\times2}\times\bbG_m,
$$
where $\SingDiag_{2\times2}$ represents the singular diagonal $2\times 2$-matrices. To ease notion, we denote the base change to $\overline{\bbF}_q$ of these $\bbZ$-schemes by the same symbols. Let $\bbT^{\vee}$ be the finite abelian dual group of $\bbT$. 
We let $R(V^{(1)}_{\mathbf{\whT}})$ be the representation ring of the 
extended monoid 
$$V^{(1)}_{\mathbf{\whT}}:=\bbT^{\vee}\times V_{\mathbf{\whT}}.$$ 
Our second main result is the existence of the {\it generic pro-$p$ Bernstein isomorphism}. 

\vskip5pt 

{\bf Theorem B.} (cf. \ref{genB1}) {\it There exists a ring isomomorphism
$$
\xymatrix{
\sB^{(1)}(\bfq):\cA^{(1)}(\bfq)\ar[r]^>>>>>{\sim} & R(V^{(1)}_{\mathbf{\whT}})
} 
$$
with the property: Restricting the isomorphism $\sB^{(1)}(\bfq)$ to the Iwahori component, its base change 
$\bbZ[\bfq]\rightarrow \bbZ[\bfq^{\pm \frac{1}{2}} ]$ recovers\footnote{By 'recovers' we mean 'coincides up to a renormalization'.} the classical Bernstein isomorphism $\tilde{\theta}$. 

}

\vskip5pt 

The extended monoid $V^{(1)}_{\mathbf{\whT}}$ has a natural $W_0$-action and the isomorphism $\sB^{(1)}(\bfq)$ is equivariant. We call the resulting ring isomorphism 

$$
\xymatrix{
\sS^{(1)}(\bfq):=\sB^{(1)}(\bfq)^{W_0}: \cA^{(1)}(\bfq)^{W_0}\ar[r]^<<<<<{\sim} & R(V^{(1)}_{\mathbf{\whT}})^{W_0}
}
$$
the \emph{generic pro-$p$-Iwahori Satake isomorphism}. Our terminology is justified by the following. Let $K=\bfG(o_F)$. 
Recall that the spherical Hecke algebra of $\mathbf{G}(F)$ 
with coefficients in any commutative ring $R$ is defined to be the convolution algebra
$$
\cH_{R}^{\sph}:=(R[K\backslash \mathbf{G}(F)/K],\star)
$$
generated by the $K$-double cosets in $\mathbf{G}(F)$. We define a {\it generic spherical Hecke algebra} $\cH^{\sph}(\bfq)$ over the ring $\bbZ[\bfq]$. Its base change $\bbZ[\bfq]\rightarrow R$, $\bfq\mapsto q$ coincides with $\cH_{R}^{\sph}$. Our third main result is the existence of the {\it generic Satake isomorphism}.

\vskip5pt 
{\bf Theorem C.} (cf. \ref{ThgenSat}) {\it There exists a ring isomorphism
$$
\xymatrix{
\sS(\bfq):\cH^{\sph}(\bfq)\ar[r]^<<<<<{\sim} & R(V_{\mathbf{\whT}})^{W_0}
}
$$
with the propery: 
Base change 
$\bbZ[\bfq]\rightarrow \bbZ[\bfq^{\pm \frac{1}{2}} ]$ and specialization $\bfq\mapsto q\in\bbC$ 
recovers\footnotemark[\value{footnote}] the classical Satake isomorphism between $\cH^{\sph}_{\bbC}$ and $R(\mathbf{\whT})_{\bbC}^{W_0}$.
}
\vskip5pt 
We emphasize that the possibility of having a generic Satake isomorphism is conceptually new and of independent interest. Its definition relies on the deep Kazhdan-Lusztig theory for the intersection cohomology on the affine flag manifold. Its proof follows from the classical case 
by specialization (to an infinite number of points $q$). The special fibre $\sS(0)$ recovers Herzig's mod $p$ Satake isomorphism \cite{H11}, by choosing Steinberg coordinates on $V_{\mathbf{\whT},0}$.

\vskip5pt 

As a corollary we obtain the \emph{generic central elements morphism} as the unique ring homomorphism
$$
\xymatrix{
\sZ(\bfq):\cH^{\sph}(\bfq)\ar[r] & \cA(\bfq)\subset\cH(\bfq)
}
$$
making the diagram 
$$
\xymatrix{
 \cA(\bfq) \ar[rr]_{\sim}^{\sB^{(1)}(\bfq)|_{\cA(\bfq)}}  && R(V_{\mathbf{\whT}}) \\
\cH^{\sph}(\bfq) \ar[u]^{\sZ(\bfq)} \ar[rr]_{\sim}^{\sS(\bfq)} && R(V_{\mathbf{\whT}})^{W_0} \ar@{^{(}->}[u]
}
$$
commutative. The morphism $\sZ(\bfq)$ is injective and has image $Z(\cH(\bfq))$. Base change 
$\bbZ[\bfq]\rightarrow \bbZ[\bfq^{\pm \frac{1}{2}} ]$ and specialization $\bfq\mapsto q\in\bbC$ recovers\footnotemark[\value{footnote}] Bernstein's classical central elements morphism. Its specialization $\bfq\mapsto q=0\in\overline{\bbF}_q$ coincides with Ollivier's construction from \cite{O14}.

\vskip5pt 

Our fourth main result is the {\it characteristic homomorphism} in the equivariant $K$-theory over the Vinberg monoid $V_{\mathbf{\whG}}$. 
The monoid $V_{\mathbf{\whG}}$ carries an action by multiplication on the right from the $\bbZ$-submonoid scheme

$$
V_{\mathbf{\whB}}:=\UpTriang_{2\times 2}\times\bbG_m\subset \Mat_{2\times 2}\times\bbG_m=V_{\mathbf{\whG}}
$$
where $\UpTriang_{2\times 2}$ represents the upper triangular $2\times 2$-matrices. One can construct (virtual) quotients in the context of semigroups and categories of equivariant vector bundles and their $K$-theory on such quotients, similar to the classical description over a groupoid, and the usual induction functor for vector bundles gives a characteristic homomorphism, which is an isomorphism in the case of monoids \cite{PS20}. Applying this general formalism, the \emph{flag variety} $V_{\mathbf{\whG}}/V_{\mathbf{\whB}}$ resp. its extended version 
$V^{(1)}_{\mathbf{\whG}}/V^{(1)}_{\mathbf{\whB}}$ is defined as a $\bbZ$-monoidoid (instead of a groupoid). 

\vskip5pt 
{\bf Theorem D.} (cf. \ref{cVGL2}) {\it Induction of equivariant vector bundles defines a characteristic isomorphism 

$$
\xymatrix{
c_{V^{(1)}_{\mathbf{\whG}}}: R(V^{(1)}_{\mathbf{\whT}}) \ar[r]^<<<<<{\sim} &  K^{V^{(1)}_{\mathbf{\whG}}}(V^{(1)}_{\mathbf{\whG}}/V^{(1)}_{\mathbf{\whB}}).
}
$$
The ring isomorphism is $R(V^{(1)}_{\mathbf{\whT}})^{W_0}=R(V^{(1)}_\mathbf{\whG})$-linear and compatible with passage to $\bfq$-fibres. 
Over the open complement $\bfq\neq 0$, its Iwahori-component coincides with the classical characteristic homomorphism $c_{\mathbf{\whG}}$ between $R(\mathbf{\whT})$ and $K^{\mathbf{\whG}}(\mathbf{\whG}/\mathbf{\whB})$.
}

\vskip5pt 

We define the {\it category of Bernstein resp. Satake parameters} $\BP_{\mathbf{\whG}}$ resp. $\SP_{\mathbf{\whG}}$ to be the category of quasi-coherent modules on the $\tilde{\bbZ}$-scheme $V^{(1)}_{\mathbf{\whT}}$ resp. 
$V^{(1)}_{\mathbf{\whT}}/W_0$. By Theorem B, restriction of scalars to the subring $\tilde{\cA}^{(1)}(\bfq)$ or $Z(\tilde{\cH}^{(1)}(\bfq))$ defines a functor $B$ resp. $P$ from the 
category of $\tilde{\cH}^{(1)}(\bfq)$-modules to the categories $\BP_{\mathbf{\whG}}$ resp. $\SP_{\mathbf{\whG}}$. 
For example, 
the Bernstein resp. Satake parameter of the spherical module $\tilde{\cM}^{(1)}$ equals the structure sheaf $\cO_{V^{(1)}_{\mathbf{\whT}}}$ resp. the quasi-coherent sheaf corresponding to the $R(V^{(1)}_{\mathbf{\whT}})^{W_0}$-module $K^{V^{(1)}_{\mathbf{\whG}}}(V^{(1)}_{\mathbf{\whG}}/V^{(1)}_{\mathbf{\whB}})$. We call $P$ the \emph{generic parametrization functor}.



\vskip5pt 
 
In the other direction, we define the \emph{generic spherical functor}
 to be the functor $\Sph:=  (\tilde{\cM}^{(1)}\otimes_{Z(\tilde{\cH}^{(1)}(\bfq))}\bullet)\circ S^{-1}$ where $S$ is the Satake equivalence between $Z(\tilde{\cH}^{(1)}(\bfq))$-modules and $\SP_{\mathbf{\whG}}$. Let $\pi: V^{(1)}_{\mathbf{\whT}}\rightarrow V^{(1)}_{\mathbf{\whT}}/W_0$ be the projection. The relation between all these functors is expressed by the commutative diagram: 
$$
\xymatrix{
&\Mod(\tilde{\cH}^{(1)}(\bfq))\ar[d]^{B} \ar[dr]^{P} &  \\
\SP_{\mathbf{\whG}} \ar[ur]^{\Sph} \ar[r]_{\pi^*}  & \BP_{\mathbf{\whG}}\ar[r]_{\pi_*}& \SP_{\mathbf{\whG}}.
}
$$
This ends our discussion of the theory in the generic setting.

\vskip5pt 

Then we pass to the special fibre, i.e. we perform the base change $\bbZ[\bfq]\rightarrow k=\overline{\bbF}_q$, $\bfq\mapsto q=0$.
Identifying the $k$-points of the $k$-scheme $V^{(1)}_{\mathbf{\whT},0}/W_0$ with the skyscraper sheaves on it, the spherical functor $\Sph$ induces a map
$$
\xymatrix{
\Sph:\big(V^{(1)}_{\mathbf{\whT},0}/W_0  \big)(k)\ar[r] & \{\textrm{left $\cH_{\overline{\bbF}_q}^{(1)}$-modules}\}.
}
$$
Considering the decomposition of $V^{(1)}_{\mathbf{\whT},0}/W_0$ into its connected components $V^{\gamma}_{\mathbf{\whT},0}/W_0$ indexed by $\gamma\in \bbT^{\vee}/W_0$, the spherical map decomposes as a disjoint union of maps 

$$
\xymatrix{
\Sph^{\gamma}:\big( V^{\gamma}_{\mathbf{\whT},0}/W_0\big)(k) \ar[r] & \{\textrm{left
$\cH_{\overline{\bbF}_q}^{\gamma}$-modules}\}.& 
}
$$
We come to our last main result, the mod $p$ dual parametrization of {\it all} irreducible $\cH_{\overline{\bbF}_q}^{(1)}$-modules via the spherical map. 
\vskip5pt 

{\bf Theorem E.} (cf. \ref{Sphreg}, \ref{Sphnonreg}) {\it 
\begin{itemize}
\item[(i)]
Let $\gamma\in \bbT^{\vee}/W_0$ regular. The spherical map induces a bijection
$$
\xymatrix{
\Sph^{\gamma}:\big(  V^{\gamma}_{\mathbf{\whT},0}/W_0\big)(k)\ar[r]^>>>>>{\sim} & \{\textrm{simple finite dimensional
 left $\cH^{\gamma}_{\overline{\bbF}_q}$-modules}\}/\sim.
}
$$
The singular locus of the parametrizing $k$-scheme
$$
V^{\gamma}_{\mathbf{\whT},0}/W_0\simeq V_{\mathbf{\whT},0}=\SingDiag_{2\times2}\times\bbG_m
$$
is given by $(0,0)\times\bbG_m\subset V_{\mathbf{\whT},0}$ in the standard coordinates, and its $k$-points correspond to the supersingular Hecke modules through the correspondence $\Sph^{\gamma}$.
\item[(ii)]
Let $\gamma\in \bbT^{\vee}/W_0$ be non-regular. Consider the decomposition 
$$
V^{\gamma}_{\mathbf{\whT},0}/W_0 =V_{\mathbf{\whT},0}/W_0\simeq \bbA^1\times\bbG_m=
D(2)_{\gamma}\cup D(1)_{\gamma}
$$
where $D(1)_{\gamma}$ is the closed subscheme defined by the parabola $z_2=z_1^2$ in the Steinberg coordinates $z_1,z_2$ and $D(2)_{\gamma}$ is the open complement. The spherical map induces bijections
$$
\xymatrix{
\Sph^{\gamma}(2):D(2)_{\gamma}(k)\ar[r]^>>>>>{\sim} & \{\textrm{simple $2$-dimensional
 left $\cH^{\gamma}_{\overline{\bbF}_q}$-modules}\}/\sim 
 }\;\;\;\;\;\;\;\;\;\;\;\;\;\;\;\;\;
 $$
 $$\xymatrix{
\Sph^{\gamma}(1):D(1)_{\gamma}(k)\ar[r]^>>>>>{\sim} & \{\textrm{spherical pairs of characters of $\cH^{\gamma}_{\overline{\bbF}_q}$}\}/\sim.
}
$$

The branch locus of the covering 
$$
V_{\mathbf{\whT},0}\lra V_{\mathbf{\whT},0}/W_0\simeq V^{\gamma}_{\mathbf{\whT},0}/W_0
$$ 
is contained in $D(2)_{\gamma}$, with equation $z_1=0$ in Steinberg coordinates, and its $k$-points correspond to the supersingular Hecke modules through the correspondence $\Sph^{\gamma}(2)$.

\end{itemize}

}



\vskip5pt 

In combination with the computation of the Satake parameter $S(\cM_{\overline{\bbF}_p}^{(1)})$ in Theorem D, we get that this dual parametrization of mod $p$ Hecke modules is realized in the equivariant $K$-theory of the dual Vinberg monoid at $\bfq=0$, whose Iwahori block is a natural specialization at $\bfq=0$ of Kazhdan-Lusztig's parametrization for $\bbC$-coefficients. This realizes the first part of a mod $p$ semisimple Langlands correspondence. We refer to \cite{PS} for the detailed relation between mod $p$ Satake parameters and mod $p$ semisimple Galois representations.

\vskip5pt 

Regarding the strategy of proofs, once the Vinberg monoid is introduced, the generic Satake isomorphism is formulated and the generic spherical mo\-du\-le is constructed, everything else follows from Vign\'eras' structure theory of the generic pro-$p$-Iwahori Hecke algebra and her classification of the irreducible representations.



\vskip5pt

 {\it Notation:} In general, the letter $F$ denotes a locally compact complete non-archimedean field with ring of integers $o_F$.
 Let $\bbF_q$ be its residue field, of characteristic $p$ and cardinality $q$. We denote by $\bfG$ the algebraic group $\mathbf{GL_2}$ over $F$ and by $G:=\bfG(F)$ its group of $F$-rational points. Let $\bfT\subset\bfG$ be the torus of diagonal matrices. Finally, $I\subset G$ denotes the upper triangular standard Iwahori subgroup and $I^{(1)}\subset I$ denotes the unique pro-$p$ Sylow subgroup of $I$.

\section{The pro-$p$-Iwahori-Hecke algebra}

\subsection{The generic pro-$p$-Iwahori Hecke algebra}
\begin{Pt*}\label{extendedWeylgroup}
We denote by $\Phi=\{\pm \alpha\}$ the root system of $(\bfG,\bfT)$. We let
$W_0=\{1,s=s_{\alpha}\}$ and $\Lambda=\bbZ\times\bbZ$
be the \emph{finite Weyl group} of $\bfG$ and the \emph{lattice of cocharacters} of $\bfT$ respectively. If
$\bbT=k^{\times}\times k^{\times}
$
denote the \emph{finite torus} $\bfT(\bbF_q)$, then $W_0$ acts naturally on $\bbT\times\Lambda$. The \emph{extended Weyl group} of $\bfG$ is 
$$
W^{(1)}=\bbT\times\Lambda\rtimes W_0.
$$
It contains the \emph{affine Weyl group} and the \emph{Iwahori-Weyl group}
$$
W_{\aff}=\bbZ(1,-1)\rtimes W_0\subseteq W=\Lambda\rtimes W_0.
$$
The affine Weyl group $W_{\aff}$ is a Coxeter group with set of simple reflexions $S_{\aff}=\{s_0,s\}$, where $s_0=(1,-1)s$. Moreover, setting $u=(1,0)s\in W$ and $\Omega=u^{\bbZ}$, we have $W=W_{\aff}\rtimes \Omega.$
The length function $\ell$ on $W_{\aff}$ can then be inflated to $W$ and $W^{(1)}$. 
\end{Pt*}

\begin{Def*} \label{defgenericprop}
Let $\bfq$ be an indeterminate. The \emph{generic pro-$p$ Iwahori Hecke algebra} is the $\bbZ[\bfq]$-algebra $\cH^{(1)}(\bfq)$ defined by generators
$$
\cH^{(1)}(\bfq):=\bigoplus_{w\in W^{(1)}} \bbZ[\bfq] T_w
$$
and relations:
\begin{itemize}
\item braid relations:
$
T_wT_{w'}=T_{ww'}\quad\textrm{for $w,w'\in W^{(1)}$ if $\ell(w)+\ell(w')=\ell(ww')$}
$ \\
\item quadratic relations:
$
T_{\ts}^2=\bfq + c_{s}T_{\ts}\quad \textrm{if $\ts\in S_{\aff}$},
$
where $c_{s}:=\sum_{t\in (1,-1)(k^{\times})}T_t$.
\end{itemize}
\end{Def*}

\begin{Pt*}\label{distinguishedelements}
The identity element is $1=T_1$. Moreover we set
$$
S:=T_s,\quad U:=T_u\quad\textrm{and}\quad S_0:=T_{s_0}=USU^{-1}.
$$
\end{Pt*}

\begin{Def*} \label{defprop}
Let $R$ be any commutative ring. The \emph{pro-$p$ Iwahori Hecke algebra of $G$ with coefficients in $R$} is defined to be the convolution algebra
$$
\cH_{R}^{(1)}:=(R[I^{(1)}\backslash G/I^{(1)}],\star)
$$
generated by the $I^{(1)}$-double cosets in $G$. 
\end{Def*}

\begin{Th*} \textbf{\emph{(Vignéras, \cite[Thm. 2.2]{V16})}} \label{presprop}
Let $\bbZ[\bfq]\ra R$ be the ring homomorphism mapping $\bfq$ to $q$. Then the $R$-linear map
$$
\xymatrix{
\cH^{(1)}(\bfq)\otimes_{\bbZ[\bfq]}R \ar[r] & \cH_{R}^{(1)}
}
$$
sending $T_w$, $w\in W^{(1)}$, to the characteristic function of the double coset $I^{(1)}\backslash w/I^{(1)}$, is an isomorphism of $R$-algebras.
\end{Th*}

\subsection{Idempotents and component algebras}

\begin{Pt*}\label{finiteT}
Recall the finite torus $\bbT=\bfT(\bbF_q)$. Let us consider its group algebra $\tilde{\bbZ}[\bbT]$ over the ring
$$
\tilde{\bbZ}:=\bbZ[\frac{1}{q-1},\mu_{q-1}].
$$
As $q-1$ is invertible in $\tilde{\bbZ}$, so is $|\bbT|=(q-1)^2$. We denote by $\bbT^{\vee}$ the set of characters 
$
\lambda: \bbT\rightarrow \mu_{q-1}\subset\widetilde{\bbZ},
$
with its natural $W_0$-action given by $^s\lambda(t_1,t_2)=\lambda(t_2,t_1)$ for
$(t_1,t_2)\in\bbT$. The set of $W_0$-orbits $\bbT^{\vee}/W_0$ has cardinality $\frac{q^2-q}{2}$. Also $W^{(1)}$ acts on
$\bbT^{\vee}$ through the canonical quotient map $W^{(1)}\rightarrow W_0$. Because of the braid relations in $\cH^{(1)}(\bfq)$, the rule $t\mapsto T_t$ induces an embedding of $\tilde{\bbZ}$-\emph{algebras}
$$
\tilde{\bbZ}[\bbT]\subset \cH_{\tilde{\bbZ}}^{(1)}(\bfq):=\cH^{(1)}(\bfq)\otimes_{\bbZ}\tilde{\bbZ}.
$$ 

\end{Pt*}

\begin{Def*}
For all $\lambda\in \bbT^{\vee}$ and for $\gamma\in \bbT^{\vee}/W_0$, we define
$$
\varepsilon_{\lambda}:=|\bbT|^{-1}\sum_{t\in\bbT}\lambda^{-1}(t)T_t\hskip10pt \text{and} \hskip10pt \varepsilon_{\gamma}:=\sum_{\lambda\in\gamma}\varepsilon_{\lambda}.
$$

\end{Def*}

\begin{Lem*} \label{decompprop}
The elements $\varepsilon_{\lambda}$, $\lambda\in \bbT^{\vee}$, are idempotent, pairwise orthogonal and their sum is equal to $1$. The elements $\varepsilon_{\gamma}$, $\gamma\in \bbT^{\vee}/W_0$, are idempotent, pairwise orthogonal, their sum is equal to $1$ and they are central in $\cH^{(1)}_{\widetilde{\bbZ}}(\bfq)$. The $\tilde{\bbZ}[\bfq]$-algebra $\cH^{(1)}_{\tilde{\bbZ}}(\bfq)$ is the direct product of its 
sub-$\tilde{\bbZ}[\bfq]$-algebras $\cH^{\gamma}_{\tilde{\bbZ}}(\bfq):=\cH^{(1)}_{\tilde{\bbZ}}(\bfq)\varepsilon_{\gamma}$:
$$
\cH^{(1)}_{\tilde{\bbZ}}(\bfq)=\prod_{\gamma\in\bbT^{\vee}/W_0}\cH^{\gamma}_{\tilde{\bbZ}}(\bfq).
$$

In particular, the category of $\cH_{\tilde{\bbZ}}^{(1)}(\bfq)$-modules decomposes into a finite product of the module categories for the individual component rings $\cH^{(1)}_{\tilde{\bbZ}}(\bfq)\varepsilon_{\gamma}$.
\end{Lem*}

\begin{proof}
The elements $\varepsilon_{\gamma}$ are central because of the relations $T_sT_t=T_{s(t)}T_s$, $T_{s_0}T_t=T_{s_0(t)}T_{s_0}$ and $T_uT_t=T_{s(t)}T_u$ for all $t\in(1,-1)k^{\times}$.
\end{proof}

\begin{Pt*}
Following the terminology of \cite{V04}, we call $|\gamma|=2$ a {\it regular} case and $|\gamma|=1$ a {\it non-regular} (or {\it Iwahori}) case.
\end{Pt*}

\subsection{The Bernstein presentation} \label{The Bernstein presentation}

The inverse image in $W^{(1)}$ of any subset of $W$ along the canonical projection $W^{(1)}\ra W$ will be denoted with a superscript 
${}^{(1)}$.

\begin{Th*} \textbf{\emph{(Vignéras \cite[Th. 2.10, Cor 5.47]{V16})}} \label{Bernsteinpresprop}
The $\bbZ[\bfq]$-algebra $\cH^{(1)}(\bfq)$ admits the following \emph{Bernstein presentation}:
$$
\cH^{(1)}(\bfq)=\bigoplus_{w\in W^{(1)}}\bbZ[\bfq] E(w)
$$
satisfying
\begin{itemize}
\item braid relations:
$
E(w)E(w')=E(ww')\quad\textrm{for $w,w'\in W_0^{(1)}$ if $\ell(w)+\ell(w')=\ell(ww')$}
$
\item quadratic relations:
$
E(\ts)^2=\bfq E(\ts^2)+c_{\ts}E(\ts) \textrm{ if $\ts=ts \in s^{(1)}$},
$
where $c_{\ts}:=T_{s(t)}c_s$ with $t\in\bbT$  
\item product formula:
$
E(\lambda)E(w)=\bfq^{\frac{\ell(\lambda)+\ell(w)-\ell(\lambda w)}{2}}E(\lambda w)\quad \textrm{for $\lambda\in\Lambda^{(1)}$ and $w\in W^{(1)}$}
$
\item Bernstein relations for $\ts\in s^{(1)}$ and $\lambda\in\Lambda^{(1)}$ : set $V:=\bbR\Phi^{\vee}$ and let
$$
\nu:\Lambda^{(1)}\ra V
$$
be the homomorphism such that $\lambda\in\Lambda^{(1)}$ acts on $V$ by translation by $\nu(\lambda)$ ; then the Bernstein element
$$
B(\lambda,\ts):=E(\ts \lambda \ts^{-1})E(\ts)-E(\ts)E(\lambda)
$$
\begin{eqnarray*}
=& 0 & \textrm{if $\lambda\in(\Lambda^s)^{(1)}$} \\
= &\sign(\alpha\circ\nu(\lambda))\sum_{k=0}^{|\alpha\circ\nu(\lambda)|-1}\mathbf{q}(k,\lambda)c(k,\lambda)E(\mu(k,\lambda)) & \textrm{if $\lambda\in\Lambda^{(1)}\setminus(\Lambda^s)^{(1)}$}
\end{eqnarray*}
where $\bfq(k,\lambda)c(k,\lambda)\in\bbZ[\bfq][\bbT]$ and $\mu(k,\lambda)\in\Lambda^{(1)}$ are explicit, cf. \cite[Th. 5.46]{V16} and references therein.
\end{itemize}
\end{Th*}

\begin{Pt*}\label{AIH}
Let 
$$
\cA(\bfq):=\bigoplus_{\lambda\in\Lambda}\bbZ[\bfq]E(\lambda)\subset \cA^{(1)}(\bfq):=\bigoplus_{\lambda\in\Lambda^{(1)}}\bbZ[\bfq]E(\lambda)\subset\cH^{(1)}(\bfq).
$$
It follows from the product formula that these are \emph{commutative sub-$\bbZ[\bfq]$-algebras of $\cH^{(1)}(\bfq)$}. Moreover, by definition \cite[5.22-5.25]{V16}, we have $E(t)=T_t$ for all $t\in\bbT$, so that $\bbZ[\bbT]\subset \cA^{(1)}(\bfq)$. Then, again by the product formula, the commutative algebra $\cA^{(1)}(\bfq)$ decomposes as the 
tensor product of the subalgebras 
$$
\cA^{(1)}(\bfq) =\bbZ[\bbT]\otimes_{\bbZ} \cA(\bfq).
$$
Also, after base extension $\bbZ\ra\tilde{\bbZ}$, we can set $\cA_{\tilde{\bbZ}}^{\gamma}(\bfq):=\cA_{\tilde{\bbZ}}^{(1)}(\bfq)\varepsilon_{\gamma}$, and obtain the decomposition
$$
\cA_{\tilde{\bbZ}}^{(1)}(\bfq)=\prod_{\gamma\in\bbT^{\vee}/W_0}\cA^{\gamma}_{\tilde{\bbZ}}(\bfq)\subset\prod_{\gamma\in\bbT^{\vee}/W_0}\cH^{\gamma}_{\tilde{\bbZ}}(\bfq)=\cH^{(1)}_{\tilde{\bbZ}}(\bfq).
$$
\end{Pt*}

\begin{Lem*}\label{presAq} Let $X$,Y$,z_2$ be indeterminates. There exists a unique ring homomorphism 
$$
\xymatrix{
\bbZ[\bfq][z_2^{\pm1}][X,Y]/(XY-\bfq z_2)\ar[r] & \cA(\bfq)
}
$$
such that
$$
X\lmapsto E(1,0),\quad Y\lmapsto E(0,1)\quad\textrm{and}\quad z_2\lmapsto E(1,1). 
$$
It is an isomorphism. Moreover, for all $\gamma\in\bbT^{\vee}/W_0$,
$$
\cA_{\tilde{\bbZ}}^{\gamma}(\bfq)=
\left\{ \begin{array}{ll}
(\tilde{\bbZ}\varepsilon_{\lambda}\times\tilde{\bbZ}\varepsilon_{\mu})\otimes_{\bbZ}\cA(\bfq) & \textrm{ if $\gamma=\{\lambda,\mu\}$ is regular} \\
\tilde{\bbZ}\varepsilon_{\lambda}\otimes_{\bbZ}\cA(\bfq)& \textrm{ if $\gamma=\{\lambda\}$ is non-regular}.
\end{array} \right.
$$
\end{Lem*}

\begin{proof}
For any $(n_1,n_2)\in\bbZ^2=\Lambda$, we have $\ell(n_1,n_2)=|n_1-n_2|$. Hence it follows from product formula that $z_2$ is invertible and $XY=\bfq z_2$, so that we get a $\bbZ[\bfq]$-algebra homomorphism
$$
\xymatrix{
\bbZ[\bfq][z_2^{\pm1}][X,Y]/(XY-\bfq z_2)\ar[r] & \cA(\bfq).
}
$$ 
Moreover it maps the $\bbZ[\bfq][z_2^{\pm1}]$-basis $\{X^n\}_{n>1}\coprod\{1\}\coprod \{Y^n\}_{n>1}$ to the $\bbZ[\bfq][z_2^{\pm1}]$-basis 
$$
\{E(n,0)\}_{n>1}\coprod\{1\}\coprod \{E(0,n)\}_{n>1},
$$ 
and hence is an isomorphism. The rest of the lemma is clear since $\cA_{\tilde{\bbZ}}^{(1)}(\bfq)=\tilde{\bbZ}[\bbT]\otimes_{\bbZ}\cA(\bfq)$ and $\tilde{\bbZ}[\bbT]=\prod_{\lambda\in\bbT^{\vee}}\tilde{\bbZ}\varepsilon_{\lambda}$.
\end{proof}
In the following, we will sometimes view the isomorphism of the lemma as an identification and write $X=E(1,0), Y=E(0,1)$ and $z_2=E(1,1)$.

\begin{Pt*} \label{centerHI1}
The rule $E(\lambda)\mapsto E(w(\lambda))$ defines an action of the finite Weyl group $W_0=\{1,s\}$ on $\cA^{(1)}(\bfq)$ by 
$\bbZ[\bfq]$-algebra homomorphisms.
By \cite[Th. 4]{V05} (see also \cite[Th. 1.3]{V14}), the subring of $W_0$-invariants is equal to the center of $\cH^{(1)}(\bfq)$, and the same is true after the scalar extension $\bbZ\ra\tilde{\bbZ}$. Now the action on $\cA_{\tilde{\bbZ}}^{(1)}(\bfq)$ stabilizes each component 
$\cA^{\gamma}_{\tilde{\bbZ}}(\bfq)$ and then the resulting subring of $W_0$-invariants is the center of  $\cH_{\tilde{\bbZ}}^{\gamma}(\bfq)$. In terms of the description of $\cA^{\gamma}_{\tilde{\bbZ}}(\bfq)$ given in Lemma \ref{presAq}, this translates into :
\end{Pt*}

\begin{Lem*} \label{centergamma}
Let $\gamma\in\bbT^{\vee}/W_0$. 
\begin{itemize}
\item If $\gamma=\{\lambda,\mu\}$ is regular, then the map
\begin{eqnarray*}
\cA_{\tilde{\bbZ}}(\bfq) & \lra & \cA_{\tilde{\bbZ}}^{\gamma}(\bfq)^{W_0}=Z(\cH_{\tilde{\bbZ}}^{\gamma}(\bfq)) \\
a & \lmapsto & a\varepsilon_{\lambda}+s(a)\varepsilon_{\mu}
\end{eqnarray*}
is an isomorphism of $\tilde{\bbZ}[\bfq]$-algebras. It depends on the choice of order $(\lambda,\mu)$ on the set $\gamma$.
\item If $\gamma=\{\lambda\}$ is non-regular, then
$$
Z(\cH_{\tilde{\bbZ}}^{\gamma}(\bfq))=\cA_{\tilde{\bbZ}}^{\gamma}(\bfq)^{W_0}=\tilde{\bbZ}[\bfq][z_2^{\pm1},z_1]\varepsilon_{\lambda}
$$
with $z_1:=X+Y$.
\end{itemize}
\end{Lem*}

\begin{Pt*} \label{BernsteinVSIM} One can express $X,Y,z_2\in\cA^{(1)}(\bfq)\subset\cH^{(1)}(\bfq)$ in terms of the distinguished elements \ref{distinguishedelements}. This is an application of \cite[Ex. 5.30]{V16}. We find:
$$
(1,0)=s_0u=us\in\Lambda\Rightarrow X:=E(1,0)=(S_0-c_{s_0})U=U(S-c_s),
$$ 
$$
(0,1)=su\in\Lambda \Rightarrow Y:=E(0,1)=SU,
$$ 
$$
(1,1)=u^2 \in\Lambda\Rightarrow z_2:=E(1,1)=U^2.
$$
Also
$$
z_1:=X+Y=U(S-c_s)+SU.
$$
\end{Pt*}

\section{The generic regular spherical representation} \label{regcomponents}

\subsection{The generic regular Iwahori-Hecke algebras} \label{genericgammalagebrareg}

Let $\gamma=\{\lambda,\mu\}\in\bbT^{\vee}/W_0$ be a regular orbit. 
We define a model $\cH_2(\bfq)$ over $\bbZ$ for the component algebra $\cH_{\tilde{\bbZ}}^{\gamma}(\bfq)\subset \cH_{\tilde{\bbZ}}^{(1)}(\bfq)$.
The algebra $\cH_2(\bfq)$ itself will not depend on $\gamma$.

\begin{Pt*}
By construction, the $\tilde{\bbZ}[\bfq]$-algebra $\cH_{\tilde{\bbZ}}^{\gamma}(\bfq)$ admits the following presentation:

$$
\cH_{\tilde{\bbZ}}^{\gamma}(\bfq) =(\tilde{\bbZ}\varepsilon_{\lambda}\times\tilde{\bbZ}\varepsilon_{\mu})\otimes_{\bbZ}'\bigoplus_{w\in W} \bbZ[\bfq] T_w,
$$
where $\otimes_{\bbZ}'$ is the tensor product $\otimes_{\bbZ}$ of $\bbZ$-modules, whose algebra structure is
\emph{twisted} by the $W$-action on $\{\lambda,\mu\}$ through the quotient map $W\ra W_0$, together with the orthogonality relation $\varepsilon_{\lambda}\varepsilon_{\mu}=0$ and the
\begin{itemize}
\item braid relations:
$
T_wT_{w'}=T_{ww'}\quad\textrm{for $w,w'\in W$ if $\ell(w)+\ell(w')=\ell(ww')$}
$
\item quadratic relations:
$
T_{\ts}^2=\bfq \quad \textrm{if $\ts\in S_{\aff}$}.
$
\end{itemize}
\end{Pt*}

\begin{Def*}
Let $\bfq$ be an indeterminate. The \emph{generic second Iwahori-Hecke algebra} is the $\bbZ[\bfq]$-algebra $\cH_2(\bfq)$ defined by generators
$$
\cH_2(\bfq):=(\bbZ\varepsilon_1\times\bbZ\varepsilon_2)\otimes_{\bbZ}'\bigoplus_{w\in W} \bbZ[\bfq] T_w,
$$
where $\otimes_{\bbZ}'$ is the tensor product $\otimes_{\bbZ}$ of $\bbZ$-modules, whose algebra structure is twisted by the $W$-action on $\{1,2\}$ through the quotient map $W\ra W_0=\mathfrak{S}_2$, together with 
$\varepsilon_{1}\varepsilon_{2}=0$, and the relations:
\begin{itemize}
\item braid relations:
$
T_wT_{w'}=T_{ww'}\quad\textrm{for $w,w'\in W$ if $\ell(w)+\ell(w')=\ell(ww')$}
$
\item quadratic relations:
$
T_{\ts}^2=\bfq \quad \textrm{if $\ts\in S_{\aff}$}.
$
\end{itemize}
\end{Def*}

\begin{Pt*}\label{presH2q}
The identity element of $\cH_2(\bfq)$ is $1=T_1$. Moreover we set in $\cH_2(\bfq)$
$$
S:=T_s,\quad U:=T_u\quad\textrm{and}\quad S_0:=T_{s_0}=USU^{-1}.
$$
Then one checks that
$$
\cH_2(\bfq)=(\bbZ\varepsilon_1\times\bbZ\varepsilon_2)\otimes_{\bbZ}'\bbZ[\bfq][S,U^{\pm 1}],\quad S^2=\bfq,\quad U^2S=SU^2
$$
is a presentation of $\cH_2(\bfq)$. Note that the element $U^2$ is invertible in $\cH_2(\bfq)$.

\end{Pt*}

\begin{Pt*}\label{H2VSHgamma} 
Choosing the ordering $(\lambda,\mu)$ on the set $\gamma=\{\lambda,\mu\}$ and mapping $\varepsilon_1\mapsto\varepsilon_{\lambda}, \varepsilon_2\mapsto \varepsilon_{\mu}$ defines an isomorphism of $\tilde{\bbZ}[\bfq]$-algebras
$$
\xymatrix{
\cH_2(\bfq)\otimes_{\bbZ}\tilde{\bbZ}\ar[r]^<<<<<{\sim} & \cH_{\tilde{\bbZ}}^{\gamma}(\bfq),
}
$$
such that $S\otimes 1\mapsto S\varepsilon_{\gamma}$, $U\otimes 1\mapsto U\varepsilon_{\gamma}$ and $S_0\otimes 1\mapsto S_0\varepsilon_{\gamma}$.
\end{Pt*}

\begin{Pt*} \label{presA2q} We identify two important commutative subrings of $\cH_2(\bfq)$.
We define $\cA_2(\bfq)\subset\cH_2(\bfq)$ to be the $\bbZ[\bfq]$-subalgebra generated by the elements $\varepsilon_1$, $\varepsilon_2$, $US$, $SU$ and $U^{\pm 2}$. Let $X,Y$ and $z_2$ be indeterminates. 
Then there is a unique $(\bbZ\varepsilon_1\times\bbZ\varepsilon_2)\otimes_{\bbZ}\bbZ[\bfq]$-algebra homomorphism 
$$
(\bbZ\varepsilon_1\times\bbZ\varepsilon_2)\otimes_{\bbZ}\bbZ[\bfq][z_2^{\pm1}][X,Y]/(XY-\bfq z_2) \lra \cA_2(\bfq)$$
such that $X\mapsto US, Y\mapsto SU, z_2\mapsto U^2$, and it is an isomorphism. In particular,  $\cA_2(\bfq)$ is a {\it commutative} subalgebra of $\cH_2(\bfq)$.
The isomorphism \ref{H2VSHgamma} identifies $\cA_2(\bfq)\otimes_{\bbZ}\tilde{\bbZ}$ with $\cA_{\tilde{\bbZ}}^{\gamma}(\bfq)$.
Moreover, permuting $\varepsilon_1$ and $\varepsilon_2$, and $X$ and $Y$, extends to an action of $W_0=\mathfrak{S}_2$ on 
$\cA_2(\bfq)$ by homomorphisms of $\bbZ[\bfq]$-algebras, whose invariants is the center $Z(\cH_2(\bfq))$ of $\cH_2(\bfq)$, and the map
\begin{eqnarray*}
\bbZ[\bfq][z_2^{\pm1}][X,Y]/(XY-\bfq z_2) & \lra & \cA_2(\bfq)^{W_0}=Z(\cH_2(\bfq)) \\
a & \lmapsto & a\varepsilon_1+s(a)\varepsilon_2
\end{eqnarray*}
is an isomorphism of $\bbZ[\bfq]$-algebras. This is a consequence of \ref{H2VSHgamma}, \ref{BernsteinVSIM}, \ref{presAq} and \ref{centergamma}.
 In the following, we will sometimes view the above isomorphisms as identifications. In particular, we will write $X=US, Y=SU$ and $z_2=U^2$.
\end{Pt*}

\subsection{The Vignéras operator}

In this subsection and the following, we will investigate the structure of the $Z(\cH_2(\bfq))$-algebra $\End_{Z(\cH_2(\bfq))}(\cA_2(\bfq))$
of $Z(\cH_2(\bfq))$-linear endomorphisms of 
$\cA_2(\bfq)$. Recall from the preceding subsection that $Z(\cH_2(\bfq))= \cA_2(\bfq)^s$ is the subring of invariants of the commutative ring 
$\cA_2(\bfq)$.

\begin{Lem*} 
We have
$$
\cA_2(\bfq)=\cA_2(\bfq)^s\varepsilon_1\oplus\cA_2(\bfq)^s\varepsilon_2
$$
as $\cA_2(\bfq)^{s}$-modules.
\end{Lem*}

\begin{proof}
This is immediate from the two isomorphisms in \ref{presA2q}.
\end{proof}
According to the lemma, we may use the $\cA_2(\bfq)^{s}$-basis $\varepsilon_1,\varepsilon_2$ to identify $\End_{Z(\cH_2(\bfq))}(\cA_2(\bfq))$ with the algebra of $2\times2$-matrices
over $\cA_2(\bfq)^{s}=\bbZ[\bfq][z_2^{\pm1}][X,Y]/(XY-\bfq z_2)$.

\begin{Def*} The endomorphism of $\cA_2(\bfq)$ corresponding to the matrix 
$$
V_s(\bfq):=\left (\begin{array}{cc}
0 & Y\varepsilon_1+X\varepsilon_2 \\
z_2^{-1}(X\varepsilon_1+Y\varepsilon_2) & 0
\end{array} \right)
$$
will be called \emph{the Vign\'eras operator on $\cA_2(\bfq)$}.
\end{Def*}

\begin{Lem*}\label{relationVsq}
We have 
$
V_s(\bfq)^2=\bfq.
$
\end{Lem*}
\begin{proof} This is a short calculation. 
\end{proof}

\subsection{The generic regular spherical representation} \label{HnilqGL2}

In the following theorem we define the generic regular spherical representation of the algebra $\cH_2(\bfq)$ on the $Z(\cH_2(\bfq))$-module
$\cA_2(\bfq)$. Note that the commutative ring
 $\cA_2(\bfq)$ is naturally a subring 
 $$\cA_2(\bfq)\subset \End_{Z(\cH_2(\bfq))}(\cA_2(\bfq)),
 $$
an element $a\in \cA_2(\bfq)$ acting by multiplication $b\mapsto ab$ on $ \cA_2(\bfq)$.

\begin{Th*} \label{sA2q}
There exists a unique $\bbZ[\bfq]$-algebra homomorphism
$$
\xymatrix{
\sA_2(\bfq):\cH_2(\bfq) \ar[r] & \End_{Z(\cH_2(\bfq))}(\cA_2(\bfq))
}
$$
such that 
\begin{itemize}
\item[(i)]
$
\sA_2(\bfq)|_{\cA_2(\bfq)}=\textrm{ the natural inclusion $\cA_2(\bfq)\subset\End_{Z(\cH_2(\bfq))}(\cA_2(\bfq))$}
$
\item[(ii)]
$
\sA_2(\bfq)(S)=V_s(\bfq).
$
\end{itemize}
\end{Th*}

\begin{proof} 
Recall that $\cH_2(\bfq)=(\bbZ\varepsilon_1\times\bbZ\varepsilon_2)\otimes_{\bbZ}'\bbZ[\bfq][S,U^{\pm1}]$ with the relations $S^2=\bfq$ and $U^2S=SU^2$. In particular $\sA_2(\bfq)(S):=V_s(\bfq)$ is well-defined thanks to \ref{relationVsq}. Now let us consider the question of finding the restriction of $\sA_2(\bfq)$ to the subalgebra $\bbZ[\bfq][S,U^{\pm1}]$. As the $\bbZ[\bfq]$-algebra $\cA_2(\bfq)\cap\bbZ[\bfq][S,U^{\pm1}]$ is generated by
$$
z_2=U^2,\quad X=US\quad\textrm{and}\quad Y=SU,
$$
such a $\bbZ[\bfq]$-algebra homomorphism exists if and only if there exists
$$
\sA_2(\bfq)(U)\in\End_{Z(\cH_2(\bfq))}(\cA_2(\bfq))
$$
satisfying
\begin{enumerate}
\item $\sA_2(\bfq)(U)^2=\sA_2(\bfq)(U^2)=\sA_2(\bfq)(z_2)=z_2\Id$ \hskip5pt (in particular $\sA_2(\bfq)(U)$ is invertible) 
\item $\sA_2(\bfq)(U)V_s(\bfq)=\textrm{ multiplication by $X$}$
\item $V_s(\bfq)\sA_2(\bfq)(U)=\textrm{ multiplication by $Y$}$.
\end{enumerate}
As before we use the $Z(\cH_2(\bfq))$-basis $\varepsilon_1,\varepsilon_2$ of $\cA_2(\bfq)$ to identify $\End_{Z(\cH_2(\bfq))}(\cA_2(\bfq))$ with the algebra of $2\times 2$-matrices over the ring $Z(\cH_2(\bfq))=\cA_2(\bfq)^s$. Then, by definition,
$$
V_s(\bfq)=
\left (\begin{array}{cc}
0 & Y\varepsilon_1+X\varepsilon_2 \\
z_2^{-1}(X\varepsilon_1+Y\varepsilon_2) & 0
\end{array} \right).
$$
Moreover, the multiplications by $X$ and by $Y$ on $\cA_2(\bfq)$ correspond then to the matrices
$$
\left (\begin{array}{cc}
X\varepsilon_1+Y\varepsilon_2 & 0 \\
0 & Y\varepsilon_1+X\varepsilon_2
\end{array} \right)
\quad\textrm{and}\quad
\left (\begin{array}{cc}
Y\varepsilon_1+X\varepsilon_2 & 0 \\
0 & X\varepsilon_1+Y\varepsilon_2
\end{array} \right).
$$
Now, writing
$$
\sA_2(\bfq)(U)=
\left (\begin{array}{cc}
a& c \\
b& d
\end{array} \right)
$$
we have:
$$
\sA_2(\bfq)(U)^2=z_2\Id \Longleftrightarrow
\left (\begin{array}{cc}
a^2+bc& c(a+d) \\
b(a+d)& d^2+bc
\end{array} \right)
=
\left (\begin{array}{cc}
z_2& 0 \\
0& z_2
\end{array} \right),
$$
$$
\sA_2(\bfq)(U)V_s(\bfq)=\textrm{ multiplication by $X$}
$$
$$
\Longleftrightarrow
\left (\begin{array}{cc}
cz_2^{-1}(X\varepsilon_1+Y\varepsilon_2)& a(Y\varepsilon_1+X\varepsilon_2) \\
dz_2^{-1}(X\varepsilon_1+Y\varepsilon_2)& b(Y\varepsilon_1+X\varepsilon_2) 
\end{array} \right)
=
\left (\begin{array}{cc}
X\varepsilon_1+Y\varepsilon_2 & 0 \\
0 & Y\varepsilon_1+X\varepsilon_2
\end{array} \right)
$$
and
$$
V_s(\bfq)\sA_2(\bfq)(U)=\textrm{ multiplication by $Y$}
$$
$$
\Longleftrightarrow
\left (\begin{array}{cc}
b(Y\varepsilon_1+X\varepsilon_2)& d(Y\varepsilon_1+X\varepsilon_2) \\
az_2^{-1}(X\varepsilon_1+Y\varepsilon_2)& cz_2^{-1}(X\varepsilon_1+Y\varepsilon_2)
\end{array} \right)
=
\left (\begin{array}{cc}
Y\varepsilon_1+X\varepsilon_2 & 0 \\
0 & X\varepsilon_1+Y\varepsilon_2
\end{array} \right).
$$
Each of the two last systems admits a unique solution, namely
$$
\sA_2(\bfq)(U)
=
\left (\begin{array}{cc}
a& c \\
b& d
\end{array} \right)
=
\left (\begin{array}{cc}
0& z_2 \\
1& 0
\end{array} \right),
$$
which is also a solution of the first one. Moreover, the determinant
$$
ad-bc=-z_2
$$
is invertible. 

Finally, $\cA_2(\bfq)$ is generated by $\cA_2(\bfq)\cap\bbZ[\bfq][S,U^{\pm1}]$ together with $\varepsilon_1$ and $\varepsilon_2$. The latter are assigned to map to the projectors
$$
\textrm{multiplication by $\varepsilon_1$}=\left (\begin{array}{cc}
1 & 0 \\
0 & 0
\end{array} \right)
\quad\textrm{and}\quad
\textrm{multiplication by $\varepsilon_2$}=\left (\begin{array}{cc}
0 & 0 \\
0 & 1
\end{array} \right).
$$
Thus it only remains to check that
$$
\left (\begin{array}{cc}
1 & 0 \\
0 & 0
\end{array} \right)
\sA_2(\bfq)(S)
=
\sA_2(\bfq)(S)
\left (\begin{array}{cc}
0 & 0 \\
0 & 1
\end{array} \right)
$$
and
$$
\left (\begin{array}{cc}
0 & 0 \\
0 & 1
\end{array} \right)
\sA_2(\bfq)(S)
=
\sA_2(\bfq)(S)
\left (\begin{array}{cc}
1 & 0 \\
0 & 0
\end{array} \right),
$$
and similarly with $\sA_2(\bfq)(U)$ in place of $\sA_2(\bfq)(U)$, which is straightforward.
\end{proof}

\begin{Rem*}
The map $\sA_2(\bfq)$, together with the fact that it is an isomorphism (see below), is a rewriting of a theorem of Vignéras, namely \cite[Cor. 2.3]{V04}. In loc. cit., the algebra $\cH_2(\bfq)$ is identified with the algebra of $2\times 2$-matrices over the ring $\bbZ[\bfq][z_2^{\pm1}][X,Y]/(XY-\bfq z_2)$. In our approach, we have replaced the \emph{abstract rank $2$ module} underlying the standard representation of this matrix algebra, by the \emph{subring $\cA_2(\bfq)$ of $\cH_2(\bfq)$} with $\{\varepsilon_1,\varepsilon_2\}$ for the canonical basis. 
\end{Rem*}

\begin{Prop*}\label{sA2sf}
The homomorphism $\sA_2(\bfq)$ is an isomorphism.
\end{Prop*}

\begin{proof}
It follows from \ref{presH2q} and \ref{presA2q} that the $\bbZ[\bfq]$-algebra $\cH_2(\bfq)$ is generated by the elements
$$
\varepsilon_1,\  \varepsilon_2,\  S,\  U,\  SU
$$
as a module over its center $Z(\cH_2(\bfq))$. Moreover, as $SU^2=U^2S=:z_2S$ and $SU=:Y$, we have
$$
S=z_2^{-1}YU=z_2^{-1}Y(\varepsilon_1U+\varepsilon_2U)=z_2^{-1}(Y\varepsilon_1+X\varepsilon_2)\varepsilon_1U+z_2^{-1}(X\varepsilon_1+Y\varepsilon_2)\varepsilon_2U,
$$
$$
 U=\varepsilon_1U+\varepsilon_2U\quad\textrm{and}\quad SU=(Y\varepsilon_1+X\varepsilon_2)\varepsilon_1+(X\varepsilon_1+Y\varepsilon_2)\varepsilon_2.
$$
Consequently $\cH_2(\bfq)$ is generated as a $Z(\cH_2(\bfq))$-module by the elements 
$$
\varepsilon_1,\  \varepsilon_2,\  z_2^{-1}\varepsilon_1U,\  \varepsilon_2U.
$$
Since 
$$
\sA_2(\bfq)(U)
:=
\left (\begin{array}{cc}
0& z_2 \\
1& 0
\end{array} \right),
$$
these four elements are mapped by $\sA_2(\bfq)$ to
$$
\left (\begin{array}{cc}
1& 0 \\
0& 0
\end{array} \right),\
\left (\begin{array}{cc}
0& 0 \\
0& 1
\end{array} \right),\
\left (\begin{array}{cc}
0& 1 \\
0& 0
\end{array} \right),\
\left (\begin{array}{cc}
0& 0 \\
1& 0
\end{array} \right).
$$
As $\sA_2(\bfq)$ indentifies $Z(\cH_2(\bfq))\subset\cH_2(\bfq)$ with the center of the matrix algebra 
$$
\End_{Z(\cH_2(\bfq))}(\cA_2(\bfq))=\End_{Z(\cH_2(\bfq))}(Z(\cH_2(\bfq))\varepsilon_1\oplus Z(\cH_2(\bfq))\varepsilon_2),
$$
it follows that the elements $\varepsilon_1$, $\varepsilon_2$,  $z_2^{-1}\varepsilon_1U$, $\varepsilon_2U$ are linearly independent over $Z(\cH_2(\bfq))$ and that  $\sA_2(\bfq)$ is an isomorphism.
\end{proof}

We record the following corollary of the proof.

\begin{Cor*}\label{dimOverCenterH2} 
The ring $\cH_2(\bfq)$ is a free $Z(\cH_2(\bfq))$-module on the basis 
$
\varepsilon_1,\varepsilon_2,z_2^{-1}\varepsilon_1U,\varepsilon_2U.
$
\end{Cor*}

\begin{Pt*}
We end this section by noting an equivariance property of $\sA_2(\bfq)$. As already noticed, the finite Weyl group $W_0$ acts on $\cA_2(\bfq)$ by $\bbZ[\bfq]$-algebra automorphisms, and the action is clearly faithful. Moreover $\cA_2(\bfq)^{W_0}=Z(\cH_2(q))$. Hence $W_0$ can be viewed as a subgroup of $\End_{Z(\cH_2(\bfq))}(\cA_2(\bfq))$, and we can let it act on $\End_{Z(\cH_2(\bfq))}(\cA_2(\bfq))$ by conjugation.
\end{Pt*}

\begin{Lem*}
The embedding $\sA_2(\bfq)|_{\cA_2(\bfq)}$ is $W_0$-equivariant.
\end{Lem*}

\begin{proof}
Indeed, for all $a,b\in\cA_2(\bfq)$ and $w\in W_0$, we have
$$
\sA_2(\bfq)(w(a))(b)=w(a)b=w(aw^{-1}(b))=(waw^{-1})(b)=(w\sA_2(\bfq)(a)w^{-1})(b).
$$
\end{proof}

\section{The generic non-regular spherical representation} \label{IwahoriKtheory}

\subsection{The generic non-regular Iwahori-Hecke algebras} \label{genericgammalagebranonreg}

Let $\gamma=\{\lambda\}\in\bbT^{\vee}/W_0$ be a non-regular orbit. 
As in the regular case, we define a model $\cH_1(\bfq)$ over $\bbZ$ for the component algebra $\cH_{\tilde{\bbZ}}^{\gamma}(\bfq)\subset \cH_{\tilde{\bbZ}}^{(1)}(\bfq)$. The algebra $\cH_1(\bfq)$ will not depend on $\gamma$.

\begin{Pt*}
By construction, the $\tilde{\bbZ}[\bfq]$-algebra $\cH_{\tilde{\bbZ}}^{\gamma}(\bfq)$ admits the following presentation:
$$
\cH_{\tilde{\bbZ}}^{\gamma}(\bfq) =\bigoplus_{w\in W} \bbZ[\bfq] T_w\varepsilon_{\lambda}, 
$$
with
\begin{itemize}
\item braid relations:
$
T_wT_{w'}=T_{ww'}\quad\textrm{for $w,w'\in W$ if $\ell(w)+\ell(w')=\ell(ww')$}
$
\item quadratic relations:
$
T_{\ts}^2=\bfq +(q-1)T_{\ts}\quad \textrm{if $\ts\in S_{\aff}$}.
$
\end{itemize}
\end{Pt*}

\begin{Def*} \label{defgeneric1}
Let $\bfq$ be an indeterminate. The \emph{generic Iwahori-Hecke algebra} is the $\bbZ[\bfq]$-algebra $\cH_1(\bfq)$ defined by generators
$$
\cH_1(\bfq):=\bigoplus_{w\in W} \bbZ[\bfq] T_w
$$
and relations:
\begin{itemize}
\item braid relations:
$
T_wT_{w'}=T_{ww'}\quad\textrm{for $w,w'\in W$ if $\ell(w)+\ell(w')=\ell(ww')$}
$
\item quadratic relations:
$
T_{\ts}^2=\bfq +(\bfq-1)T_{\ts} \quad \textrm{if $\ts\in S_{\aff}$}.
$
\end{itemize}
\end{Def*}

\begin{Pt*}\label{presH1q}
The identity element of $\cH_1(\bfq)$ is $1=T_1$. Moreover we set in $\cH_1(\bfq)$
$$
S:=T_s,\quad U:=T_u\quad\textrm{and}\quad S_0:=T_{s_0}=USU^{-1}.
$$
Then one checks that
$$
\cH_1(\bfq)=\bbZ[\bfq][S,U^{\pm 1}],\quad S^2=\bfq+(\bfq-1)S,\quad U^2S=SU^2
$$
is a presentation of $\cH_1(\bfq)$. Note that the element $U^2$ is invertible in $\cH_1(\bfq)$.
\end{Pt*}

\begin{Pt*}\label{H1VSHgamma}
Sending $1$ to $\varepsilon_{\gamma}$ defines an isomorphism of $\tilde{\bbZ}[\bfq]$-algebras
$$
\xymatrix{
\cH_1(\bfq)\otimes_{\bbZ}\tilde{\bbZ}\ar[r]^<<<<<{\sim} & \cH_{\tilde{\bbZ}}^{\gamma}(\bfq),
}
$$
such that $S\otimes 1\mapsto S\varepsilon_{\gamma}$, $U\otimes 1\mapsto U\varepsilon_{\gamma}$ and $S_0\otimes 1\mapsto S_0\varepsilon_{\gamma}$.
\end{Pt*}

\begin{Pt*} \label{presA1q}
We define $\cA_1(\bfq)\subset\cH_1(\bfq)$ to be the $\bbZ[\bfq]$-subalgebra generated by the elements 
$(S_0-(\bfq-1))U$, $SU$ and $ U^{\pm 2}$. Let $X,Y$ and $z_2$ be indeterminates. 
Then there is a unique $\bbZ[\bfq]$-algebra homomorphism 
$$
\bbZ[\bfq][z_2^{\pm1}][X,Y]/(XY-\bfq z_2) \lra \cA_1(\bfq)
$$
such that 
$X\mapsto (S_0-(\bfq-1))U$, $Y\mapsto SU$, $z_2\mapsto U^2$, and it is an isomorphism. In particular,  $\cA_1(\bfq)$ is a {\it commutative} subalgebra of $\cH_1(\bfq)$. The isomorphism \ref{H1VSHgamma} identifies $\cA_1(\bfq)\otimes_{\bbZ}\tilde{\bbZ}$ with $\cA_{\tilde{\bbZ}}^{\gamma}(\bfq)$. Moreover, permuting $X$ and $Y$ extends to an action of $W_0=\mathfrak{S}_2$ on $\cA_1(\bfq)$ by homomorphisms of $\bbZ[\bfq]$-algebras, whose invariants is the center $Z(\cH_1(\bfq))$ of $\cH_1(\bfq)$ and 
$$
\bbZ[\bfq][z_2^{\pm1}][z_1]\xrightarrow{\sim}\cA_1(\bfq)^{W_0}=Z(\cH_1(\bfq))
$$
with $z_1:=X+Y$. This is a consequence of \ref{H1VSHgamma}, \ref{BernsteinVSIM}, \ref{presAq} and \ref{centergamma}.
 In the following, we will sometimes view the above isomorphisms as identifications. In particular, we will write 
 $$
X=(S_0-(\bfq-1))U=U(S-(\bfq-1)),\quad Y=SU\quad\textrm{and}\quad z_2=U^2\quad\textrm{in}\quad\cH_1(\bfq).
$$
\end{Pt*}
\begin{Pt*}
It is well-known that the generic Iwahori-Hecke algebra $\cH_1(\bfq)$ is a $\bfq$-deformation of the group ring $\bbZ[W]$ of the Iwahori-Weyl group $W=\Lambda\rtimes W_0$. More precisely, specializing the chain of inclusions $\cA_1(\bfq)^{W_0}\subset \cA_1(\bfq) \subset \cH_1(\bfq)$ at $\bfq=1$, yields the chain of inclusions 
$ \bbZ[\Lambda]^{W_0}\subset \bbZ[\Lambda] \subset \bbZ[W].$
\end{Pt*}

\subsection{The Kazhdan-Lusztig-Ginzburg operator}

As in the regular case, we will study the $Z(\cH_1(\bfq))$-algebra $\End_{Z(\cH_1(\bfq))}(\cA_1(\bfq))$
of $Z(\cH_1(\bfq))$-linear endomorphisms of 
$\cA_1(\bfq)$. Recall that $Z(\cH_1(\bfq))= \cA_1(\bfq)^s$ is the subring of invariants of the commutative ring 
$\cA_1(\bfq)$.

\begin{Lem*} \label{1Ybasis}
We have
$$
\cA_1(\bfq)=\cA_1(\bfq)^{s}X\oplus\cA_1(\bfq)^s=\cA_1(\bfq)^{s}\oplus\cA_1(\bfq)^sY
$$
as $\cA_1(\bfq)^{s}$-modules.
\end{Lem*}

\begin{proof}
Applying $s$, the two decompositions are equivalent; so it suffices to check that $\bbZ[z_2^{\pm1}][X,Y]$ is free of rank 2 with basis $1,Y$ over the subring of symmetric polynomials $\bbZ[z_2^{\pm1}][X+Y,XY]$. First if $P=QY$ with $P$ and $Q$ symmetric, then applying $s$ we get $P=QX$ and hence $Q(X-Y)=0$ which implies $P=Q=0$. It remains to check that any monomial $X^iY^j$, $i,j\in\bbN$, belongs to
$$
\bbZ[z_2^{\pm1}][X+Y,XY]+\bbZ[z_2^{\pm1}][X+Y,XY]Y.
$$ 
As $X=(X+Y)-Y$ and $Y^2=-XY+(X+Y)Y$, the later is stable under multiplication by $X$ and $Y$; as it contains $1$, the result follows.
\end{proof}

\begin{Rem*}
The basis $\{1,Y\}$ 
specializes at $\bfq=1$ to the so-called {\it Pittie-Steinberg basis} \cite{St75}  of $\bbZ[\Lambda]$ over 
$\bbZ[\Lambda]^{W_0}$.
\end{Rem*}

\begin{Def*} We let 
$$
D_s:=\textrm{ projector on $\cA_1(\bfq)^sY$ along $\cA_1(\bfq)^{s}$ }
$$
$$
D_s':=\textrm{ projector on $\cA_1(\bfq)^s$ along $\cA_1(\bfq)^{s}X$ }
$$
$$
D_s(\bfq):=D_s-\bfq D_s'.
$$
\end{Def*}

\begin{Rem*}
The operators $D_s$ and $D_s'$ specialize at $\bfq=1$ to the \emph{Demazure operators} on $\bbZ[\Lambda]$, as introduced in \cite{D73,D74}.

\end{Rem*}

\begin{Lem*}\label{relationDsq}
We have 
$$
D_s(\bfq)^2=(1-\bfq)D_s(\bfq)+\bfq.
$$
\end{Lem*}

\begin{proof}
Noting that $Y=z_1-X$, we have
$$
D_s(\bfq)^2(1)=D_s(\bfq)(-\bfq)=\bfq^2=(1-\bfq)(-\bfq)+\bfq=((1-\bfq)D_s(\bfq)+\bfq)(1)
$$
and
\begin{eqnarray*}
D_s(\bfq)^2(Y)&=&D_s(\bfq)(Y-\bfq z_1)\\
&=&Y-\bfq z_1-\bfq z_1(-\bfq)\\
&=&(1-\bfq)(Y-\bfq z_1)+\bfq Y\\
&=&((1-\bfq)D_s(\bfq)+\bfq)(Y).
\end{eqnarray*}
\end{proof}

\subsection{The generic non-regular spherical representation} \label{HIqGL2}

We define the generic non-regular spherical representation of the algebra $\cH_1(\bfq)$ on the $Z(\cH_1(\bfq))$-module
$\cA_1(\bfq)$. The commutative ring
 $\cA_1(\bfq)$ is naturally a subring 
 $$\cA_1(\bfq)\subset \End_{Z(\cH_1(\bfq))}(\cA_1(\bfq)),
 $$
an element $a\in \cA_1(\bfq)$ acting by multiplication $b\mapsto ab$ on $ \cA_1(\bfq)$.

\begin{Th*} \label{sA1q}
There exists a unique $\bbZ[\bfq]$-algebra homomorphism
$$
\xymatrix{
\sA_1(\bfq):\cH_1(\bfq) \ar[r] & \End_{Z(\cH_1(\bfq))}(\cA_1(\bfq))
}
$$
such that 
\begin{itemize}
\item[(i)]
$
\sA_1(\bfq)|_{\cA_1(\bfq)}=\textrm{ the natural inclusion $\cA_1(\bfq)\subset\End_{Z(\cH_1(\bfq))}(\cA_1(\bfq))$}
$
\item[(ii)]
$
\sA_1(\bfq)(S)=-D_s(\bfq).
$
\end{itemize}
\end{Th*}

\begin{proof} 
Recall that $\cH_1(\bfq)=\bbZ[\bfq][S,U^{\pm1}]$ with the relations $S^2=(\bfq-1)S+\bfq$ and $U^2S=SU^2$. In particular
$\sA_1(\bfq)(S):=-D_s(\bfq)$ is well-defined thanks to \ref{relationDsq}. On the other hand, the $\bbZ[\bfq]$-algebra $\cA_1(\bfq)$ is generated by
$$
z_2=U^2,\quad X=US+(1-\bfq)U\quad\textrm{and}\quad Y=SU.
$$
Consequently, there exists a $\bbZ[\bfq]$-algebra homomorphism $\sA_1(\bfq)$ as in the statement of the theorem if and only if there exists
$$
\sA_1(\bfq)(U)\in\End_{Z(\cH_1(\bfq))}(\cA_1(\bfq))
$$
satisfying
\begin{enumerate}
\item $\sA_1(\bfq)(U)^2=\sA_1(\bfq)(U^2)=\sA_1(\bfq)(z_2)=z_2\Id$  \hskip5pt (in particular $\sA_1(\bfq)(U)$ is invertible)
\item $\sA_1(\bfq)(U)(-D_s(\bfq))+(1-\bfq)\sA_1(\bfq)(U)=\textrm{ multiplication by $X$}$
\item $-D_s(\bfq)\sA_1(\bfq)(U)=\textrm{ multiplication by $Y$}$.
\end{enumerate}
Let us use the $Z(\cH_1(\bfq))$-basis $1,Y$ of $\cA_1(\bfq)$ to identify $\End_{Z(\cH_1(\bfq))}(\cA_1(\bfq))$ with the algebra of $2\times 2$-matrices over the ring $Z(\cH_1(\bfq))=\cA_1(\bfq)^s$. Then, by definition,
$$
-D_s(\bfq)=
\left (\begin{array}{cc}
0& 0 \\
0& -1
\end{array} \right)
+
\bfq
\left (\begin{array}{cc}
1& z_1 \\
0& 0
\end{array} \right)
=
\left (\begin{array}{cc}
\bfq& \bfq z_1 \\
0& -1
\end{array} \right).
$$
Moreover, as $X=z_1-Y$, $XY=\bfq z_2$ and $Y^2=-XY+(X+Y)Y=-\bfq z_2+z_1Y$, the multiplications by $X$ and by $Y$ on 
$\cA_1(\bfq)$ get identified with the matrices
$$
\left (\begin{array}{cc}
z_1& \bfq z_2 \\
-1& 0
\end{array} \right)
\quad\textrm{and}\quad
\left (\begin{array}{cc}
0& -\bfq z_2 \\
1& z_1
\end{array} \right).
$$
Now, writing
$$
\sA_1(\bfq)(U)=
\left (\begin{array}{cc}
a& c \\
b& d
\end{array} \right)
$$
we have:
$$
\sA_1(\bfq)(U)^2=z_2\Id \Longleftrightarrow
\left (\begin{array}{cc}
a^2+bc& c(a+d) \\
b(a+d)& d^2+bc
\end{array} \right)
=
\left (\begin{array}{cc}
z_2& 0 \\
0& z_2
\end{array} \right),
$$
$$
\sA_1(\bfq)(U)(-D_s(\bfq))+(1-\bfq)\sA_1(\bfq)(U)=\textrm{ multiplication by $X$}
$$
$$
\Longleftrightarrow
\left (\begin{array}{cc}
a& \bfq(a z_1-c) \\
b&  \bfq(b z_1-d) 
\end{array} \right)
=
\left (\begin{array}{cc}
z_1& \bfq z_2 \\
-1& 0
\end{array} \right)
$$
and
$$
-D_s(\bfq)\sA_1(\bfq)(U)=\textrm{ multiplication by $Y$}
$$
$$
\Longleftrightarrow
\left (\begin{array}{cc}
\bfq(a+z_1b)& \bfq(c+z_1d) \\
-b& -d
\end{array} \right)
=
\left (\begin{array}{cc}
0& -\bfq z_2 \\
1& z_1
\end{array} \right).
$$
Each of the two last systems admits a unique solution, namely
$$
\sA_1(\bfq)(U)=
\left (\begin{array}{cc}
a& c \\
b& d
\end{array} \right)
=
\left (\begin{array}{cc}
z_1& z_1^2-z_2 \\
-1& -z_1
\end{array} \right),
$$
which is also a solution of the first one. Moreover, the determinant
$$
ad-bc=-z_1^2+(z_1^2-z_2)=-z_2
$$
is invertible. 
\end{proof}

\begin{Pt*}
The relation between our generic non-regular representation $\sA_1(\bfq)$ and the theory of Kazhdan-Lusztig \cite{KL87}, and Ginzburg \cite{CG97}, is the following. Introducing a square root $\bfq^{\frac{1}{2}}$ of $\bfq$ and extending scalars along $\bbZ[\bfq]\subset\bbZ[\bfq^{\pm\frac{1}{2}}]$, we obtain the Hecke algebra 
$\cH_1(\bfq^{\pm\frac{1}{2}})$ together with its commutative subalgebra $\cA_1(\bfq^{\pm\frac{1}{2}})$. The latter contains the elements 
$\tilde{\theta}_{\lambda}$, $\lambda\in\Lambda$, introduced by Bernstein and Lusztig, which are defined as follows: writing 
$\lambda=\lambda_1-\lambda_2$ with $\lambda_1,\lambda_2$ antidominant,  one has 
$$
\tilde{\theta}_{\lambda}:=\tilde{T}_{e^{\lambda_1}}\tilde{T}_{e^{\lambda_2}}^{-1}:=\bfq^{-\frac{\ell(\lambda_1)}{2}}\bfq^{\frac{\ell(\lambda_2)}{2}}
T_{e^{\lambda_1}}T_{e^{\lambda_2}}^{-1}.
$$
They are related to the Bernstein basis $\{E(w),\ w\in W\}$ of $\cH_1(\bfq)$ introduced by Vignéras (which is analogous to the Bernstein basis of $\cH^{(1)}(\bfq)$ which we have recalled in \ref {Bernsteinpresprop}) by the formula:
$$
\forall \lambda\in\Lambda,\ \forall w\in W_0,\quad E(e^{\lambda}w)=\bfq^{\frac{\ell(e^{\lambda}w)-\ell(w)}{2}}\tilde{\theta}_{\lambda}T_w\quad\in \cH_1(\bfq)\subset\cH_1(\bfq^{\pm\frac{1}{2}}).
$$
In particular $E(e^{\lambda})=\bfq^{\frac{\ell(e^{\lambda})}{2}}\tilde{\theta}_{\lambda}$, and by the product formula (analogous to the 
product formula for $\cH^{(1)}(\bfq)$, cf. \ref {Bernsteinpresprop}), the $\bbZ[\bfq^{\pm\frac{1}{2}}]$-linear isomorphism
\begin{eqnarray*}
\tilde{\theta}: \bbZ[\bfq^{\pm\frac{1}{2}}][\Lambda] & \iso & \cA_1(\bfq^{\pm\frac{1}{2}})\\
e^{\lambda} & \lmapsto & \tilde{\theta}_{\lambda}
\end{eqnarray*}
is in fact multiplicative, i.e. it is an isomorphism of $\bbZ[\bfq^{\pm\frac{1}{2}}]$-algebras. 

Consequently, if we base change our action map $\sA_1(\bfq)$ to $\bbZ[\bfq^{\pm\frac{1}{2}}]$, we get a representation
$$
\xymatrix{
\sA_1(\bfq^{\pm\frac{1}{2}}):\cH_1(\bfq^{\pm\frac{1}{2}}) \ar[r] & \End_{Z(\cH_1(\bfq^{\pm{\frac{1}{2}}}))}(\cA_1(\bfq^{\pm\frac{1}{2}})) \simeq  \End_{\bbZ[\bfq^{\pm\frac{1}{2}}][\Lambda]^{W_0}}(\bbZ[\bfq^{\pm\frac{1}{2}}][\Lambda]),
}
$$
which coincides with the natural inclusion $\bbZ[\bfq^{\pm\frac{1}{2}}][\Lambda]\subset\End_{\bbZ[\bfq^{\pm\frac{1}{2}}][\Lambda]^{W_0}}(\bbZ[\bfq^{\pm\frac{1}{2}}][\Lambda])$ when restricted to  $\cA_1(\bfq^{\pm\frac{1}{2}})\simeq\bbZ[\bfq^{\pm\frac{1}{2}}][\Lambda]$, and 
which sends $S$ to the opposite $-D_s(\bfq)$ of the $\bfq$-deformed Demazure operator. Hence, modulo our choice of antidominant orientation, this is the spherical representation defined by Kazhdan-Lusztig \cite[Lem. 3.9] {KL87} and Ginzburg \cite[7.6] {CG97}.\footnote{Moreover, it can be checked, in analogy to loc.cit., that the $\cH_1(\bfq)$-module $\sA_1(\bfq)$ is isomorphic to the induction of the trivial character of the finite Hecke (sub)algebra $\bbZ[\bfq][S]$. But we will not make use of this in the following.}

\vskip5pt

In particular, $\sA_1(1)$ is the usual action of the Iwahori-Weyl group $W=\Lambda\rtimes W_0$ on $\Lambda$, and $\sA_1(0)$ can be thought of as a degeneration of the latter.

\end{Pt*}




\begin{Prop*} \label{sAqinj} \label{injIGL2} \label{sA1sf}
The homomorphism $\sA_1(\bfq)$ is injective.
\end{Prop*}

\begin{proof}
It follows from \ref{presH1q} and \ref{presA1q} that the ring $\cH_1(\bfq)$ is generated by the elements
$$
1,\  S,\  U,\  SU
$$
as a module over its center $Z(\cH_1(\bfq))=\bbZ[\bfq][z_1,z_2^{\pm1}]$. As the latter is mapped isomorphically to the center of the matrix algebra $\End_{Z(\cH_1(\bfq))}(\cA_1(\bfq))$ by $\sA_1(\bfq)$, it suffices to check that the images
$$
1,\ \sA_1(\bfq)(S),\ \sA_1(\bfq)(U),\ \sA_1(\bfq)(SU)
$$
of $1, S, U, SU$ by $\sA_1(\bfq)$ are free over $Z(\cH_1(\bfq))$. So let $\alpha,\beta,\gamma,\delta\in Z(\cH_1(\bfq))$ (which is an integral domain) be such that
$$
\alpha
\left (\begin{array}{cc}
1& 0 \\
0& 1
\end{array} \right)
+
\beta
\left (\begin{array}{cc}
\bfq& \bfq z_1 \\
0& -1
\end{array} \right)
+
\gamma
\left (\begin{array}{cc}
z_1& z_1^2-z_2 \\
-1& -z_1
\end{array} \right)
+
\delta
\left (\begin{array}{cc}
0& -\bfq z_2 \\
1& z_1
\end{array} \right)
=0.
$$
Then
$$
\left\{ \begin{array}{lll}
\alpha+\beta\bfq+\gamma z_1&=&0\\
-\gamma+\delta&=&0\\
\beta\bfq z_1+\gamma(z_1^2-z_2)-\delta\bfq z_2&=&0 \\
\alpha-\beta+(\delta-\gamma)z_1&=&0.
\end{array} \right.
$$
We obtain $\delta=\gamma$, $\alpha=\beta$ and
$$
\left\{ \begin{array}{lll}
\alpha(1+\bfq)+\gamma z_1&=&0\\
\alpha\bfq z_1+\gamma(z_1^2-z_2-\bfq z_2)&=&0. \\
\end{array} \right.
$$
The latter system has determinant
$$
(1+\bfq)(z_1^2-z_2-\bfq z_2)-\bfq z_1^2=z_1^2-z_2-2\bfq z_2-\bfq^2z_2
$$
which is nonzero (its specialisation at $\bfq=0$ is equal to $z_1^2-z_2\neq 0$), whence $\alpha=\gamma=0=\beta=\delta$. 
\end{proof}

We record the following two corollaries of the proof.

\begin{Cor*}\label{dimOverCenterH1} 
The ring $\cH_1(\bfq)$ is a free $Z(\cH_1(\bfq))$-module on the basis $1, S, U, SU$.
\end{Cor*}

\begin{Cor*}\label{faithfulatzero} 
The homomorphism $\sA_1(0)$ is injective.
\end{Cor*}

\begin{Pt*}
We end this section by noting an equivariance property of $\sA_1(\bfq)$. As already noticed, the finite Weyl group $W_0$ acts on $\cA_1(\bfq)$ by $\bbZ[\bfq]$-algebra automorphisms, and the action is clearly faithful. Moreover $\cA_1(\bfq)^{W_0}=Z(\cH_1(q))$. Hence $W_0$ can be viewed as a subgroup of $\End_{Z(\cH_1(\bfq))}(\cA_1(\bfq))$, and we can let it act on $\End_{Z(\cH_1(\bfq))}(\cA_1(\bfq))$ by conjugation.
\end{Pt*}

\begin{Lem*}
The embedding $\sA_1(\bfq)|_{\cA_1(\bfq)}$ is $W_0$-equivariant.
\end{Lem*}

\begin{proof}
Indeed, for all $a,b\in\cA_1(\bfq)$ and $w\in W_0$, we have
$$
\sA_1(\bfq)(w(a))(b)=w(a)b=w(aw^{-1}(b))=(waw^{-1})(b)=(w\sA_1(\bfq)(a)w^{-1})(b).
$$
\end{proof}

\section{$K$-theory of the dual flag variety}

\subsection{The Vinberg monoid of the dual group $\mathbf{\whG}=\mathbf{GL_2}$}\label{Vinsubsection}

\begin{Pt*}
The Langlands dual group over $k:=\overline{\bbF}_q$ of the connected reductive algebraic group $GL_2$ over $F$ is $\mathbf{\whG}=\mathbf{GL_2}$. We recall the $k$-monoid scheme introduced by Vinberg in \cite{V95}, in the particular case of $\mathbf{GL_2}$. It is in fact defined over 
$\bbZ$, as the group $\mathbf{GL_2}$. In the following, all the fiber products are taken over the base ring $\bbZ$. \end{Pt*}

\begin{Def*}
Let $\Mat_{2\times 2}$ be the $\bbZ$-monoid scheme of $2\times 2$-matrices (with usual matrix multiplication as operation). The \emph{Vinberg monoid for $\mathbf{GL_2}$} is the 
$\bbZ$-monoid scheme
$$
V_{\mathbf{GL_2}}:=\Mat_{2\times 2}\times \bbG_m.
$$ 
\end{Def*}

\begin{Pt*}
The group $\mathbf{GL_2}\times\bbG_m$ is recovered from the monoid $V_{\mathbf{GL_2}}$ as its group of units. The group 
$\mathbf{GL_2}$ itself is recovered as follows. Denote by $z_2$ the canonical coordinate on $\bbG_m$. Then \emph{let $\bfq$ be the homomorphism from $V_{\mathbf{GL_2}}$ to the multiplicative monoid $(\bbA^1,\cdot )$ defined by $(f,z_2)\mapsto \det(f)z_2^{-1}$}:
$$
\xymatrix{
V_{\mathbf{GL_2}} \ar[d]_{\bfq} \\
\bbA^1.
}
$$
Then $\mathbf{GL_2}$ is recovered as the fiber at $\bfq=1$, canonically: 
$$ 
\bfq^{-1}(1)= \{ (f,z_2) : \det(f)=z_2 \} \iso \mathbf{GL_2},\quad (f,z_2)\mapsto f.
$$

The fiber at $\bfq=0$ is the $\bbZ$-semigroup scheme
$$
V_{\mathbf{GL_2},0} := \bfq^{-1}(0)= \Sing_{2\times 2}\times\bbG_m,
$$
where $\Sing_{2\times 2}$ represents the singular $2\times 2$-matrices. Note that it has no identity element, i.e. it is a semigroup which is not a monoid.  
\end{Pt*}

\begin{Pt*}\label{VT}
Let $\Diag_{2\times 2}\subset \Mat_{2\times 2}$ be the submonoid scheme of diagonal $2\times 2$-matrices, and set
$$
V_{\mathbf{\whT}}:=\Diag_{2\times 2}\times\bbG_m\subset V_{\mathbf{GL_2}}= \Mat_{2\times 2}\times\bbG_m.
$$
This is a diagonalizable $\bbZ$-monoid scheme with character monoid 
$$
\bbX^{\bullet}(V_{\mathbf{\whT}})=\bbN(1,0)\oplus\bbN(0,1)\oplus\bbZ \subset\bbZ(1,0)\oplus\bbZ(0,1)\oplus\bbZ=\Lambda\oplus \bbZ=\bbX^{\bullet}(\mathbf{\whT})\oplus\bbX^{\bullet}(\bbG_m).
$$
In particular, setting $X:=e^{(1,0)}$ and $Y:=e^{(0,1)}$ in the group ring $\bbZ[\Lambda]$, we have
$$
\mathbf{\whT}=\Spec(\bbZ[X^{\pm1},Y^{\pm1}])\subset \Spec(\bbZ[z_2^{\pm1}][X,Y])=V_{\mathbf{\whT}}.
$$
Again, this closed subgroup is recovered as the fiber at $\bfq=1$ of the fibration $\bfq|_{V_{\mathbf{\whT}}} : V_{\mathbf{\whT}}\ra \bbA^1$, and the fiber at $\bfq=0$ is the $\bbZ$-semigroup scheme $\SingDiag_{2\times2}\times\bbG_m$ where $\SingDiag_{2\times2}$ represents the singular diagonal $2\times 2$-matrices:
$$
\xymatrix{
\mathbf{\whT}\ar@{^{(}->}[r] \ar[d] & V_{\mathbf{\whT}} \ar[d]_{\bfq} & \ar@{_{(}->}[l]\SingDiag_{2\times2}\times\bbG_m \ar[d]\\
\Spec(\bbZ)\ar@{^{(}->}[r]^<<<<<1 &\bbA^1 & \ar@{_{(}->}[l]_>>>>>>>>>0  \Spec(\bbZ).
}
$$
In terms of equations, the $\bbA^1$-family 
$$
\xymatrix{
\bfq:V_{\mathbf{\whT}}=\Diag_{2\times 2}\times\bbG_m=\Spec(\bbZ[z_2^{\pm1}][X,Y])\ar[r] &\bbA^1
}
$$ 
is given by the formula $\bfq(\diag(x,y),z_2)=\det(\diag(x,y))z_2^{-1}=xyz_2^{-1}$. Hence, after fixing $z_2\in\bbG_m$, the fiber over a point $\bfq\in\bbA^1$ is the hyperbola $xy=\bfq z_2$, which is non-degenerate if $\bfq \neq 0$, and is the union of the two coordinate axis if $\bfq =0$. 
\end{Pt*}

\subsection{The associated flag variety and its equivariant $K$-theory}

\begin{Pt*}\label{ref_to_virtual}
Let $\mathbf{\whB}\subset\mathbf{GL_2}$ be the Borel subgroup of upper triangular matrices, let $\UpTriang_{2\times 2}$ be the $\bbZ$-monoid scheme representing the upper triangular $2\times 2$-matrices, and set
$$
V_{\mathbf{\whB}}:=\UpTriang_{2\times 2}\times\bbG_m\subset \Mat_{2\times 2}\times\bbG_m=:V_{\mathbf{GL_2}}.
$$
Then we can apply to this inclusion of $\bbZ$-monoid schemes the general formalism developed in \cite{PS20}. In particular, the \emph{flag variety $V_{\mathbf{GL_2}}/V_{\mathbf{\whB}}$} is defined as a $\bbZ$-monoidoid. Moreover, after base changing along 
$\bbZ\ra k$, we have defined a ring $K^{V_{\mathbf{GL_2}}}(V_{\mathbf{GL_2}}/V_{\mathbf{\whB}})$ of $V_{\mathbf{GL_2}}$-equivariant $K$-theory on the flag variety, together with an induction isomorphism 
$$
\xymatrix{
\cI nd_{V_{\mathbf{\whB}}}^{V_{\mathbf{GL_2}}}:R(V_{\mathbf{\whB}}) \ar[r]^<<<<<{\sim} & K^{V_{\mathbf{GL_2}}}(V_{\mathbf{GL_2}}/V_{\mathbf{\whB}})
}
$$
from the ring $R(V_{\mathbf{\whB}})$ of right representations of the $k$-monoid scheme $V_{\mathbf{\whB}}$ on finite dimensional $k$-vector spaces.
\end{Pt*}

\begin{Pt*}\label{retraction}
Now, we have the inclusion of monoids $V_{\mathbf{\whT}}=\Diag_{2\times 2}\times\bbG_m\subset V_{\mathbf{\whB}}=\UpTriang_{2\times 2}\times\bbG_m$, which admits the retraction
\begin{eqnarray*}
V_{\mathbf{\whB}} & \lra & V_{\mathbf{\whT}} \\
\bigg(\left (\begin{array}{cc}
x & c\\
0 & y
\end{array} \right),\ z_2\bigg)
&
\lmapsto 
&
\bigg(\left (\begin{array}{cc}
x & 0\\
0 & y
\end{array} \right),\ z_2\bigg).
\end{eqnarray*}
Let $\Rep(V_{\mathbf{\whT}})$ be the category of representations of the commutative $k$-monoid scheme $V_{\mathbf{\whT}}$ on finite dimensional $k$-vector spaces. The above preceding inclusion and retraction define a \emph{restriction functor} and an \emph{inflation functor}
$$
\xymatrix{
\Res_{V_{\mathbf{\whT}}}^{V_{\mathbf{\whB}}}:\Rep(V_{\mathbf{\whB}}) \ar@<1ex>[r] & \Rep(V_{\mathbf{\whT}}) :\Infl_{V_{\mathbf{\whT}}}^{V_{\mathbf{\whB}}}. \ar@<1ex>[l] 
}
$$
These functors are exact and compatible with the tensors products and units.
\end{Pt*}

\begin{Lem*}\label{ResInfl}
The ring homomorphisms
$$
\xymatrix{
\Res_{V_{\mathbf{\whT}}}^{V_{\mathbf{\whB}}}:R(V_{\mathbf{\whB}}) \ar@<1ex>[r] & R(V_{\mathbf{\whT}}) :\Infl_{V_{\mathbf{\whT}}}^{V_{\mathbf{\whB}}} \ar@<1ex>[l]
}
$$
are isomorphisms, which are inverse one to the other.
\end{Lem*}

\begin{proof}
We have $\Res_{V_{\mathbf{\whT}}}\circ \Infl_{V_{\mathbf{\whT}}}^{V_{\mathbf{\whB}}}=\Id$ by construction. Conversely, let $M$ be an object of $\Rep(V_{\mathbf{\whB}})$. The solvable subgroup $\mathbf{\whB}\times\bbG_m\subset V_{\mathbf{\whB}}$ stabilizes a line $L\subseteq M$. As $\mathbf{\whB}\times\bbG_m$ is dense in $V_{\mathbf{\whB}}$, the line $L$ is automatically $V_{\mathbf{\whB}}$-stable. Moreover the unipotent radical $\mathbf{\whU}\subset \mathbf{\whB}$ acts trivially on $L$, so that $\mathbf{\whB}\times\bbG_m$ acts on $L$ through the quotient $\mathbf{\whT}\times\bbG_m$. Hence, by density again, $V_{\mathbf{\whB}}$ acts on $L$ through the retraction $V_{\mathbf{\whB}}\ra V_{\mathbf{\whT}}$.
This shows that any irreducible $M$ is a character inflated from a character of $V_{\mathbf{\whT}}$. In particular, the 
map $R(V_{\mathbf{\whT}}) \rightarrow R(V_{\mathbf{\whB}})$ is surjective and hence bijective.
\end{proof}

\begin{Cor*}\label{cVGL2}
We have a ring isomorphism
$$
\xymatrix{
c_{V_{\mathbf{GL_2}}}:=\cI nd_{V_{\mathbf{\whB}}}^{V_{\mathbf{GL_2}}}\circ\Infl_{V_{\mathbf{\whT}}}^{V_{\mathbf{\whB}}}:\bbZ[X,Y,z_2^{\pm1}]\cong R(V_{\mathbf{\whT}}) \ar[r]^<<<<<{\sim} & K^{V_{\mathbf{GL_2}}}(V_{\mathbf{GL_2}}/V_{\mathbf{\whB}}),
}
$$
that we call \emph{the characteristic isomorphism in the equivariant $K$-theory of the flag variety $V_{\mathbf{GL_2}}/V_{\mathbf{\whB}}$}.
\end{Cor*}

\begin{Pt*}
We have a commutative diagram \emph{specialization at $\bfq=1$}
$$
\xymatrix{
\bbZ[X,Y,z_2^{\pm1}] \ar[rr]^{c_{V_{\mathbf{GL_2}}}}_{\sim}   \ar@{->>}[d] && K^{V_{\mathbf{GL_2}}}(V_{\mathbf{GL_2}}/V_{\mathbf{\whB}})  \ar@{->>}[d] \\
\bbZ[X^{\pm1},Y^{\pm1}]\ar[rr]^{c_{\mathbf{GL_2}}}_{\sim}   && K^{\mathbf{GL_2}}(\mathbf{GL_2}/\mathbf{\whB}).
}
$$
The vertical map on the left-hand side is given by specialization $\bfq = 1$, i.e. by the surjection 
$$
\bbZ[X,Y,z_2^{\pm 1 }] = \bbZ[\bfq][X,Y,z_2^{\pm 1 }] / (XY-\bfq z_2) \longrightarrow \bbZ[X,Y,z_2^{\pm 1 }] / (XY-z_2) = \bbZ[X^{\pm 1},Y^{\pm 1}].
$$
The vertical map on the right-hand side is given by restricting equivariant vector bundles to the $1$-fiber of $\bfq: V_{\mathbf{GL_2}}\ra\bbA^1$, thereby recovering the classical theory. 

\end{Pt*}

\begin{Pt*}
Let $\Rep(V_{\mathbf{GL_2}})$ be the category of right representations of the $k$-monoid scheme $V_{\mathbf{GL_2}}$ on finite dimensional $k$-vector spaces. The inclusion $V_{\mathbf{\whB}}\subset V_{\mathbf{GL_2}}$ defines a restriction functor 
$$
\xymatrix{
\Res^{V_{\mathbf{GL_2}}}_{V_{\mathbf{\whB}}}:\Rep(V_{\mathbf{GL_2}}) \ar[r] & \Rep(V_{\mathbf{\whB}}),
}
$$
whose composition with $\Res^{V_{\mathbf{\whB}}}_{V_{\mathbf{\whT}}}$ is the restriction from $V_{\mathbf{GL_2}}$ to $V_{\mathbf{\whT}}$:
$$
\xymatrix{
\Res^{V_{\mathbf{GL_2}}}_{V_{\mathbf{\whT}}}=\Res^{V_{\mathbf{\whB}}}_{V_{\mathbf{\whT}}}\circ \Res^{V_{\mathbf{GL_2}}}_{V_{\mathbf{\whB}}}:\Rep(V_{\mathbf{GL_2}}) \ar[r] & \Rep(V_{\mathbf{\whT}}).
}
$$
These restriction functors are exact and compatible with the tensors products and units.
\end{Pt*}

\begin{Pt*}
The action of the Weyl group $W_0$ on $\bbX^{\bullet}(\mathbf{\whT})\oplus\bbX^{\bullet}(\bbG_m)$ (trivial on $\bbX^{\bullet}(\bbG_m)$) stabilizes $\bbX^{\bullet}(V_{\mathbf{\whT}})$. Consequently $W_0$ acts on $V_{\mathbf{\whT}}$ and the inclusion $\mathbf{\whT}\subset V_{\mathbf{\whT}}$ is $W_0$-equivariant. Explicitly, $W_0=\{1,s\}$ and $s$ acts on $V_{\mathbf{\whT}}=\Diag_{2\times 2}\times\bbG_m$ by permuting the two diagonal entries and trivially on the $\bbG_m$-factor.
\end{Pt*}

\begin{Lem*}\label{charVin}
The ring homomorphism 
$$
\xymatrix{
\Res^{V_{\mathbf{GL_2}}}_{V_{\mathbf{\whT}}}:R(V_{\mathbf{GL_2}}) \ar[r] & R(V_{\mathbf{\whT}})
}
$$
is injective, with image the subring $R(V_{\mathbf{\whT}})^{W_0}\subset R(V_{\mathbf{\whT}})$ of $W_0$-invariants. The resulting ring isomorphism
$$
\xymatrix{
\chi_{V_{\mathbf{GL_2}}}:R(V_{\mathbf{GL_2}}) \ar[r]^<<<<<{\sim} & R(V_{\mathbf{\whT}})^{W_0}
}
$$
is the \emph{character isomorphism of $V_{\mathbf{GL_2}}$}.
\end{Lem*}

\begin{proof}
This is a general result on the representation theory of $V_{\mathbf{\whG}}$. Note that in the case of $\mathbf{\whG}=\mathbf{GL_2}$, we have
$$
R(V_{\mathbf{\whT}})^{W_0}=\bbZ[X+Y,XYz_2^{-1}=:\bfq,z_2^{\pm1}]\subset \bbZ[X,Y,z_2^{\pm1}]=R(V_{\mathbf{\whT}}).
$$
\end{proof}

\section{Dual parametrization of generic Hecke modules}

We keep all the notations introduced in the preceding section. In particular, $k=\overline{\bbF}_q$. 

\subsection{The generic Bernstein isomorphism}

Recall from \ref{AIH} the subring $\cA(\bfq)\subset\cH^{(1)}(\bfq)$ and the remarkable Bernstein basis elements $E(1,0)$, $E(0,1)$ and $E(1,1)$. Also recall from \ref{VT} the representation ring $R(V_{\mathbf{\whT}})=\bbZ[X,Y,z_2^{\pm1}]$ of the diagonalizable $k$-submonoid scheme $V_{\mathbf{\whT}}\subset V_{\mathbf{\whG}}$ of the Vinberg $k$-monoid scheme of the Langlands dual $k$-group $\mathbf{\whG}=\mathbf{GL_2}$ of $GL_{2,F}$.

\begin{Th*}\label{genB}
There exists a unique ring homomorphism
$$
\xymatrix{
\sB(\bfq):\cA(\bfq) \ar[r]& R(V_{\mathbf{\whT}})
} 
$$
such that
$$
\sB(\bfq)(E(1,0))=X,\quad \sB(\bfq)(E(0,1))=Y,\quad \sB(\bfq)(E(1,1))=z_2\quad\textrm{and}\quad\sB(\bfq)(\bfq)=XYz_2^{-1}.
$$
It is an isomorphism.
\end{Th*}

\begin{proof}
This is a reformulation of the first part of \ref{presAq}.
\end{proof}

\begin{Pt*}\label{Repring_of_bbT_dual}
Then recall from \ref{AIH} the subring $\cA^{(1)}(\bfq)=\bbZ[\bbT]\otimes_{\bbZ}\cA(\bfq)\subset\cH^{(1)}(\bfq)$ where $\bbT$ is the finite abelian group $\bfT(\bbF_q)$. Let $\bbT^{\vee}$ be the finite abelian dual group of $\bbT$. As $\bbT^{\vee}$ has order prime to $p$, it defines a constant finite diagonalizable $k$-group scheme, whose group of characters is $\bbT$, and hence whose representation ring $R(\bbT^{\vee})$ identifies with $\bbZ[\bbT]$: $t\in\bbT\subset \bbZ[\bbT]$ corresponds to the character $\ev_t$ of $\bbT^{\vee}$ given by evaluation at $t$. 
Set
$$V^{(1)}_{\mathbf{\whT}}:=\bbT^{\vee}\times V_{\mathbf{\whT}}.$$ 
\end{Pt*}

\begin{Cor*}\label{genB1}
There exists a unique ring homomorphism
$$
\xymatrix{
\sB^{(1)}(\bfq):\cA^{(1)}(\bfq)\ar[r] & R(V^{(1)}_{\mathbf{\whT}})
} 
$$
such that
$$
\sB^{(1)}(\bfq)(E(1,0))=X,\quad \sB^{(1)}(\bfq)(E(0,1))=Y,\quad \sB^{(1)}(\bfq)(E(1,1))=z_2,\quad\sB^{(1)}(\bfq)(\bfq)=XYz_2^{-1}
$$
$$
\textrm{and}\quad\forall t\in\bbT,\ \sB^{(1)}(\bfq)(T_t)=\ev_t.
$$
It is an isomorphism, that we call the \emph{generic (pro-$p$) Bernstein isomorphism}.
\end{Cor*}

\begin{Pt*}\label{B1qchar}
Also, setting $V^{(1)}_{\mathbf{\whB}}:=\bbT^{\vee}\times V_{\mathbf{\whB}}$, we have from \ref{ResInfl} the ring isomorphism
$$
\xymatrix{
\Infl^{V^{(1)}_{\mathbf{\whB}}}_{V^{(1)}_{\mathbf{\whT}}}=\Id_{\bbZ[\bbT]}\otimes_{\bbZ}\Res^{V_{\mathbf{\whB}}}_{V_{\mathbf{\whT}}}:
R(V^{(1)}_{\mathbf{\whT}})=\bbZ[\bbT]\otimes_{\bbZ}R(V_{\mathbf{\whT}})\ar[r]^<<<<<{\sim} &
R(V^{(1)}_{\mathbf{\whB}})=\bbZ[\bbT]\otimes_{\bbZ}R(V_{\mathbf{\whB}}),
}
$$
and setting $V^{(1)}_{\mathbf{\whG}}:=\bbT^{\vee}\times V_{\mathbf{\whG}}$, we have from \cite[2.5.2]{PS20}, the ring isomorphism
$$
\xymatrix{
\cI nd_{V^{(1)}_{\mathbf{\whB}}}^{V^{(1)}_{\mathbf{\whG}}}:R(V^{(1)}_{\mathbf{\whB}}) \ar[r]^<<<<<{\sim} & K^{V^{(1)}_{\mathbf{\whG}}}(V^{(1)}_{\mathbf{\whG}}/V^{(1)}_{\mathbf{\whB}});
}
$$
hence by composition we get the \emph{characteristic isomorphism} 
$$
\xymatrix{
c_{V^{(1)}_{\mathbf{\whG}}}:R(V^{(1)}_{\mathbf{\whT}}) \ar[r]^<<<<<{\sim} & K^{V^{(1)}_{\mathbf{\whG}}}(V^{(1)}_{\mathbf{\whG}}/V^{(1)}_{\mathbf{\whB}}).
}
$$
Whence a ring isomorphism
$$
\xymatrix{
c_{V^{(1)}_{\mathbf{\whG}}}\circ \sB^{(1)}(\bfq):\cA^{(1)}(\bfq) \ar[r]^<<<<<{\sim} & K^{V^{(1)}_{\mathbf{\whG}}}(V^{(1)}_{\mathbf{\whG}}/V^{(1)}_{\mathbf{\whB}}).
} 
$$
\end{Pt*}

\begin{Pt*}\label{SpecB1q}
The representation ring $R(V_{\mathbf{\whT}})$ is canonically isomorphic to the ring $\bbZ[V_{\mathbf{\whT}}]$ of regular functions of $V_{\mathbf{\whT}}$ considered now as a diagonalizable monoid scheme over $\bbZ$. Also recall from \ref{finiteT} the ring extension $\bbZ\subset\tilde{\bbZ}$, and denote by $\tilde{\bullet}$ the base change functor from $\bbZ$ to $\tilde{\bbZ}$. For example, we will from now on write $\tilde{\cA}^{(1)}(\bfq)$ instead of $\cA_{\tilde{\bbZ}}^{(1)}(\bfq)$. We have the constant finite diagonalizable $\tilde{\bbZ}$-group scheme $\bbT^{\vee}$, whose group of characters is $\bbT$, and whose ring of regular functions is 
$$
\tilde{\bbZ}[\bbT]=\prod_{\lambda\in \bbT^{\vee}}\tilde{\bbZ}\varepsilon_{\lambda}.
$$
Hence applying the functor $\Spec$ to $\tilde{\sB}^{(1)}(\bfq)$, we obtain the commutative diagram of 
$\tilde{\bbZ}$-schemes
$$
\xymatrix{
\Spec(\tilde{\cA}^{(1)}(\bfq)) \ar[dr]_{\pi_0\times\bfq} &  & V^{(1)}_{\mathbf{\whT}}=\bbT^{\vee}\times V_{\mathbf{\whT}} \ar[dl]^{\Id\times\bfq} \ar[ll]_>>>>>>>>>>>>>>>>>>>>>>>>{\Spec(\tilde{\sB}^{(1)}(\bfq))}^>>>>>>>>>>>>>>>>>>>>>>>>{\sim} \\
&(\bbA^1)^{(1)}:= \bbT^{\vee}\times\bbA^1 &
}
$$
where $\pi_0:\Spec(\tilde{\cA}^{(1)}(\bfq))\ra \bbT^{\vee}$ is the decomposition of $\Spec(\tilde{\cA}^{(1)}(\bfq))$ into its connected components. 
In particular, for each $\lambda\in \bbT^{\vee}$, we have the subring $\tilde{\cA}^{\lambda}(\bfq)= \tilde{\cA}^{(1)}(\bfq)\varepsilon_\lambda$ of 
$\tilde{\cA}^{(1)}(\bfq)$ and the isomorphism 
$$
\xymatrix{
\Spec(\tilde{\cA}^{\lambda}(\bfq))  & &\{\lambda\} \times V_{\mathbf{\whT}} \ar[ll]_>>>>>>>>>>>{\Spec(\tilde{\sB}^{\lambda}(\bfq))}^>>>>>>>>>>>{\sim}
}
$$
of $\tilde{\bbZ}$-schemes over $\{\lambda\}\times\bbA^1$. In turn, each of these isomorphisms admits a model over $\bbZ$, obtained by applying $\Spec$ to the ring isomorphism in \ref{presA1q}
$$
\xymatrix{
\sB_1(\bfq):\cA_1(\bfq)\ar[r]^<<<<<{\sim} & R(V_{\mathbf{\whT}}).
} 
$$
\end{Pt*}

\subsection{The generic Satake isomorphism}

Recall part of our notation: $\mathbf{G}$ is the algebraic group $\mathbf{GL_2}$ (which is defined over $\bbZ$), $F$ is a local field and $G:=\mathbf{G}(F)$. We have denoted by $o_F$ the ring of integers of $F$. Now we set $K:=\mathbf{G}(o_F)$.

\begin{Def*} \label{defsph}
Let $R$ be any commutative ring. The \emph{spherical Hecke algebra of $G$ with coefficients in $R$} is defined to be the convolution algebra
$$
\cH_{R}^{\sph}:=(R[K\backslash G/K],\star)
$$
generated by the $K$-double cosets in $G$. 
\end{Def*}

\begin{Pt*}
By the work of Kazhdan and Lusztig, the $R$-algebra $\cH_{R}^{\sph}$ depends on $F$ only through the cardinality $q$ of its residue field. Indeed, choose a uniformizer $\varpi\in o_F$. For a dominant cocharacter $\lambda\in\Lambda^+$ of $\mathbf{T}$, let $\mathbbm{1}_{\lambda}$ be the characteristic function of the double coset $K\lambda(\varpi)K$. Then $(\mathbbm{1}_{\lambda})_{\lambda\in \Lambda^+}$ is an $R$-basis of $\cH_{R}^{\sph}$. Moreover, for all 
$\lambda,\mu,\nu\in\Lambda^+$, there exist polynomials
$$
N_{\lambda,\mu;\nu}(\bfq)\in\bbZ[\bfq]
$$
depending only on the triple $(\lambda,\mu,\nu)$, such that
$$
\mathbbm{1}_{\lambda}\star\mathbbm{1}_{\mu}=\sum_{\nu\in\Lambda^+}N_{\lambda,\mu;\nu}(q)\mathbbm{1}_{\nu}
$$
where $N_{\lambda,\mu;\nu}(q)\in\bbZ\subset R$ is the value of $N_{\lambda,\mu;\nu}(\bfq)$ at $\bfq=q$. These polynomials are uniquely determined by this property since when $F$ vary, the corresponding integers $q$ form an infinite set. Their existence can be deduced from the theory of the spherical algebra with coefficients in $\bbC$, as $\cH_{R}^{\sph}=R\otimes_{\bbZ}\cH_{\bbZ}^{\sph}$ and $\cH_{\bbZ}^{\sph}\subset \cH_{\bbC}^{\sph}$ (e.g. using arguments similar to those in the proof of \ref{ThgenSat} below). 
\end{Pt*}

\begin{Def*}
Let $\bfq$ be an indeterminate. The \emph{generic spherical Hecke algebra} is the $\bbZ[\bfq]$-algebra $\cH^{\sph}(\bfq)$ defined by generators
$$
\cH^{\sph}(\bfq):=\oplus_{\lambda\in\Lambda^+}\bbZ[\bfq]T_{\lambda}
$$
and relations:
$$
 T_{\lambda}T_{\mu}=\sum_{\nu\in\Lambda^+}N_{\lambda,\mu;\nu}(\bfq)T_{\nu}\quad\textrm{for all $\lambda,\mu\in\Lambda^+$}.
$$
\end{Def*}

\begin{Th*}\label{ThgenSat}
There exists a unique ring homomorphism
$$
\xymatrix{
\sS(\bfq):\cH^{\sph}(\bfq)\ar[r] & R(V_{\mathbf{\whT}})
}
$$
such that
$$
\sS(\bfq)(T_{(1,0)})=X+Y,\quad \sS(\bfq)(T_{(1,1)})=z_2\quad\textrm{and}\quad\sS(\bfq)(\bfq)=XYz_2^{-1}.
$$
It is an isomorphism onto the subring $R(V_{\mathbf{\whT}})^{W_0}$ of $W_0$-invariants
$$
\xymatrix{
\sS(\bfq):\cH^{\sph}(\bfq)\ar[r]^>>>>>{\sim} & R(V_{\mathbf{\whT}})^{W_0} \subset R(V_{\mathbf{\whT}}).
}
$$
In particular, the algebra $\cH^{\sph}(\bfq)$ is commutative. 
\end{Th*}

\begin{proof} 
Let
$$
\xymatrix{
\sS_{\cl}:\cH_{\bbC}^{\sph}\ar[r]^{\sim} & \bbC[\bbX^{\bullet}(\mathbf{\whT})]^{W_0}
}
$$
be the `classical' isomorphism constructed by Satake \cite{S63}. We use \cite{Gr98} as a reference. 

For $\lambda\in\Lambda^+$, let $\chi_{\lambda}\in \bbZ[\bbX^{\bullet}(\mathbf{\whT})]^{W_0}$ be the character of the irreducible representation of $\mathbf{\whG}$ of highest weight $\lambda$. Then $(\chi_{\lambda})_{\lambda\in\Lambda^+}$ is a $\bbZ$-basis of 
$\bbZ[\bbX^{\bullet}(\mathbf{\whT})]^{W_0}$. Set $f_{\lambda}:=\sS_{\cl}^{-1}(q^{\lan\rho,\lambda\ran}\chi_{\lambda})$, where 
$2\rho=\alpha:=(1,-1)$. Then for each $\lambda,\mu\in\Lambda^+$, there exist polynomials $d_{\lambda,\mu}(\bfq)\in\bbZ[\bfq]$ such that
$$
f_{\lambda}=\mathbbm{1}_{\lambda}+\sum_{\mu<\lambda}d_{\lambda,\mu}(q)\mathbbm{1}_{\mu}\in \cH_{\bbC}^{\sph},
$$
where $d_{\lambda,\mu}(q)\in\bbZ$ is the value of $d_{\lambda,\mu}(\bfq)$ at $\bfq=q$; the polynomial $d_{\lambda,\mu}(\bfq)$ depends only on the couple $(\lambda,\mu)$, in particular it is uniquely determined by this property. As $(\mathbbm{1}_{\lambda})_{\lambda\in \Lambda^+}$ is a $\bbZ$-basis of $\cH_{\bbZ}^{\sph}$, so is $(f_{\lambda})_{\lambda\in\Lambda^+}$. Then let us set
$$
f_{\lambda}(\bfq):=T_{\lambda}+\sum_{\mu<\lambda}d_{\lambda,\mu}(\bfq)T_{\mu}\in \cH^{\sph}(\bfq).
$$
As $(T_{\lambda})_{\lambda\in \Lambda^+}$ is a $\bbZ[\bfq]$-basis of $\cH^{\sph}(\bfq)$, so is $(f_{\lambda}(\bfq))_{\lambda\in\Lambda^+}$.

Next consider the following $\bbZ[\bfq^{\frac{1}{2}}]$-linear map:
\begin{eqnarray*}
\sS_{\cl}(\bfq):\bbZ[\bfq^{\frac{1}{2}}]\otimes_{\bbZ[\bfq]}\cH^{\sph}(\bfq)&\lra& \bbZ[\bfq^{\frac{1}{2}}]\otimes_{\bbZ}\bbZ[\bbX^{\bullet}(\mathbf{\whT})]=\bbZ[\bfq^{\frac{1}{2}}][\bbX^{\bullet}(\mathbf{\whT})] \\
1\otimes f_{\lambda}(\bfq) & \lmapsto & \bfq^{\lan\rho,\lambda\ran}\chi_{\lambda}.
\end{eqnarray*}
We claim that it is a ring homomorphism. Indeed, for $h_1(\bfq),h_2(\bfq)\in \bbZ[\bfq^{\frac{1}{2}}]\otimes_{\bbZ[\bfq]}\cH^{\sph}(\bfq)$, we need to check the identity 
$$
\sS_{\cl}(\bfq)(h_1(\bfq)h_2(\bfq))=\sS_{\cl}(\bfq)(h_1(\bfq))\sS_{\cl}(\bfq)(h_2(\bfq))\in \bbZ[\bfq^{\frac{1}{2}}][\bbX^{\bullet}(\mathbf{\whT})].
$$
Projecting in the $\bbZ[\bfq^{\frac{1}{2}}]$-basis $\bbX^{\bullet}(\mathbf{\whT})$, the latter corresponds to (a finite number of) identities in the ring $\bbZ[\bfq^{\frac{1}{2}}]$ of polynomials in the variable $\bfq^{\frac{1}{2}}$. Now, by construction and because $\sS_{\cl}$ is a ring homomorphism, the desired identities hold after specialyzing $\bfq$ to any power of a prime number; hence they hold in $\bbZ[\bfq^{\frac{1}{2}}]$. Also note that $\sS_{\cl}(\bfq)$ maps $1=T_{(0,0)}$ to $1=\chi_{(0,0)}$ by definition. 

It can also be seen that $\sS_{\cl}(\bfq)$ is injective using a specialization argument: if $h(\bfq)\in\bbZ[\bfq^{\frac{1}{2}}]\otimes_{\bbZ[\bfq]}\cH^{\sph}(\bfq)$ satisfies $\sS_{\cl}(\bfq)(h(\bfq))=0$, then the coordinates of $h(\bfq)$ (in the basis $(1\otimes f_{\lambda}(\bfq))_{\lambda\in\Lambda^+}$ say, one can also use the basis $(1\otimes T_{\lambda})_{\lambda\in\Lambda^+}$) are polynomials in the variable $\bfq^{\frac{1}{2}}$ which must vanish for an infinite number of values of $\bfq$, and hence they are identically zero. 

Let us describe the image of $\cH^{\sph}(\bfq)\subset \bbZ[\bfq^{\frac{1}{2}}]\otimes_{\bbZ[\bfq]}\cH^{\sph}(\bfq)$ under the ring embedding $\sS_{\cl}(\bfq)$. By construction, we have
$$
\sS_{\cl}(\bfq)(\cH^{\sph}(\bfq))=\bigoplus_{\lambda\in\Lambda^+}\bbZ[\bfq]\bfq^{\lan\rho,\lambda\ran}\chi_{\lambda}.
$$
Explicitly,  
$$
\Lambda^+=\bbN(1,0)\oplus\bbZ(1,1)\subset\bbZ(1,0)\oplus\bbZ(0,1)=\Lambda,
$$
so that
$$
\sS_{\cl}(\bfq)(\cH^{\sph}(\bfq))=\bigg(\bigoplus_{n\in\bbN}\bbZ[\bfq]\bfq^{\frac{n}{2}}\chi_{(n,0)}\bigg)\otimes_{\bbZ}\bbZ[\chi_{(1,1)}^{\pm1}].
$$
On the other hand, recall that the ring of symmetric polynomials in the two variables $e^{(1,0)}$ and $e^{(0,1)}$ is a graded ring generated the two characters $\chi_{(1,0)}=e^{(1,0)}+e^{(0,1)}$ and $\chi_{(1,1)}=e^{(1,0)}e^{(0,1)}$:
$$
\bbZ[e^{(1,0)},e^{(0,1)}]^{s}=\bigoplus_{n\in\bbN}\bbZ[e^{(1,0)},e^{(0,1)}]_n^s=\bbZ[\chi_{(1,0)},\chi_{(1,1)}].
$$
As $\chi_{(1,0)}$ is homogeneous of degree 1 and $\chi_{(1,1)}$ is homogeneous of degree 2, this implies that
$$
\bbZ[e^{(1,0)},e^{(0,1)}]_n^s=\bigoplus_{\substack{(a,b)\in\bbN^2\\a+2b=n}}\bbZ\chi_{(1,0)}^a\chi_{(1,1)}^b.
$$
Now if $a+2b=n$, then $\bfq^{\frac{n}{2}}\chi_{(1,0)}^a\chi_{(1,1)}^b=(\bfq^{\frac{1}{2}}\chi_{(1,0)})^a(\bfq\chi_{(1,1)})^b$. As the symmetric polynomial $\chi_{(n,0)}$ is homogeneous of degree $n$, we get the inclusion
$$
\sS_{\cl}(\bfq)(\cH^{\sph}(\bfq))\subset\bbZ[\bfq][\bfq^{\frac{1}{2}}\chi_{(1,0)},\bfq\chi_{(1,1)}]\otimes_{\bbZ}\bbZ[\chi_{(1,1)}^{\pm1}]
=\bbZ[\bfq][\bfq^{\frac{1}{2}}\chi_{(1,0)},\chi_{(1,1)}^{\pm1}].
$$
Since by definition of $\sS_{\cl}(\bfq)$ we have $\sS_{\cl}(\bfq)(f_{(1,0)}(\bfq))=\bfq^{\frac{1}{2}}\chi_{(1,0)}$, $\sS_{\cl}(\bfq)(f_{(1,1)}(\bfq))=\chi_{(1,1)}$ and $\sS_{\cl}(\bfq)(f_{(-1,-1)}(\bfq))=\chi_{(-1,-1)}=\chi_{(1,1)}^{-1}$, this inclusion is an equality. We have thus obtained the $\bbZ[\bfq]$-algebra isomorphism:
$$
\xymatrix{
\sS_{\cl}(\bfq)|_{\cH^{\sph}(\bfq)}:\cH^{\sph}(\bfq) \ar[r]^<<<<<{\sim} & \bbZ[\bfq][\bfq^{\frac{1}{2}}\chi_{(1,0)},\chi_{(1,1)}^{\pm1}].
}
$$
Also note that $T_{(1,0)}\mapsto \bfq^{\frac{1}{2}}\chi_{(1,0)}$ and $T_{(1,1)}\mapsto \chi_{(1,1)}$ since $T_{(1,0)}=f_{(1,0)}(\bfq)$ and $T_{(1,1)}=f_{(1,1)}(\bfq)$.

Finally, recall that $V_{\widehat{\bfT}}$ being the diagonalizable $k$-monoid scheme $\Spec(k[X,Y,z_2^{\pm1}])$, we have
$$
R(V_{\widehat{\bfT}})^{W_0}=\bbZ[X,Y,z_2^{\pm1}]^{W_0}=\bbZ[X+Y,XY,z_2^{\pm1}]=\bbZ[X+Y,XYz_2^{-1},z_2^{\pm1}].
$$
Hence we can define a ring isomorphism
$$
\xymatrix{
\iota:\bbZ[\bfq][\bfq^{\frac{1}{2}}\chi_{(1,0)},\chi_{(1,1)}^{\pm1}]\ar[r]^<<<<<{\sim} & R(V_{\widehat{\bfT}})^{W_0}
}
$$
by $\iota(\bfq):=XYz_2^{-1}$, $\iota(\bfq^{\frac{1}{2}}\chi_{(1,0)})=X+Y$ and $\iota(\chi_{(1,1)})=z_2$. Composing, we get the desired isomorphism
$$
\xymatrix{
\sS(\bfq):=\iota\circ\sS_{\cl}(\bfq)|_{\cH^{\sph}(\bfq)}:\cH^{\sph}(\bfq) \ar[r]^<<<<<{\sim} & R(V_{\widehat{\bfT}})^{W_0}.
}
$$
Note that $\sS(\bfq)(T_{(1,0)})=X+Y$, $\sS(\bfq)(T_{(1,1)})=z_2$, $\sS(\bfq)(\bfq)=XYz_2^{-1}$, and that $\sS(\bfq)$ is uniquely determined by these assignments since the ring $\cH^{\sph}(\bfq)$ is the polynomial ring in the variables $\bfq$, $T_{(1,0)}$ and $T_{(1,1)}^{\pm 1}$, thanks to the isomorphism $\sS_{\cl}(\bfq)|_{\cH^{\sph}(\bfq)}$.
\end{proof}

\begin{Rem*}\label{Remiota}
The choice of the isomorphism $\iota$ in the preceding proof may seem \emph{ad hoc}. However, it is natural from the point of view of the Vinberg fibration $\bfq:V_{\mathbf{\whT}}\ra \bbA^1$. 

First, as pointed out by Herzig in \cite[\S 1.2]{H11}, one can make the classical complex Satake transform $\sS_{\cl}$ integral, by removing the factor $\delta^{\frac{1}{2}}$ from its definition, where $\delta$ is the modulus character of the Borel subgroup. Doing so produces a ring embedding
$$
\xymatrix{
\cS':\cH_{\bbZ}^{\sph}\ar@{^{(}->}[r] & \bbZ[\bbX^{\bullet}(\mathbf{\whT})].
}
$$
The image of $\cS'$ is not contained in the subring $\bbZ[\bbX^{\bullet}(\mathbf{\whT})]^{W_0}$ of $W_0$-invariants. In fact, 
$$
\cS'(T_{(1,0)})=qe^{(1,0)}+e^{(0,1)}\quad\textrm{and}\quad\cS'(T_{(1,1)})=e^{(1,1)},
$$
so that 
$$
\xymatrix{
\cS':\cH_{\bbZ}^{\sph}\ar[r]^<<<<<{\sim} & \bbZ[(qe^{(1,0)}+e^{(0,1)}),e^{\pm (1,1)}]\subset \bbZ[\bbX^{\bullet}(\mathbf{\whT})].
}
$$
Now, 
$$
\bbZ[\bbX^{\bullet}(\mathbf{\whT})]=\bbZ[\mathbf{\whT}]=\bbZ[V_{\mathbf{\whT},1}],
$$
where $\mathbf{\whT}\cong V_{\mathbf{\whT},1}$ is the fiber \emph{at $1$} of the fibration $\bfq:V_{\mathbf{\whT}}\ra \bbA^1$ considered over 
$\bbZ$. But the algebra $\cH_{\bbZ}^{\sph}$ is the specialisation \emph{at $q$} of the generic algebra $\cH^{\sph}(\bfq)$. From this perspective, the morphism $\cS'$ is unnatural, since it mixes a $1$-fiber with a $q$-fiber. To restore the $\bfq$-compatibility, one must consider the composition of $\bbQ\otimes_{\bbZ}\cS'$ with the isomorphism 
\begin{eqnarray*}
\bbQ[V_{\mathbf{\whT},1}]=\bbQ[X,Y,z_2^{\pm1}]/(XY-z_2)&\xrightarrow{\sim} &\bbQ[V_{\mathbf{\whT},q}]=\bbQ[X,Y,z_2^{\pm1}]/(XY-qz_2)\\
X & \mapsto & q^{-1}X \\
Y & \mapsto & Y \\
z_2 & \mapsto & z_2.
\end{eqnarray*}
But then one obtains the formulas
\begin{eqnarray*}
\cH_{\bbQ}^{\sph}&\xrightarrow{\sim} &\bbQ[V_{\mathbf{\whT},q}]=\bbQ[X,Y,z_2^{\pm1}]/(XY-qz_2)\\
T_{(1,0)} & \mapsto & X+Y \\
T_{(1,1)} & \mapsto & z_2 \\
q & \mapsto & XYz_2^{\pm1}.
\end{eqnarray*}
This composed map is defined over $\bbZ$, its image is the subring $\bbZ[V_{\mathbf{\whT},q}]^{W_0}$ of $W_0$-invariants, and it is precisely the specialisation $\bfq=q$ of the isomorphism $\sS(\bfq)$ from \ref{ThgenSat}.
\end{Rem*}

\begin{Def*}\label{DefgenSat}
We call 
$$
\xymatrix{
\sS(\bfq):\cH^{\sph}(\bfq)\ar[r]^<<<<<{\sim} & R(V_{\mathbf{\whT}})^{W_0}
}
$$
the \emph{generic Satake isomorphism}. 
\end{Def*}

\begin{Pt*}
Composing with the inverse of the character isomorphism $\chi_{V_{\mathbf{\whG}}}^{-1}: R(V_{\mathbf{\whT}})^{W_0}\iso R(V_{\mathbf{\whG}})$
from \ref{charVin}, we arrive at an isomorphism
$$
\xymatrix{
\chi_{V_{\mathbf{\whG}}}^{-1}\circ\sS(\bfq):\cH^{\sph}(\bfq)\ar[r]^<<<<<{\sim} & R(V_{\mathbf{\whG}}).
}
$$
\end{Pt*}

\begin{Pt*}
Next, recall the generic Iwahori-Hecke algebra $\cH_1(\bfq)$ \ref{defgeneric1}, and the commutative subring $\cA_1(\bfq)\subset \cH_1(\bfq)$ \ref{presA1q} together with the isomorphism $\sB_1(\bfq)$ in \ref{SpecB1q}.
\end{Pt*}

\begin{Def*}\label{sZ1def}
The \emph{generic central elements morphism} is the unique ring homomrphism
$$
\xymatrix{
\sZ_1(\bfq):\cH^{\sph}(\bfq)\ar[r] & \cA_1(\bfq)\subset\cH_1(\bfq)
}
$$
making the diagram 
$$
\xymatrix{
 \cA_1(\bfq) \ar[rr]_{\sim}^{\sB_1(\bfq)}  && R(V_{\mathbf{\whT}}) \\
\cH^{\sph}(\bfq) \ar[u]^{\sZ_1(\bfq)} \ar[rr]_{\sim}^{\sS(\bfq)} && R(V_{\mathbf{\whT}})^{W_0} \ar@{^{(}->}[u]
}
$$
commutative.
\end{Def*}

\begin{Pt*}\label{sZ1iso}
By construction, the morphism $\sZ_1(\bfq)$ is injective, and is uniquely determined by the following equalities in $\cA_1(\bfq)$:
$$
\sZ_1(\bfq)(T_{(1,0)})=z_1,\quad \sZ_1(\bfq)(T_{(1,1)})=z_2\quad\textrm{and}\quad\sZ_1(\bfq)(\bfq)=\bfq.
$$
Moreover the group $W_0$ acts on the ring  $\cA_1(\bfq)$ and the invariant subring $\cA_1(\bfq)^{W_0}$ is equal to the center $Z(\cH_1(\bfq))\subset \cH_1(\bfq)$. As the isomorphism $\sB_1(\bfq)$ is $W_0$-equivariant by construction, we obtain that the image of $\sZ_1(\bfq)$ indeed is equal to the \emph{center of the generic Iwahori-Hecke algebra $\cH_1(\bfq)$}:
$$
\xymatrix{
\sZ_1(\bfq):\cH^{\sph}(\bfq)\ar[r]^>>>>>{\sim} & Z(\cH_1(\bfq))\subset \cA_1(\bfq)\subset \cH_1(\bfq).
}
$$
\end{Pt*}

\begin{Pt*}\label{compgenBSiso}
Under the identification $R(V_{\mathbf{\whT}})=\bbZ[V_{\mathbf{\whT}}]$ of \ref{SpecB1q}, the elements $\sS(\bfq)(T_{(1,0)})=X+Y$, 
$\sS(\bfq)(\bfq)=\bfq$, $\sS(\bfq)(T_{(1,1)})=z_2$,  correspond to the \emph{Steinberg choice of coordinates} $z_1$, $\bfq$, $z_2$ on the affine $\bbZ$-scheme $V_{\mathbf{\whT}}/W_0=\Spec(\bbZ[V_{\mathbf{\whT}}]^{W_0})$. On the other hand, the \emph{Trace of representations morphism}
$\Tr:R(V_{\mathbf{\whG}})\ra  \bbZ[V_{\mathbf{\whG}}]^{\mathbf{\whG}}$ fits into the commutative diagram 
$$
\xymatrix{
R(V_{\mathbf{\whT}})^{W_0} \ar@{=}[d] &\ar[l]^{\sim}_{\chi_{V_{\mathbf{\whG}}}} R(V_{\mathbf{\whG}}) \ar[d]^{\Tr} \\
\bbZ[V_{\mathbf{\whT}}]^{W_0} &\ar[l]^{\sim}_{\Ch}  \bbZ[V_{\mathbf{\whG}}]^{\mathbf{\whG}}
}
$$
where $\chi_{V_{\mathbf{\whG}}}$ is the character isomorphism of \ref{charVin}, and $\Ch$ is the \emph{Chevalley isomorphism} which is constructed for the Vinberg monoid $V_{\mathbf{\whG}}$ by Bouthier in \cite[Prop. 1.7]{Bo15}. So we have the following commutative diagram of $\bbZ$-schemes
$$
\xymatrix{
\Spec(\cA_1(\bfq)) \ar@{->>}[d]_{\Spec(\sZ_1(\bfq))} &  & V_{\mathbf{\whT}} \ar@{->>}[d]  \ar[ll]_>>>>>>>>>>>>>>>>>>>{\Spec(\sB_1(\bfq))}^>>>>>>>>>>>>>>>>>>>{\sim} \ar@{^{(}->}[rr] && V_{\mathbf{\whG}}\ar@{->>}[d] \\
\Spec(\cH^{\sph}(\bfq)) \ar[dr]^{\sim}_{(T_{(1,0)},\bfq,T_{(1,1)})} & & V_{\mathbf{\whT}}/W_0 \ar[dl]^{(z_1,\bfq,z_2)}_{\sim} \ar[ll]^>>>>>>>>>>>>>>>>>>>{\sim}_>>>>>>>>>>>>>>>>>>>{\Spec(\sS(\bfq))} \ar[rr]_{\sim}^{\Spec(\Ch)} && V_{\mathbf{\whG}}//\mathbf{\whG}\\
&\bbA^2\times\bbG_m. &
}
$$
Note that for $\mathbf{\whG}=\mathbf{GL_2}$, the  composed \emph{Chevalley-Steinberg map} $V_{\mathbf{\whG}}\ra \bbA^2\times\bbG_m$ is given explicitly by attaching to a $2\times 2$ matrix its characteristic polynomial (when $z_2=1$).
\end{Pt*}

\begin{Pt*}
We have recalled that for the generic pro-$p$-Iwahori-Hecke algebra $\cH^{(1)}(\bfq)$ too, the center can be described in terms of $W_0$-invariants, namely $Z(\cH^{(1)}(\bfq))=\cA^{(1)}(\bfq)^{W_0}$, cf. \ref{centerHI1}. As the generic Bernstein isomorphism 
$\sB^{(1)}(\bfq)$ is $W_0$-equivariant by construction, cf. \ref{genB1}, we can make the following definition. 
\end{Pt*}

\begin{Def*}\label{DefgenSat1}
We call 
$$
\xymatrix{
\sS^{(1)}(\bfq):=\sB^{(1)}(\bfq)^{W_0}:\cA^{(1)}(\bfq)^{W_0}\ar[r]^<<<<<{\sim} & R(V^{(1)}_{\mathbf{\whT}})^{W_0}
}
$$
the \emph{generic pro-$p$-Iwahori Satake isomorphism}.
\end{Def*}

\begin{Pt*}\label{SpecS1q}
Note that with $V^{\gamma}_{\mathbf{\whT}}:=\coprod_{\lambda\in\gamma} V_{\mathbf{\whT}}$ we have $V^{(1)}_{\mathbf{\whT}}=\bbT^{\vee}\times V_{\mathbf{\whT}}=\coprod_{\gamma\in \bbT^{\vee}/W_0}V^{\gamma}_{\mathbf{\whT}} $ and the $W_0$-action on this scheme respects these $\gamma$-components. We obtain the decomposition into connected components 

 $$V^{(1)}_{\mathbf{\whT}}/W_0 = \coprod_{\gamma\in \bbT^{\vee}/W_0} (\coprod_{\lambda\in\gamma} V_{\mathbf{\whT}})/W_0 =\coprod_{\gamma\in \bbT^{\vee}/W_0} V^{\gamma}_{\mathbf{\whT}}/W_0
 $$
 If $\gamma$ is regular, then $V^{\gamma}_{\mathbf{\whT}}/W_0\simeq V_{\mathbf{\whT}}$, the isomorphism depending on a choice of order on the set $\gamma$, cf. \ref{centergamma}. Hence, passing to $\tilde{\bbZ}$ as in \ref{SpecB1q}, with $\tilde{\cH}^{(1)}(\bfq):=\cH_{\tilde{\bbZ}}^{(1)}(\bfq)$, we obtain the following commutative diagram of $\tilde{\bbZ}$-schemes.

$$
\xymatrix{
\Spec(\tilde{\cA}^{(1)}(\bfq)) \ar@{->>}[d]  && V^{(1)}_{\mathbf{\whT}} \ar@{->>}[d]  \ar[ll]_{\Spec(\tilde{\sB}^{(1)}(\bfq))}^>>>>>>>>>>>>>>>>>>>{\sim}  \\
\Spec(Z(\tilde{\cH}^{(1)}(\bfq))) \ar[d]^{\rotatebox{90}{$\sim$}} && V^{(1)}_{\mathbf{\whT}}/W_0 \ar[ll]^{\sim}_{\Spec(\tilde{\sS}^{(1)}(\bfq))}  \ar[d]_{\rotatebox{90}{$\sim$}}^{\ref{centergamma}} \\
 (\bbA^2\times\bbG_m)^{\bbT^{\vee}/W_0}
&&\coprod_{(\bbT^{\vee}/W_0)_{\reg}}V_{\mathbf{\whT}}\coprod_{(\bbT^{\vee}/W_0)_{\nonreg}} V_{\mathbf{\whT}}/W_0, \ar[ll]^>>>>>>>>>>>{\sim}
}
$$
where the bottom isomorphism of the diagram is given by the standard coordinates $(x,y,z_2)$ on the regular components and by the Steinberg coordinates $(z_1,\bfq,z_2)$ on the non-regular components. 
\end{Pt*}

\subsection{The generic parametrization} \label{genericpara}

We keep the notation $\bbZ\subset \tilde{\bbZ}$ for the ring extension of \ref{finiteT}. Then we have defined the $\tilde{\bbZ}$-scheme 
$V^{(1)}_{\mathbf{\whT}}$ in \ref{SpecB1q}, and we have considered in \ref{SpecS1q} its quotient by the natural $W_0$-action.
Also recall that $\mathbf{\whG}=\mathbf{GL_2}$ is the Langlands dual $k$-group of $GL_{2,F}$.

\begin{Def*}
The category of quasi-coherent modules on the $\tilde{\bbZ}$-scheme $V^{(1)}_{\mathbf{\whT}}/W_0$ will be called the \emph{category of Satake parameters}, and denoted by $\SP_{\mathbf{\whG}}$:
$$
\SP_{\mathbf{\whG}}:=\QCoh\big(V^{(1)}_{\mathbf{\whT}}/W_0\big).
$$

For $\gamma \in \bbT^{\vee}/W_0$, we also define $\SP^{\gamma}_{\mathbf{\whG}}:=\QCoh\big(
V^{\gamma}_{\mathbf{\whT}}/W_0 \big)$, where as above $V^{\gamma}_{\mathbf{\whT}}=\coprod_{\lambda\in\gamma} V_{\mathbf{\whT}}.$
\end{Def*}

\begin{Pt*}
Now, over $\tilde{\bbZ}$, we have the isomorphism 
$$
\xymatrix{
i_{\tilde{\sS}^{(1)}(\bfq)}:=\Spec(\tilde{\sS}^{(1)}(\bfq)):V^{(1)}_{\mathbf{\whT}}/W_0\ar[r]^<<<<<{\sim} & \Spec(Z(\tilde{\cH}^{(1)}(\bfq)))
}
$$
from the scheme $V^{(1)}_{\mathbf{\whT}}/W_0$ to the spectrum of the center $Z(\tilde{\cH}^{(1)}(\bfq))$ of the generic pro-$p$-Iwahori Hecke algebra $\tilde{\cH}^{(1)}(\bfq)$, cf. \ref{SpecS1q}.  
\end{Pt*}

\begin{Cor*}
The category of modules over $Z(\tilde{\cH}^{(1)}(\bfq))$ is equivalent to the category of Satake parameters:
$$
\xymatrix{
S:=(i_{\tilde{\sS}^{(1)}(\bfq)})^*:\Mod(Z(\tilde{\cH}^{(1)}(\bfq))) \ar@<1ex>[r]_<<<<<{\sim} & \SP_{\mathbf{\whG}}:(i_{\tilde{\sS}^{(1)}(\bfq)})_*. \ar@<1ex>[l] 
}
$$ 
The equivalence $S$ will be referred to as the \emph{functor of Satake parameters}.\footnote{We hope that there is only little risk of confusing the notation $S$ with the Hecke operator introduced in \ref{distinguishedelements}.} The quasi-inverse $(i_{\tilde{\sS}^{(1)}(\bfq)})_*$ will be denoted by $S^{-1}$.
\end{Cor*}

\begin{Pt*}
Still from \ref{SpecS1q}, these categories decompose as products over $\bbT^{\vee}/W_0$ (considered as a finite set), compatibly with the equivalences: for all $\gamma\in \bbT^{\vee}/W_0$, 
$$
\xymatrix{
S^{\gamma}:=(i_{\tilde{\sS}^{\gamma}(\bfq)})^*:\Mod(Z(\tilde{\cH}^{\gamma}(\bfq))) \ar@<1ex>[r]_<<<<<{\sim} & \SP_{{\mathbf{\whG}}}^{\gamma}:(i_{\tilde{\sS}^{\gamma}(\bfq)})_*, \ar@<1ex>[l] 
}
$$ 
where
$$
 \SP_{{\mathbf{\whG}}}^{\gamma}\simeq
\left\{ \begin{array}{ll}
\QCoh(V_{\mathbf{\whT}}) & \textrm{ if $\gamma$ is regular} \\
\QCoh(V_{\mathbf{\whT}}/W_0)& \textrm{ if $\gamma$ is non-regular}.
\end{array} \right.
$$
In the regular case, the latter isomorphism depends on a choice of order on the set $\gamma$.
\end{Pt*}

\begin{Pt*}\label{Striv}
In particular, we have the trivial orbit $\gamma:=\{1\}$. The corresponding component $\tilde{\cH}^{\{1\}}(\bfq)$ of $\tilde{\cH}^{(1)}(\bfq)$ is canonically isomorphic to the $\tilde{\bbZ}$-base change of the  generic non-regular Iwahori-Hecke algebra $\cH_1(\bfq)$. Hence from \ref{sZ1iso} we have an isomorphism
$$
\xymatrix{
\tilde{\sZ}_1(\bfq):\tilde{\cH}^{\sph}(\bfq)\ar[r]^>>>>>{\sim} & Z(\tilde{\cH}^{\{1\}}(\bfq))\subset \tilde{\cA}^{\{1\}}(\bfq)\subset \tilde{\cH}^{\{1\}}(\bfq)\subset \tilde{\cH}^{(1)}(\bfq).
}
$$
Using these identifications, the equivalence $S^{\gamma}$ for $\gamma:=\{1\}$ can be rewritten as
$$
\xymatrix{
S^{\{1\}}:\Mod(\tilde{\cH}^{\sph}(\bfq)) \ar[r]^<<<<<{\sim} & \SP_{{\mathbf{\whG}}}^{\{1\}}.
}
$$ 
\end{Pt*}

\begin{Def*}
The category of quasi-coherent modules on the $\tilde{\bbZ}$-scheme $V^{(1)}_{\mathbf{\whT}}$ will be called the \emph{category of Bernstein parameters}, and denoted by $\BP_{\mathbf{\whG}}$:
$$
\BP_{\mathbf{\whG}}:=\QCoh\big(V^{(1)}_{\mathbf{\whT}}\big).
$$
\end{Def*}

\begin{Pt*}
Over $\tilde{\bbZ}$, we have the isomorphism 
$$
\xymatrix{
i_{\tilde{\sB}^{(1)}(\bfq)}:=\Spec(\tilde{\sB}^{(1)}(\bfq)):V^{(1)}_{\mathbf{\whT}}\ar[r]^<<<<<{\sim} & \Spec(\tilde{\cA}^{(1)}(\bfq)))
}
$$
from the scheme $V^{(1)}_{\mathbf{\whT}}$ to the spectrum of the commutative subring $\tilde{\cA}^{(1)}(\bfq)$ of the generic pro-$p$-Iwahori Hecke algebra $\tilde{\cH}^{(1)}(\bfq)$, cf. \ref{SpecB1q}.  Also we have the \emph{restriction functor}
$$
\xymatrix{
\Res^{\tilde{\cH}^{(1)}(\bfq)}_{\tilde{\cA}^{(1)}(\bfq)}:\Mod(\tilde{\cH}^{(1)}(\bfq)) \ar[r] & \Mod(\tilde{\cA}^{(1)}(\bfq))\cong\QCoh(\Spec(\tilde{\cA}^{(1)}(\bfq)))
}
$$
from the category of left $\tilde{\cH}^{(1)}(\bfq)$-modules to the one of  $\tilde{\cA}^{(1)}(\bfq)$-modules, equivalently of quasi-coherent modules on $\Spec(\tilde{\cA}^{(1)}(\bfq))$.
\end{Pt*}

\begin{Def*}
The \emph{functor of Bernstein parameters} is the composed functor
$$
\xymatrix{
B:=(i_{\tilde{\sB}^{(1)}(\bfq)})^*\circ\Res^{\tilde{\cH}^{(1)}(\bfq)}_{\tilde{\cA}^{(1)}(\bfq)}:\Mod(\tilde{\cH}^{(1)}(\bfq)) \ar[r] & \BP_{\mathbf{\whG}}.
}
$$ 
\end{Def*}

\begin{Pt*}
Still from \ref{SpecB1q}, the category $\BP_{{\mathbf{\whG}}}$ decomposes as a product over the finite group $\bbT^{\vee}$: 
$$
 \BP_{{\mathbf{\whG}}}\cong\prod_{\lambda\in\bbT^{\vee}} \BP_{{\mathbf{\whG}}}^{\lambda},
\quad\textrm{where}\quad
\forall \lambda\in\bbT^{\vee},\ \BP_{{\mathbf{\whG}}}^{\lambda}\simeq\QCoh(V_{\mathbf{\whT}}).
$$
\end{Pt*}

\begin{Pt*}\label{pi}
Denoting by $\pi: V^{(1)}_{\mathbf{\whT}}\rightarrow  V^{(1)}_{\mathbf{\whT}}/W_0$
the canonical projection, the compatibilty between the functors $S$ and $B$ of Satake and Bernstein parameters is expressed by the commutativity of the diagram
$$
\xymatrix{
\Mod(\tilde{\cH}^{(1)}(\bfq)) \ar[r]^<<<<<<B \ar[d]_{\Res_{Z(\tilde{\cH}^{(1)}(\bfq))}^{\tilde{\cH}^{(1)}(\bfq)}}  & \BP_{\mathbf{\whG}}\ar[d]^{\pi_*} \\
\Mod(Z(\tilde{\cH}^{(1)}(\bfq))) \ar[r]^<<<<<S_>>>>>{\sim} & \SP_{\mathbf{\whG}}. 
}
$$ 
\end{Pt*}

\begin{Def*}
The \emph{generic parametrization functor} is the functor
$$
P:=S\circ \Res_{Z(\tilde{\cH}^{(1)}(\bfq))}^{\tilde{\cH}^{(1)}(\bfq)} = \pi_*\circ B :
$$

$$
\xymatrix{
\Mod(\tilde{\cH}^{(1)}(\bfq))\ar[d] & \\
\SP_{\mathbf{\whG}}.
}
$$
\end{Def*}

\begin{Pt*}
It follows from the definitions that for all $\gamma\in \bbT^{\vee}/W_0$, the fiber of $P$ over the direct factor $\SP_{\mathbf{\whG}}^{\gamma}\subset \SP_{\mathbf{\whG}}$ is the direct factor $\Mod(\tilde{\cH}^{\gamma}(\bfq))\subset  \Mod(\tilde{\cH}^{(1)}(\bfq))$:
$$
P^{-1}(\SP_{\mathbf{\whG}}^{\gamma})=\Mod(\tilde{\cH}^{\gamma}(\bfq))\subset  \Mod(\tilde{\cH}^{(1)}(\bfq)).
$$
Accordingly the parametrization functor $P$ decomposes as the product over the finite set $\bbT^{\vee}/W_0$ of functors
$$
\xymatrix{
P^{\gamma}:\Mod(\tilde{\cH}^{\gamma}(\bfq))\ar[r] & \SP_{\mathbf{\whG}}^{\gamma}.
}
$$
\end{Pt*}

\begin{Pt*}\label{Strivv}
In the case of the trivial orbit $\gamma:=\{1\}$, it follows from \ref{Striv} that $P^{\{1\}}$ factors as
$$
\xymatrix{
 \Mod(\tilde{\cH}^{\{1\}}(\bfq)) \ar[d]_{\Res^{\tilde{\cH}^{\{1\}}(\bfq)}_{\tilde{\cH}^{\sph}(\bfq)}} \ar[dr]^{P^{\{1\}}}& \\
 \Mod(\tilde{\cH}^{\sph}(\bfq)) \ar[r]_<<<<<<{\sim}^>>>>{S^{\{1\}}} & \SP_{\mathbf{\whG}}^{\{1\}}.
}
$$
\end{Pt*}

\subsection{The generic spherical module}

Recall the generic regular and non-regular spherical representations $\sA_2(\bfq)$ \ref{sA2q} and $\sA_1(\bfq)$ \ref{sA1q} of $\cH_2(\bfq)$ and $\cH_1(\bfq)$. Thanks to \ref{H2VSHgamma} and \ref{H1VSHgamma}, they are models over $\bbZ$ of representations 
$\tilde{\sA}^{\gamma}(\bfq)$ of the regular and non-regular components $\tilde{\sA}^{\gamma}(\bfq)$, $\gamma\in\bbT^{\vee}/W_0$, of the generic pro-$p$-Iwahori Hecke algebra $\tilde{\cH}^{(1)}(\bfq)$ over $\tilde{\bbZ}$, cf. \ref{decompprop} and \ref{AIH}. Taking the product over 
$\bbT^{\vee}/W_0$ of these representations, we obtain a representation
$$
\xymatrix{
\tilde{\sA}^{(1)}(\bfq):\tilde{\cH}^{(1)}(\bfq) \ar[r] & \End_{Z(\tilde{\cH}^{(1)}(\bfq))}(\tilde{\cA}^{(1)}(\bfq)).
}
$$
By construction, the representation $\tilde{\sA}^{(1)}(\bfq)$ depends on a choice of order on each regular orbit $\gamma$.

\begin{Def*}
We call $\tilde{\sA}^{(1)}(\bfq)$ the \emph{generic spherical representation}, and the corresponding left $\tilde{\cH}^{(1)}(\bfq)$-module 
$\tilde{\cM}^{(1)}$ the \emph{generic spherical module}.
\end{Def*}


\begin{Prop*}\label{propantisph} \begin{enumerate}
\item The generic spherical representation is faithful. 
\item The Bernstein parameter of the spherical module is the structural sheaf:
$$
B(\cM^{(1)})=\cO_{V^{(1)}_{\mathbf{\whT}}}.
$$
\item The Satake parameter of the spherical module is the $\tilde{R}(V^{(1)}_{\mathbf{\whG}})$-module of $V^{(1)}_{\mathbf{\whG}}$-equivariant $K$-theory of the flag variety of $V^{(1)}_{\mathbf{\whG}}$:
$$
\tilde{c}_{V^{(1)}_{\mathbf{\whG}}}:S(\cM^{(1)})\stackrel{\sim}{\longrightarrow} \tilde{K}^{V^{(1)}_{\mathbf{\whG}}}(V^{(1)}_{\mathbf{\whG}}/V^{(1)}_{\mathbf{\whB}}).
$$
\end{enumerate}
\end{Prop*}

\begin{proof}
Part 1. follows from \ref{sA2sf} and \ref{sA1sf}, part 2. from the property \emph{(i)} in \ref{sA2q} and \ref{sA1q}, and part 3. from the characteristic isomorphism in \ref{B1qchar}.
\end{proof}

\begin{Pt*}
Now, being a left $\tilde{\cH}^{(1)}(\bfq)$-module, the spherical module  $\tilde{\cM}^{(1)}$ defines a functor
$$
\xymatrix{
 \tilde{\cM}^{(1)}\otimes_{Z(\tilde{\cH}^{(1)}(\bfq))}\bullet:\Mod(Z(\tilde{\cH}^{(1)}(\bfq)))\ar[r] & \Mod(\tilde{\cH}^{(1)}(\bfq)).
}
$$
On the other hand, recall the canonical projection 
$\pi:V^{(1)}_{\mathbf{\whT}}\ra V^{(1)}_{\mathbf{\whT}}/W_0$
from \ref{pi}.
Then point 2. of \ref{propantisph} has the following consequence.
\end{Pt*}

\begin{Cor*}
The diagram 
$$
\xymatrix{
\Mod(\tilde{\cH}^{(1)}(\bfq)) \ar[r]^<<<<<<{B} & \BP_{\mathbf{\whG}} \\
\Mod(Z(\tilde{\cH}^{(1)}(\bfq))) \ar[r]^<<<<{S}_>>>>>>{\sim} \ar[u]^{ \tilde{\cM}^{(1)}\otimes_{Z(\tilde{\cH}^{(1)}(\bfq))}\bullet} & \SP_{\mathbf{\whG}} \ar[u]_{\pi^*}
}
$$ 
is commutative.
\end{Cor*}

\begin{Def*}
The \emph{generic spherical functor} is the functor
$$
\Sph:=  (\tilde{\cM}^{(1)}\otimes_{Z(\tilde{\cH}^{(1)}(\bfq))}\bullet)\circ S^{-1}:
$$
$$
\xymatrix{
\SP_{\mathbf{\whG}}\ar[r] & \Mod(\tilde{\cH}^{(1)}(\bfq)).
}
$$

\end{Def*}

\begin{Cor*}
The diagram 
$$
\xymatrix{
&&\Mod(\tilde{\cH}^{(1)}(\bfq))\ar[d]^{P}  \\
\SP_{\mathbf{\whG}} \ar[urr]^{\Sph} \ar[r]_{\pi^*}  & \BP_{\mathbf{\whG}}\ar[r]_{\pi_*}& \SP_{\mathbf{\whG}}
}
$$
is commutative.
\end{Cor*}
\begin{proof} One has $P\circ \Sph = \pi_* \circ (B \circ \Sph ) =  \pi_* \circ \pi^*$ by the preceding corollary.
\end{proof}

\begin{Pt*}
By construction, the spherical functor $\Sph$ decomposes as a product of functors $\Sph^{\gamma}$ for $\gamma\in \bbT^{\vee}/W_0$, and accordingly the previous diagram decomposes over $\bbT^{\vee}/W_0$. 
\end{Pt*}

\begin{Pt*}
In particular for $\gamma=\{1\}$ we have the commutative diagram
$$
\xymatrix{
&&\Mod(\tilde{\cH}^{\{1\}}(\bfq))\ar[d]^{P^{\{1\}}}  \ar[dr]^{\Res^{\tilde{\cH}^{\{1\}}(\bfq)}_{\tilde{\cH}^{\sph}(\bfq)}} \\
\SP_{\mathbf{\whG}}^{\{1\}} \ar[urr]^{\Sph^{\{1\}}} \ar[r]_{\pi^*}  & \BP_{\mathbf{\whG}}^{\{1\}}\ar[r]_{\pi_*}& \SP_{\mathbf{\whG}}^{\{1\}}
& \Mod(\tilde{\cH}^{\sph}(\bfq)). \ar[l]^>>>>>>>>{S^{\{1\}}}_>>>>>>>>{\sim}
}
$$
\end{Pt*}

\section{The theory at $\bfq=q=0$}

We keep all the notations introduced in the preceding section. In particular, $k=\overline{\bbF}_q$.

\subsection{$K$-theory of the dual flag variety at $\bfq=0$}

\begin{Pt*}
Recall from \ref{Vinsubsection} the $k$-semigroup scheme
$$
V_{\mathbf{GL_2},0}=\Sing_ {2\times 2}\times\bbG_m,
$$
which can even be defined over $\bbZ$, and which is obtained as the $0$-fiber of 
$$
\xymatrix{
V_{\mathbf{GL_2}}\ar[d]^{\bfq} \\
\bbA^1.
}
$$
\end{Pt*}

\begin{Pt*}
It admits 
$$
V_{\mathbf{\whT},0}=\SingDiag_ {2\times 2}\times\bbG_m
$$
as a commutative subsemigroup scheme. The latter has the following structure: it is the pinching of the monoids
$$
\bbA_X^1\times\bbG_m:=\Spec(k[X,z_2^{\pm1}]) \quad \textrm{and} \quad \bbA_Y^1\times\bbG_m:=\Spec(k[Y,z_2^{\pm1}]) 
$$ 
along the sections $X=0$ and $Y=0$. These monoids are semisimple, with representation rings
$$
R(\bbA_X^1\times\bbG_m)=\bbZ[X,z_2^{\pm1}] \quad \textrm{and} \quad R(\bbA_Y^1\times\bbG_m)=\bbZ[Y,z_2^{\pm1}]. 
$$

There are three remarkable elements in $V_{\mathbf{\whT},0}$, namely
$$
\varepsilon_X:=(\diag(1,0),1),\quad \varepsilon_Y:=(\diag(0,1),1)\quad\textrm{and}\quad\varepsilon_0:=(\diag(0,0),1).
$$
They are idempotents. Now let $M$ be a finite dimensional $k$-representation of $V_{\mathbf{\whT},0}$. The idempotents act on $M$ as projectors, and as the semigroup $V_{\mathbf{\whT},0}$ is commutative, the $k$-vector space $M$ decomposes as a direct sum 
$$
M=\bigoplus_{(\lambda_X,\lambda_Y,\lambda_0)\in\{0,1\}^3} M(\lambda_X,\lambda_Y,\lambda_0)
$$
where 
$$
M(\lambda_X,\lambda_Y,\lambda_0)=\{m\in M\ |\ m\varepsilon_X=\lambda_Xm,\ m\varepsilon_Y=\lambda_Ym,\ m\varepsilon_0=\lambda_0m\}.
$$
Moreover, since $V_{\mathbf{\whT},0}$ is commutative, each of these subspaces is in fact a subrepresentation of $M$.

As $\varepsilon_X\varepsilon_Y=\varepsilon_0\in V_{\mathbf{\whT},0}$, we have $M(\lambda_X,\lambda_Y,\lambda_0)\neq 0\implies
\lambda_X\lambda_Y=\lambda_0$. 
Consequently
$$
M=M(1,0,0)\bigoplus M(0,1,0)\bigoplus M(1,1,1)\bigoplus M(0,0,0).
$$
The restriction $\Res^{V_{\mathbf{\whT},0}}_{\bbA_X^1}M(1,0,0)$ is a representation of the monoid $\bbA_X^1$ where $0$ acts by $0$, and 
 $\Res^{V_{\mathbf{\whT},0}}_{\bbA_Y^1}M(1,0,0)$ is the null representation. Hence, if for $n>0$ we still denote by $X^n$ the character of $V_{\mathbf{\whT},0}$ which restricts to the character $X^n$ of $\bbA_X^1\times\bbG_m$ and the null map of $\bbA_Y^1\times\bbG_m$, then $M(1,0,0)$ decomposes as a sum of weight spaces
$$
M(1,0,0)=\oplus_{n>0}M(X^n):=\oplus_{n>0,m\in\bbZ}M(X^nz_2^m).
$$
Similarly
$$
M(0,1,0)=\oplus_{n>0}M(Y^n):=\oplus_{n>0,m\in\bbZ}M(Y^nz_2^m).
$$
Finally, $V_{\mathbf{\whT},0}$ acts through the projection $V_{\mathbf{\whT},0}\ra\bbG_m$ on
$$
M(1,1,1)=:M(1)=\oplus_{m\in\bbZ}M(z_2^m),
$$
and by $0$ on 
$$
M(0,0,0)=:M(0).
$$
Thus we have obtained the following
\end{Pt*}

\begin{Lem*}\label{RVT0}
The category $\Rep(V_{\mathbf{\whT},0})$ is semisimple, and there is a ring isomorphism
$$
R(V_{\mathbf{\whT},0})\cong\big(\bbZ[X,Y,z_2^{\pm1}]/(XY)\big)\times\bbZ.
$$
\end{Lem*}

\begin{Pt*}
Next let  
$$
V_{\mathbf{\whB},0}=\SingUpTriang_ {2\times 2}\times\bbG_m\subset V_{\mathbf{GL_2},0}=\Sing_{2\times 2}\times \bbG_m
$$
be the subsemigroup scheme of singular upper triangular $2\times 2$-matrices. It contains $V_{\mathbf{\whT},0}$, and the inclusion $V_{\mathbf{\whT},0}\subset V_{\mathbf{\whB},0}$ admits a retraction $V_{\mathbf{\whB},0}\ra V_{\mathbf{\whT},0}$, namely the specialisation at $\bfq=0$ of the retraction \ref{retraction}.

Let $M$ be an object of $\Rep(V_{\mathbf{\whB},0})$. Write 
$$
\Res^{V_{\mathbf{\whB},0}}_{V_{\mathbf{\whT},0}}M=M(1,0,0)\oplus M(0,1,0)\oplus M(1)\oplus M(0).
$$

\noindent For a subspace $N\subset M$, consider the following property:

\medskip

(P${}_N$) \emph{the subspace $N\subset M$ is a subrepresentation, and $V_{\mathbf{\whB},0}$ acts on $N$ through the retraction of $k$-semigroup schemes $V_{\mathbf{\whB},0}\ra V_{\mathbf{\whT},0}$}. 

\medskip

\noindent Let us show that (P${}_{M(0,1,0)}$) is true. Indeed for $m\in M(0,1,0)=\oplus_{n>0}M(Y^n)$, we have
$$
m\left (\begin{array}{cc}
x & c\\
0 & 0
\end{array} \right)
=
(m\varepsilon_Y)\left (\begin{array}{cc}
x & c\\
0 & 0
\end{array} \right)
=m\varepsilon_0=0
=m\left (\begin{array}{cc}
x & 0\\
0 & 0
\end{array} \right)
$$
and
$$
m\left (\begin{array}{cc}
0 & c\\
0 & y
\end{array} \right)
=
(m\varepsilon_Y)\left (\begin{array}{cc}
0 & c\\
0 & y
\end{array} \right)
=m\left (\begin{array}{cc}
0 & 0\\
0 & y
\end{array} \right).
$$

\noindent Next assume $M(0,1,0)=0$, and let us show that in this case (P${}_{M(0)}$) is true. Indeed for $m\in M(0)$, we have
$$
m\left (\begin{array}{cc}
x & c\\
0 & 0
\end{array} \right)
=
m\bigg(\varepsilon_X\left (\begin{array}{cc}
x & c\\
0 & 0
\end{array} \right)\bigg)
=(m\varepsilon_X)\left (\begin{array}{cc}
x & c\\
0 & 0
\end{array} \right)
=0,
$$
and if we decompose 
$$
m':=m\left (\begin{array}{cc}
0 & c\\
0 & y
\end{array} \right)
=m_{(1,0,0)}'+m_1'+m_0'\in M(1,0,0)\oplus M(1)\oplus M(0),
$$
then by applying $\varepsilon_X$ on the right we see that $0=m_{(1,0,0)}'+m_1'$ so that $m'\in M(0)$ and hence
$$
m\left (\begin{array}{cc}
0 & c\\
0 & y
\end{array} \right)
=
m\bigg(\left (\begin{array}{cc}
0 & c\\
0 & y
\end{array} \right)\varepsilon_Y\bigg)
=m'\varepsilon_Y=0.
$$

\noindent Next assume $M(0,1,0)=M(0)=0$, and let us show that in this case (P${}_{M(1,0,0)}$) is true. Indeed, let $m\in M(1,0,0)=\oplus_{n>0}M(X^n)$. Then for any $c\in k$,
$$
m':=m\left (\begin{array}{cc}
0 & c\\
0 & 0
\end{array} \right)
$$
satisfies $m'\varepsilon_X=0$, $m'\varepsilon_Y=m'$, $m'\varepsilon_0=0$, i.e. $m'\in M(0,1,0)$, and hence is equal to $0$ by our assumption. It follows that
$$
m\left (\begin{array}{cc}
0 & c\\
0 & y
\end{array} \right)
=
(m\varepsilon_X)\left (\begin{array}{cc}
0 & c\\
0 & y
\end{array} \right)
=
m\bigg(\varepsilon_X\left (\begin{array}{cc}
0 & c\\
0 & y
\end{array} \right)\bigg)
=m\left (\begin{array}{cc}
0 & c\\
0 & 0
\end{array} \right)
=0
=m\left (\begin{array}{cc}
0 & 0\\
0 & y
\end{array} \right).
$$
On the other hand, if we decompose
$$
m':=m
\left (\begin{array}{cc}
x & c\\
0 & 0
\end{array} \right)
=m_{(1,0,0)}'+m_1'\in M(1,0,0)\oplus M(1),
$$
then by applying $\varepsilon_0$ on the right we find $0=m_1'$, i.e. $m'\in M(1,0,0)$ and hence
$$
m
\left (\begin{array}{cc}
x & c\\
0 & 0
\end{array} \right)
=m'=m'\varepsilon_X=m\bigg(\left (\begin{array}{cc}
x & c\\
0 & 0
\end{array} \right)\varepsilon_X\bigg)
=
m\left (\begin{array}{cc}
x & 0\\
0 & 0
\end{array} \right).
$$
Finally assume $M(0,1,0)=M(0)=M(1,0,0)=0$, and let us show that in this case (P${}_{M(1)}$) is true, i.e. that $V_{\mathbf{\whB},0}$ acts through the projection $V_{\mathbf{\whB},0}\ra\bbG_m$ on $M=M(1)$. Indeed for any $m$ we have
$$
m
\left (\begin{array}{cc}
x & c\\
0 & 0
\end{array} \right)
=
\bigg(m
\left (\begin{array}{cc}
x & c\\
0 & 0
\end{array} \right)\bigg)\varepsilon_0
=
m\bigg(
\left (\begin{array}{cc}
x & c\\
0 & 0
\end{array} \right)\varepsilon_0\bigg)
=m\varepsilon_0=m
$$
and
$$
m
\left (\begin{array}{cc}
0 & c\\
0 & y
\end{array} \right)
=
\bigg(m
\left (\begin{array}{cc}
0 & c\\
0 & y
\end{array} \right)\bigg)\varepsilon_0
=
m\bigg(
\left (\begin{array}{cc}
0 & c\\
0 & y
\end{array} \right)\varepsilon_0\bigg)
=m\varepsilon_0=m.
$$

It follows from the preceding discussion that the irreducible representations of $V_{\mathbf{\whB},0}$ are the characters, which are inflated from those of $V_{\mathbf{\whT},0}$ through the retraction $V_{\mathbf{\whB},0}\ra V_{\mathbf{\whT},0}$. As a consequence, considering the \emph{restriction} and \emph{inflation} functors
$$
\xymatrix{
\Res_{V_{\mathbf{\whT},0}}^{V_{\mathbf{\whB},0}}:\Rep(V_{\mathbf{\whB},0}) \ar@<1ex>[r] & \Rep(V_{\mathbf{\whT},0}) :\Infl_{V_{\mathbf{\whT},0}}^{V_{\mathbf{\whB},0}}, \ar@<1ex>[l] 
}
$$
which are exact and compatible with tensor products and units, we get:
\end{Pt*}

\begin{Lem*}\label{ResInfl0}
The ring homomorphisms
$$
\xymatrix{
\Res_{V_{\mathbf{\whT},0}}^{V_{\mathbf{\whB},0}}:R(V_{\mathbf{\whB},0}) \ar@<1ex>[r] & R(V_{\mathbf{\whT},0}) :\Infl_{V_{\mathbf{\whT},0}}^{V_{\mathbf{\whB},0}}, \ar@<1ex>[l] 
}
$$
are isomorphisms, which are inverse one to the other.
\end{Lem*}

\begin{Pt*}
Finally, note that $\varepsilon_0\in V_{\mathbf{GL_2}}(k)$ belongs to all the left $V_{\mathbf{GL_2}}(k)$-cosets in $V_{\mathbf{GL_2}}(k)$. Hence, by \cite[2.4.3]{PS20}, the catgory $\Rep(V_{\mathbf{\whB},0})$ is equivalent to the one of induced vector bundles on the semigroupoid flag variety $V_{\mathbf{GL_2},0}/V_{\mathbf{\whB},0}$:
$$
\xymatrix{
\cI nd_{V_{\mathbf{\whB},0}}^{V_{\mathbf{GL_2},0}}:\Rep(V_{\mathbf{\whB},0}) \ar[r]^<<<<<{\sim} & \cC_{\cI nd}^{V_{\mathbf{GL_2},0}}(V_{\mathbf{GL_2},0}/V_{\mathbf{\whB},0}) \subset \cC^{V_{\mathbf{GL_2},0}}(V_{\mathbf{GL_2},0}/V_{\mathbf{\whB},0}).
}
$$
\end{Pt*}

\begin{Cor*}
We have a ring isomorphism
$$
\xymatrix{
\cI nd_{V_{\mathbf{\whB},0}}^{V_{\mathbf{GL_2},0}}\circ \Infl_{V_{\mathbf{\whT},0}}^{V_{\mathbf{\whB},0}} : R(V_{\mathbf{\whT},0}) \ar[r]^<<<<<{\sim} & K_{\cI nd}^{V_{\mathbf{GL_2},0}}(V_{\mathbf{GL_2},0}/V_{\mathbf{\whB},0}).
}
$$
\end{Cor*}

\begin{Def*}\label{defrel}
We call \emph{relevant} the full subcategory  
$$
\Rep(V_{\mathbf{\whT},0})^{\rel}\subset  \Rep(V_{\mathbf{\whT},0})
$$
 whose objects $M$ satisfy $M(0)=0$. Correspondingly, we have relevant full subcategories 
 $$
 \Rep(V_{\mathbf{\whB},0})^{\rel}\subset  \Rep(V_{\mathbf{\whB},0})\quad\textrm{and}\quad \cC_{\cI nd}^{V_{\mathbf{GL_2},0}}(V_{\mathbf{GL_2},0}/V_{\mathbf{\whB},0})^{\rel}\subset\cC_{\cI nd}^{V_{\mathbf{GL_2},0}}(V_{\mathbf{GL_2},0}/V_{\mathbf{\whB},0}).
 $$ 
\end{Def*}

\begin{Cor*}
We have a ring isomorphism
$$
\xymatrix{
c_{V_{\mathbf{GL_2},0}}:=\bbZ[X,Y,z_2^{\pm1}]/(XY)\cong R(V_{\mathbf{\whT},0})^{\rel} \ar[r]^<<<<<{\sim} & K_{\cI nd}^{V_{\mathbf{GL_2},0}}(V_{\mathbf{GL_2},0}/V_{\mathbf{\whB},0})^{\rel},
}
$$
that we call \emph{the characteristic isomorphism in the equivariant $K$-theory of the flag variety $V_{\mathbf{GL_2},0}/V_{\mathbf{\whB},0}$}.
\end{Cor*}

\begin{Pt*}
We have a commutative diagram \emph{specialization at $\bfq=0$}
$$
\xymatrix{
\bbZ[X,Y,z_2^{\pm1}] \ar[rr]^{c_{V_{\mathbf{GL_2}}}}_{\sim} \ar@{->>}[d] && K^{V_{\mathbf{GL_2}}}(V_{\mathbf{GL_2}}/V_{\mathbf{\whB}}) \ar@{->>}[d]\\
\bbZ[X,Y,z_2^{\pm1}]/(XY)\ar[rr]^{c_{V_{\mathbf{GL_2},0}}}_{\sim}  && K_{\cI nd}^{V_{\mathbf{GL_2},0}}(V_{\mathbf{GL_2},0}/V_{\mathbf{\whB},0})^{\rel},
}
$$
where the vertical right-hand side map is given by restricting equivariant vector bundles to the $0$-fiber of $\bfq: V_{\mathbf{GL_2}}\ra\bbA^1$. 
\end{Pt*}

\subsection{The mod $p$ Satake and Bernstein isomorphisms}

\begin{Not*}
In the sequel, we will denote by $(\bullet)_{\overline{\bbF}_q}$ the \emph{specialization at $\bfq=q=0$}, i.e. the base change functor along the ring morphism
\begin{eqnarray*}
\bbZ[\bfq] & \lra & \overline{\bbF}_q=:k\\
\bfq & \lmapsto & 0.
\end{eqnarray*}
Also we fix an embedding $\mu_{q-1}\subset\overline{\bbF}_q^{\times}$, so that the above morphism factors through the inclusion $\bbZ[\bfq] \subset \tilde{\bbZ}[\bfq]$, where $\bbZ\subset\tilde{\bbZ}$ is the ring extension considered in \ref{finiteT}.
\end{Not*} 

\begin{Pt*}
\textbf{The mod $p$ Satake and pro-$p$-Iwahori Satake isomorphisms.} Specializing \ref{DefgenSat}, we get an isomorphism of $\overline{\bbF}_q$-algebras
$$
\xymatrix{
\sS_{\overline{\bbF}_q}:\cH_{\overline{\bbF}_q}^{\sph}\ar[r]^<<<<<{\sim} & \overline{\bbF}_q[V_{\mathbf{\whT},0}]^{W_0}=\big(\overline{\bbF}_q[X,Y,z_2^{\pm1}]/(XY)\big)^{W_0}.
}
$$
In \cite{H11}, Herzig constructed an isomorphism
$$
\xymatrix{
\sS_{\Her}:\cH_{\overline{\bbF}_q}^{\sph}\ar[r]^<<<<<{\sim} & \overline{\bbF}_q[\bbX^{\bullet}(\mathbf{\whT})_-]=\overline{\bbF}_q[e^{(0,1)},e^{\pm(1,1)}]
}
$$
(this is $\overline{\bbF}_q\otimes_{\bbZ}\cS'$, with the notation $\cS'$ from \ref{Remiota}). They are related by the Steinberg choice of coordinates $z_1:=X+Y$ and $z_2$ on the quotient $V_{\mathbf{\whT},0}/W_0$, cf. \ref{compgenBSiso}, i.e. by the following commutative diagram
$$
\xymatrix{
\cH_{\overline{\bbF}_q}^{\sph}\ar[rr]_<<<<<<<<<<<<<<<<<<<<<<{\sim}^<<<<<<<<<<<<<<<<<<<<<<{\sS_{\overline{\bbF}_q}} \ar[dr]_{\sS_{\Her}}^{\sim} && \big(\overline{\bbF}_q[X,Y,z_2^{\pm1}]/(XY)\big)^{W_0} \\
&\overline{\bbF}_q[e^{(0,1)},e^{\pm(1,1)}] \ar[ur]_{\quad \quad e^{(0,1)}\mapsto z_1,\ e^{(1,1)}\mapsto z_2}^{\sim}. &
}
$$

Specializing \ref{DefgenSat1} and using $R(\bbT^{\vee})=\bbZ[\bbT]$, cf. \ref{Repring_of_bbT_dual}, we get an isomorphism of $\overline{\bbF}_q$-algebras
$$
\xymatrix{
\sS_{\overline{\bbF}_q}^{(1)}:(\cA_{\overline{\bbF}_q}^{(1)})^{W_0}\ar[r]^<<<<<{\sim} & \overline{\bbF}_q[V^{(1)}_{\mathbf{\whT},0}]^{W_0}=\big(\overline{\bbF}_q[\bbT][X,Y,z_2^{\pm1}]/(XY)\big)^{W_0}.
}
$$
\end{Pt*}

\begin{Pt*} \label{B1qcharmodp}
\textbf{The mod $p$ Bernstein isomorphism.} Specializing \ref{genB1}, we get an isomorphism of $\overline{\bbF}_q$-algebras
$$
\xymatrix{
\sB_{\overline{\bbF}_q}^{(1)}:\cA_{\overline{\bbF}_q}^{(1)}\ar[r]^<<<<<{\sim} & \overline{\bbF}_q[V^{(1)}_{\mathbf{\whT},0}]=\overline{\bbF}_q[\bbT][X,Y,z_2^{\pm1}]/(XY).
}
$$
Moreover, similarly as in \ref{B1qchar} but here using \ref{ResInfl0} and \cite[2.5.1]{PS20}, we get the \emph{characteristic isomorphism}
$$
\xymatrix{
c_{V^{(1)}_{\mathbf{\whG},0}}:R(V^{(1)}_{\mathbf{\whT},0}) \ar[r]^<<<<<{\sim} & K_{\cI nd}^{V^{(1)}_{\mathbf{\whG},0}}(V^{(1)}_{\mathbf{\whG},0}/V^{(1)}_{\mathbf{\whB},0}).
} 
$$
Whence by \ref{RVT0} (and recalling \ref{defrel}) an isomorphism
$$
\xymatrix{
c_{V^{(1)}_{\mathbf{\whG},0},\overline{\bbF}_q}^{\rel}\circ \sB_{\overline{\bbF}_q}^{(1)}:\cA_{\overline{\bbF}_q}^{(1)} \ar[r]^<<<<<{\sim} & K_{\cI nd,\overline{\bbF}_q}^{V^{(1)}_{\mathbf{\whG},0}}(V^{(1)}_{\mathbf{\whG},0}/V^{(1)}_{\mathbf{\whB},0})^{\rel}.
} 
$$

Also, specializing \ref{SpecB1q}, $\sB_{\overline{\bbF}_q}^{(1)}$ splits as a product over $\bbT^{\vee}$ of $\overline{\bbF}_q$-algebras isomorphisms $\sB_{\overline{\bbF}_q}^{\lambda}$, each of them being of the form
$$
\xymatrix{
\sB_{1,\overline{\bbF}_q}:\cA_{1,\overline{\bbF}_q}\ar[r]^<<<<<{\sim} & \overline{\bbF}_q[V_{\mathbf{\whT},0}]=\overline{\bbF}_q[X,Y,z_2^{\pm1}]/(XY).
}
$$
\end{Pt*}

\begin{Pt*}\label{sZ1isomodp}
\textbf{The mod $p$ central elements embedding.} Specializing \ref{sZ1def}, we get an embedding of $\overline{\bbF}_q$-algebras
$$
\xymatrix{
\sZ_{1,\overline{\bbF}_q}:\cH_{\overline{\bbF}_q}^{\sph}\ar[r]^<<<<<{\sim} & Z(\cH_{1,\overline{\bbF}_q}) \subset \cA_{1,\overline{\bbF}_q} \subset \cH_{1,\overline{\bbF}_q}
}
$$
making the diagram
$$
\xymatrix{
 \cA_{1,\overline{\bbF}_q} \ar[rr]_<<<<<<<<<<<{\sim}^<<<<<<<<<<<{\sB_{1,\overline{\bbF}_q}}  && \overline{\bbF}_q[V_{\mathbf{\whT},0}]=\overline{\bbF}_q[X,Y,z_2^{\pm1}]/(XY) \\
\cH^{\sph}_{\overline{\bbF}_q} \ar@{^{(}->}[u]^{\sZ_{1,\overline{\bbF}_q}} \ar[rr]_<<<<<<<<<{\sim}^<<<<<<<<<{\sS_{\overline{\bbF}_q} } && \overline{\bbF}_q[V_{\mathbf{\whT},0}]^{W_0}=\big(\overline{\bbF}_q[X,Y,z_2^{\pm1}]/(XY)\big)^{W_0}\ar@{^{(}->}[u]
}
$$
commutative. Then $\sZ_{1,\overline{\bbF}_q}$ coincides with the central elements construction of Ollivier \cite[Th. 4.3]{O14} for the case of 
$\mathbf{GL_2}$.
This follows from the explicit formulas for the values of $\sZ_{1}(\bfq)$ on $T_{ (1,0) }$ and 
$T_{(1,1) }$, cf. \ref{sZ1iso}.

\end{Pt*}

\subsection{The mod $p$ parametrization}

\begin{Def*}
The category of quasi-coherent modules on the $k$-scheme $V^{(1)}_{\mathbf{\whT},0}/W_0$ will be called the \emph{category of mod $p$ Satake parameters}, and denoted by $\SP_{\mathbf{\whG},0}$:
$$
\SP_{\mathbf{\whG},0}:=\QCoh\big(V^{(1)}_{\mathbf{\whT},0}/W_0\big).
$$

For $\gamma \in \bbT^{\vee}/W_0$, we also define $\SP^{\gamma}_{\mathbf{\whG},0}:=
\QCoh\big(V^{\gamma}_{\mathbf{\whT},0}/W_0 \big)$, where $ V^{\gamma}_{\mathbf{\whT},0}=\coprod_{\lambda\in\gamma} V_{\mathbf{\whT},0}$.
\end{Def*}

\begin{Pt*} 
Similarly to the generic case \ref{genericpara}, the mod $p$ pro-$p$-Iwahori Satake isomorphism induces an equivalence of categories
$$
\xymatrix{
S:\Mod(Z(\cH_{\overline{\bbF}_q}^{(1)})) \ar[r]^<<<<<{\sim} & \SP_{\mathbf{\whG},0},
}
$$
that will be referred to as the \emph{functor of mod $p$ Satake parameters}, and which decomposes as a product over the finite set $\bbT^{\vee}/W_0$:
$$
\xymatrix{
S=\prod_{\gamma}S^{\gamma}:\prod_{\gamma}\Mod(Z(\cH_{\overline{\bbF}_q}^{\gamma}))\ar[r]^<<<<<{\sim} & \prod_{\gamma}\SP_{\mathbf{\whG},0}^{\gamma}\simeq \prod_{\gamma\ \reg}\QCoh(V_{\mathbf{\whT},0})\prod_{\gamma\ \nonreg}\QCoh(V_{\mathbf{\whT},0}/W_0).
}
$$

For $\gamma=\{1\}$ and using \ref{sZ1isomodp} we get an equivalence
$$
\xymatrix{
S^{\{1\}}:\Mod(\cH_{\overline{\bbF}_q}^{\sph}) \ar[r]^<<<<<{\sim} & \SP_{\mathbf{\whG},0}^{\{1\}}=\QCoh(V_{\mathbf{\whT},0}/W_0).
}
$$
Note that under this equivalence, the characters $\cH_{\overline{\bbF}_q}^{\sph}\ra \overline{\bbF}_q$ correspond to the skyscraper sheaves on $V_{\mathbf{\whT},0}/W_0$, and hence to its $k$-points. Choosing the Steinberg coordinates $(z_1,z_2)$ on the $k$-scheme $V_{\mathbf{\whT},0}/W_0$, they may also be regarded as the $k$-points of $\Spec(k[\bbX^{\bullet}(\mathbf{\whT})_-])$, which are precisely the mod $p$ Satake parameters defined by Herzig in \cite{H11}.
\end{Pt*}

\begin{Def*}
The category of quasi-coherent modules on the $k$-scheme $V^{(1)}_{\mathbf{\whT},0}$ will be called the \emph{category of mod $p$ Bernstein parameters}, and denoted by $\BP_{\mathbf{\whG},0}$:
$$
\BP_{\mathbf{\whG},0}:=\QCoh\big(V^{(1)}_{\mathbf{\whT},0}\big).
$$
\end{Def*}

\begin{Pt*}
Similarly to the generic case \ref{genericpara}, the inclusion $\cH^{(1)}_{\overline{\bbF}_q}\supset \cA^{(1)}_{\overline{\bbF}_q}$ together with the mod $p$ Bernstein isomorphism define a  \emph{functor of mod $p$ Bernstein parameters} 
$$
\xymatrix{
B:\Mod(\cH^{(1)}_{\overline{\bbF}_q}) \ar[r] & \BP_{\mathbf{\whG},0}.
}
$$ 
Moreover the category $\BP_{\mathbf{\whG},0}$ decomposes as a product over the finite group $\bbT^{\vee}$:
$$
 \BP_{{\mathbf{\whG}},0}=\prod_{\lambda} \BP_{{\mathbf{\whG}},0}^{\lambda}=\prod_{\lambda}\QCoh(V_{\mathbf{\whT},0}).
$$
\end{Pt*}

\begin{Not*}
Let $\pi:V^{(1)}_{\mathbf{\whT},0}\ra V^{(1)}_{\mathbf{\whT},0}/W_0$ be the canonical projection.
\end{Not*}

\begin{Def*}\label{modpP}
The \emph{mod $p$ parametrization functor} is the functor
$$
P:=S\circ \Res_{Z(\cH^{(1)}_{\overline{\bbF}_q})}^{\cH^{(1)}_{\overline{\bbF}_q}} = \pi_*\circ B:
$$
$$
\xymatrix{
\Mod(\cH^{(1)}_{\overline{\bbF}_q})\ar[d] & \\
\SP_{\mathbf{\whG},0}.
}
$$
\end{Def*}

\begin{Pt*}
The functor $P$ decomposes as a product over the finite set $\bbT^{\vee}/W_0$:
$$
\xymatrix{
P=\prod_{\gamma}P^{\gamma}:\prod_{\gamma}\Mod(\cH_{\overline{\bbF}_q}^{\gamma})\ar[r]^<<<<<{\sim} & \prod_{\gamma}\SP_{\mathbf{\whG},0}^{\gamma}.
}
$$
In the case of the trivial orbit $\gamma:=\{1\}$, $P^{\{1\}}$ factors as
$$
\xymatrix{
 \Mod(\cH^{\{1\}}_{\overline{\bbF}_q}) \ar[d]_{\Res^{\cH^{\{1\}}_{\overline{\bbF}_q}}_{\cH^{\sph}_{\overline{\bbF}_q}}} \ar[dr]^{P^{\{1\}}}& \\
 \Mod(\cH^{\sph}_{\overline{\bbF}_q}) \ar[r]_<<<<<{\sim}^>>>>>{S^{\{1\}}} & \SP_{\mathbf{\whG},0}^{\{1\}}.
}
$$
\end{Pt*}

\subsection{The mod $p$ spherical module}

\begin{Def*}
We call 
$$
\xymatrix{
\sA_{\overline{\bbF}_q}^{(1)}:\cH^{(1)}_{\overline{\bbF}_q} \ar[r] & \End_{Z(\cH^{(1)}_{\overline{\bbF}_q})}(\cA^{(1)}_{\overline{\bbF}_q})
}
$$
the \emph{mod $p$ spherical representation}, and the corresponding left $\cH^{(1)}_{\overline{\bbF}_q}$-module $\cM^{(1)}_{\overline{\bbF}_q}$ the \emph{mod $p$ spherical module}.
\end{Def*}

\begin{Prop*}\label{propantisphmodp}
\begin{enumerate}
\item The mod $p$ spherical representation is faithful. 
\item The mod $p$ Bernstein parameter of the spherical module is the structural sheaf:
$$
B(\cM^{(1)}_{\overline{\bbF}_q})=\cO_{V^{(1)}_{\mathbf{\whT},0}}.
$$
\item The mod $p$ Satake parameter of the spherical module is the $R_{\overline{\bbF}_q}(V^{(1)}_{\mathbf{\whT},0})^{\rel,W_0}$-module of the relevant induced $V^{(1)}_{\mathbf{\whG},0}$-equivariant $K_{\overline{\bbF}_q}$-theory of the flag variety of $V^{(1)}_{\mathbf{\whG},0}$:
$$
c_{V^{(1)}_{\mathbf{\whG},0},\overline{\bbF}_q}^{\rel}:S(\cM^{(1)}_{\overline{\bbF}_q})\xrightarrow{\sim} K^{V^{(1)}_{\mathbf{\whG},0}}_{\cI nd,\overline{\bbF}_q}(V^{(1)}_{\mathbf{\whG},0}/V^{(1)}_{\mathbf{\whB},0})^{\rel}.
$$
\end{enumerate}
\end{Prop*}

\begin{proof}
Part 1. follows from \ref{sA2sf} and \ref{faithfulatzero}, 
part 2. from the property \emph{(i)} in \ref{sA2q} and \ref{sA1q}, and part 3. from the characteristic isomorphism in \ref{B1qcharmodp}.
\end{proof}

\begin{Cor*}
The diagram 
$$
\xymatrix{
\Mod(\cH^{(1)}_{\overline{\bbF}_q}) \ar[r]^<<<<<<{B} & \BP_{\mathbf{\whG},0} \\
\Mod(Z(\cH^{(1)}_{\overline{\bbF}_q})) \ar[r]^<<<<{S}_>>>>>>{\sim} \ar[u]^{ \cM^{(1)}_{\overline{\bbF}_q}\otimes_{Z(\cH^{(1)}_{\overline{\bbF}_q})}\bullet} & \SP_{\mathbf{\whG},0} \ar[u]_{\pi^*}
}
$$ 
is commutative.
\end{Cor*}

\begin{Def*}
The \emph{mod $p$ spherical functor} is the functor
$$
\Sph:=  (\cM^{(1)}_{\overline{\bbF}_q}\otimes_{Z(\cH^{(1)}_{\overline{\bbF}_q})}\bullet)\circ S^{-1}:
$$
$$
\xymatrix{
\SP_{\mathbf{\whG},0}\ar[r] & \Mod(\cH^{(1)}_{\overline{\bbF}_q}).
}
$$
\end{Def*}

\begin{Cor*}\label{PSph}
The diagram 
$$
\xymatrix{
&&\Mod(\cH^{(1)}_{\overline{\bbF}_q})\ar[d]^{P}  \\
\SP_{\mathbf{\whG},0} \ar[urr]^{\Sph} \ar[r]_{\pi^*}  & \BP_{\mathbf{\whG},0}\ar[r]_{\pi_*}& \SP_{\mathbf{\whG},0}
}
$$
is commutative.
\end{Cor*}

\begin{Pt*}
The spherical functor $\Sph$ decomposes as a product of functors $\Sph^{\gamma}$ for $\gamma\in \bbT^{\vee}/W_0$, and accordingly the previous diagram decomposes over $\bbT^{\vee}/W_0$.  In particular for $\gamma=\{1\}$ we have the commutative diagram
$$
\xymatrix{
&&\Mod(\cH^{\{1\}}_{\overline{\bbF}_q})\ar[d]^{P^{\{1\}}}  \ar[dr]^{\Res^{\cH^{\{1\}}_{\overline{\bbF}_q}}_{\cH^{\sph}_{\overline{\bbF}_q}}} \\
\SP_{\mathbf{\whG},0}^{\{1\}} \ar[urr]^{\Sph^{\{1\}}} \ar[r]_{\pi^*}  & \BP_{\mathbf{\whG},0}^{\{1\}}\ar[r]_{\pi_*}& \SP_{\mathbf{\whG},0}^{\{1\}}
& \Mod(\cH^{\sph}_{\overline{\bbF}_q}). \ar[l]^>>>>>>>>{S^{\{1\}}}_>>>>>>>>{\sim}
}
$$
\end{Pt*}

\begin{Pt*}\label{connected_comp}
Now, identifying the $k$-points of the $k$-scheme $V^{(1)}_{\mathbf{\whT},0}/W_0$ with the skyscraper sheaves on it, the spherical functor $\Sph$ induces a map
$$
\xymatrix{
\Sph:\big(V^{(1)}_{\mathbf{\whT},0}/W_0\big)(k)\ar[r] & \{\textrm{left $\cH^{(1)}_{\overline{\bbF}_q}$-modules}\}.
}
$$
Considering the decomposition of $V^{(1)}_{\mathbf{\whT},0}/W_0$ into its connected components, cf. \ref{SpecS1q},
$$
V^{(1)}_{\mathbf{\whT},0}/W_0=\coprod_{\gamma\in (\bbT^{\vee}/W_0)} V^{\gamma}_{\mathbf{\whT},0}/W_0\simeq\coprod_{\gamma\in (\bbT^{\vee}/W_0)_{\reg}}V_{\mathbf{\whT},0}\coprod_{\gamma\in (\bbT^{\vee}/W_0)_{\nonreg}}V_{\mathbf{\whT},0}/W_0,
$$
the spherical map decomposes as a disjoint union of maps
$$
\xymatrix{
\Sph^{\gamma}:\big(V^{\gamma}_{\mathbf{\whT},0}/W_0\big)(k)\simeq V_{\mathbf{\whT},0}(k)\ar[r] & \{\textrm{left
$\cH^{\gamma}_{\overline{\bbF}_q}$-modules}\}& \textrm{for $\gamma$ regular},
}
$$
$$
\xymatrix{
\Sph^{\gamma}:\big(V^{\gamma}_{\mathbf{\whT},0}/W_0\big)(k)\simeq(V_{\mathbf{\whT},0}/W_0)(k)\ar[r] & \{\textrm{left 
$\cH^{\gamma}_{\overline{\bbF}_q}$-modules}\}& \textrm{for $\gamma$ non-regular}.
}
$$
\end{Pt*}

\begin{Pt*} \label{regularcase}
In the regular case, we make the standard choice of coordinates
$$
V_{\mathbf{\whT},0}(k)=\bigg(\{(x,0)\ |\ x\in k\}\coprod_{(0,0)}\{(0,y)\ |\ y\in k\}\bigg)\times \{z_2\in k^{\times}\}
$$
and we identify $\cH^{\gamma}_{\overline{\bbF}_q}$ with $\cH_{2,\overline{\bbF}_q}$ using \ref{H2VSHgamma}. A point $v\in V_{\mathbf{\whT},0}(k)$ corresponds by \ref{presA2q} to a character
$$
\theta_v:Z(\cH_{2,\overline{\bbF}_q})\simeq\overline{\bbF}_q[X,Y,z_2^{\pm1}]/(XY)\lra\overline{\bbF}_q,
$$ 
and then $\Sph^{\gamma}(v)$ identifies with the central reduction
$$
\cA_{2,\theta_v}:=\cA_{2,\overline{\bbF}_q}\otimes_{Z(\cH_{2,\overline{\bbF}_q}),\theta_v}\overline{\bbF}_q
$$
of the mod $p$ regular spherical representation $\sA_{2,\overline{\bbF}_q}$ specializing \ref{sA2q}. The latter being an isomorphism by \ref{sA2sf}, so is
$$
\xymatrix{
\sA_{2,\theta_v}:\cH_{2,\theta_v} \ar[r]^<<<<{\sim} & \End_{\overline{\bbF}_q}(\cA_{2,\theta_v}).
}
$$
Consequently $\cH_{2,\theta_v}$ is a matrix algebra and $\cA_{2,\theta_v}$ is the unique simple finite dimensional left
$\cH_{2, \overline{\bbF_q}}$-module with central character $\theta_v$, up to isomorphism. It is the \emph{standard module with character 
$\theta_v$}, with \emph{standard basis} $\{\varepsilon_1,\varepsilon_2\}$ (in particular its $\overline{\bbF}_q$-dimension is $2$).
Conversely, any simple finite dimensional $\cH_{2, \overline{\bbF}_q}$-module has a central character, by Schur's lemma.

Following \cite{V04}, a central character $\theta$ is called {\it supersingular} if $\theta(X+Y)=0$, and the standard module with character 
$\theta$ is called {\it supersingular} if $\theta$ is. Since $XY=0$, one has $\theta(X+Y)=0$ if and only if $\theta(X)=\theta(Y)=0$. 
\end{Pt*}

\begin{Th*}\label{Sphreg}
Let $\gamma\in \bbT^{\vee}/W_0$ regular. Then the spherical map induces a bijection
$$
\xymatrix{
\Sph^{\gamma}:\big(V^{\gamma}_{\mathbf{\whT},0}/W_0\big)(k)\ar[r]^>>>>>{\sim} & \{\textrm{simple finite dimensional
 left $\cH^{\gamma}_{\overline{\bbF}_q}$-modules}\}/\sim.
}
$$

The singular locus of the parametrizing $k$-scheme
$V^{\gamma}_{\mathbf{\whT},0}/W_0
$
is given by $(0,0)\times\bbG_m\subset V_{\mathbf{\whT},0}$ in the standard coordinates, and its $k$-points correspond to the supersingular Hecke modules through the correspondence $\Sph^{\gamma}$.
\end{Th*}

\begin{Pt*}\label{non-regularcase}
In the non-regular case, we make the Steinberg choice of coordinates
$$
(V_{\mathbf{\whT},0}/W_0)(k)=\{z_1\in k\}\times \{z_2\in k^{\times}\}
$$
and we identify $\cH^{\gamma}_{\overline{\bbF}_q}$ with $\cH_{1,\overline{\bbF}_q}$ using \ref{H1VSHgamma}. A point 
$v\in (V_{\mathbf{\whT},0}/W_0)(k)$ corresponds to a character
$$
\theta_v:Z(\cH_{1,\overline{\bbF}_q})\simeq\overline{\bbF}_q[z_1,z_2^{\pm1}]\lra\overline{\bbF}_q,
$$ 
and then $\Sph^{\gamma}(v)$ identifies with the central reduction
$$
\cA_{1,\theta_v}:=\cA_{1,\overline{\bbF}_q}\otimes_{Z(\cH_{1,\overline{\bbF}_q}),\theta_v}\overline{\bbF}_q
$$
of the mod $p$ non-regular spherical representation $\sA_{1,\overline{\bbF}_q}$ specializing \ref{sA1q}. 

Now recall from \cite[1.4]{V04} the classification of the simple finite dimensional $\cH_{1,\overline{\bbF}_q}$-modules: they are the characters and the simple standard modules. The characters  
$$ 
\cH_{1,\overline{\bbF}_q}= \overline{\bbF}_q [S,U^{\pm 1}]\lra \overline{\bbF}^\times_q
$$
are parametrized by the set $\{0,-1\}\times \overline{\bbF}^\times_q$ via evaluation on the elements $S$ and $U$. On the other hand, given
$v=(z_1,z_2)\in k\times k^{\times}=\overline{\bbF}_q \times  \overline{\bbF}_q^{\times}$, a {\it standard module with character 
$\theta_v$} over $\cH_{1,\overline{\bbF}_q}$ is defined to be a module of type
$$ M_2(z_1,z_2):=\overline{\bbF}_q m \oplus \overline{\bbF}_q Um,\hskip15pt Sm=-m, \hskip15pt SUm=z_1 m, \hskip15pt U^2m=z_2 m $$
(in particular its $\overline{\bbF}_q$-dimension is $2$). The center $Z(\cH_{1,\overline{\bbF}_q})$ acts on $M_2(z_1,z_2)$ by the character $\theta_v$. In particular such a module is uniquely determined by its central character. It is simple if and only if $z_2\neq z_1^2$. It is called {\it supersingular} if $z_1=0$.
\end{Pt*}

\begin{Lem*} Set 
$$
\xymatrix{
\sA_{1,\theta_v}:=\sA_{1,\overline{\bbF}_q}\otimes_{Z(\cH_{1,\overline{\bbF}_q}),\theta_v}\overline{\bbF}_q:
\cH_{1,\theta_v}\ar[r] & \End_{\overline{\bbF}_q}(\cA_{1,\theta_v}).
}
$$
\begin{itemize}
\item Assume $z_2\neq z_1^2$. Then $\sA_{1,\theta_v}$ is an isomorphism, and the $\cH_{1,\overline{\bbF}_q}$-module $\cA_{1,\theta_v}$ is isomorphic to the simple standard module $M_2(z_1,z_2)$.
\item Assume $z_2=z_1^2$. Then $\sA_{1,\theta_v}$ has a $1$-dimensional kernel, and the $\cH_{1,\overline{\bbF}_q}$-module 
$\cA_{1,\theta_v}$ is a non-split extension of the character $(0,z_1)$ by the character $(-1,-z_1)$.
\end{itemize}
\end{Lem*}

\begin{proof}
The proof of Proposition \ref{sAqinj} shows that $\cH_{1,\theta_v}$ has an $\overline{\bbF}_q$-basis given by the elements $1, S, U, SU$, and that their images
$$
1,\ \sA_{1,\theta_v} (S),\ \sA_{1,\theta_v} (U),\ \sA_{1,\theta_v}(S)\sA_{1,\theta_v}(U)
$$
by $\sA_{1,\theta_v}$ are linearly independent over $\overline{\bbF}_q$ if and only if $z_1^2-z_2\neq 0$. 

If $z_2\neq z_1^2$, then $\sA_{1,\theta_v}$ is injective, and hence bijective since $\dim_{\overline{\bbF}_q} \cA_{1,\theta_v}=2$ from \ref{1Ybasis}. Moreover $S\cdot Y=-Y$ and $U\cdot Y=(z_1^2-z_2)-z_1Y$ and so $SUY=S ((z_1^2-z_2 ) -z_1Y) = S(-z_1Y)=z_1Y$, so that
$$
\cA_{1,\theta_v}=\overline{\bbF}_qY\oplus \overline{\bbF}_qU\cdot Y=M_2(z_1,z_2).
$$
If $z_2=z_1^2$, then the proof of Proposition \ref{sAqinj} shows that $\sA_{1,\theta_v}$ has a $1$-dimensional kernel which is the 
$\overline{\bbF}_q$-line generated by $-z_1(1+S)+U+SU$. Moreover $\overline{\bbF}_qY\subset\cA_{1,\theta_v}$ realizes the character 
$(-1,-z_1)$ of $\cH_{1,\overline{\bbF}_q}$, and $\cA_{1,\theta_v}/\overline{\bbF}_qY\simeq\overline{\bbF}_q1$ realizes the character $(0,z_1)$. Finally the $0$-eigenspace of $S$ in $\cA_{1,\theta_v}$ is $\overline{\bbF}_q1$, which is not $U$-stable, so that the character $(0,z_1)$ does not lift in $\cA_{1,\theta_v}$.
\end{proof}

\begin{Rem*}
Geometrically, the function $z_2-z_1^2$ on $V_{\mathbf{\whT},0}/W_0$ defines a family of parabolas
$$
\xymatrix{
V_{\mathbf{\whT},0}/W_0\ar[d]^{z_2-z_1^2},\\
\bbA^1
}
$$
whose parameter is $4\Delta$, where $\Delta$ is the discriminant of the parabola. Then the locus of $V_{\mathbf{\whT},0}/W_0$ where $z_2=z_1^2$ corresponds to the parabola at $0$, having vanishing discriminant (at least if $p\neq 2$).
\end{Rem*}

\begin{Def*}
We will say that a pair of characters of $\cH_{1,\overline{\bbF}_q}=\overline{\bbF}_q [S,U^{\pm 1}]\ra \overline{\bbF}_q^{\times}$ is \emph{antispherical} if there exists $z_1\in \overline{\bbF}^\times_q$ such that, after evaluating on $(S,U)$, it is equal to
$$
\{(0,z_1),(-1,-z_1)\}.
$$
\end{Def*}

\begin{Pt*}
Note that the set of characters $\cH_{1,\overline{\bbF}_q}\ra \overline{\bbF}_q^{\times}$ is the disjoint union of the spherical pairs, by the very definition.
\end{Pt*}

\begin{Th*}\label{Sphnonreg}
Let $\gamma\in \bbT^{\vee}/W_0$ non-regular. Consider the decomposition 
$$
V^{\gamma}_{\mathbf{\whT},0}/W_0=D(2)_{\gamma}\cup D(1)_{\gamma}
$$
where $D(1)_{\gamma}$ is the closed subscheme defined by the parabola $z_2=z_1^2$ in the Steinberg coordinates $z_1,z_2$ and $D(2)_{\gamma}$ is the open complement. Then the spherical map induces bijections
$$
\xymatrix{
\Sph^{\gamma}(2):D(2)_{\gamma}(k)\ar[r]^>>>>>{\sim} & \{\textrm{simple $2$-dimensional
 left $\cH^{\gamma}_{\overline{\bbF}_q}$-modules}\}/\sim
}
$$
$$
\xymatrix{
\Sph^{\gamma}(1):D(1)_{\gamma}(k)\ar[r]^>>>>>{\sim} & \{\textrm{spherical pairs of characters of $\cH^{\gamma}_{\overline{\bbF}_q}$}\}/\sim.
}
$$

The branch locus of the covering 
$$
V_{\mathbf{\whT},0}\lra V_{\mathbf{\whT},0}/W_0\simeq V^{\gamma}_{\mathbf{\whT},0}/W_0
$$ 
is contained in $D(2)_{\gamma}$, with equation $z_1=0$ in Steinberg coordinates, and its $k$-points correspond to the supersingular Hecke modules through the correspondence $\Sph^{\gamma}(2)$.
\end{Th*}

\begin{Rem*}
The matrices of $S$, $U$ and $S_0=USU^{-1}$ in the $\overline{\bbF}_q$-basis $\{1,Y\}$ of the supersingular module $\cA_{1,\theta_v}\cong M_2(0,z_2)$ are
$$
S=
\left (\begin{array}{cc}
0 & 0\\
0 & -1
\end{array} \right),\hskip15pt
U=
\left (\begin{array}{cc}
0 & -z_2\\
-1 & 0
\end{array} \right),
\hskip15pt
S_0=
\left (\begin{array}{cc}
-1 & 0\\
0 & 0
\end{array} \right).
$$
The two characters of the \emph{finite subalgebra} $\overline{\bbF}_q[S]$ corresponding to $S\mapsto 0$ and $S\mapsto -1$
are realized by $1$ and $Y$. From the matrix of $S_0$, we see in fact that the whole \emph{affine subalgebra} $\overline{\bbF}_q[S_0,S]$ acts on $1$ and $Y$ via the two \emph{supersingular affine characters}, which by definition are the characters different from the trivial character $(S_0,S)\mapsto(0,0)$ and the sign character $(S_0,S)\mapsto(-1,-1)$.
\end{Rem*}

\begin{Pt*}
Finally, let $v$ be any $k$-point of the parametrizing space $V^{(1)}_{\mathbf{\whT},0}/W_0$. As a particular case of \ref{PSph}, the Bernstein parameter of the spherical module $\Sph(v)$ is the structure sheaf of the fiber of the quotient map $\pi$ at $v$, and its Satake parameter is the underlying $k$-vector space:
$$
B(\Sph(v))=\cO_{\pi^{-1}(v)}\quad\textrm{and}\quad S(\Sph(v))=\pi_*\cO_{\pi^{-1}(v)}.
$$
\end{Pt*}

\subsection{Central characters}\label{Sat_par_with_cc}

In this final subsection, we show that the dual parametrization \ref{Sphnonreg} behaves naturally with respect to central characters. 

\begin{Pt*} \label{Fq_action} Let $\omega : \bbF_q^{\times}\ra k^\times$ be induced by the inclusion $\bbF_q \subset k$. 
Then $(\bbF_q^{\times})^{\vee}=\langle \omega \rangle$ is a cyclic group of order $q-1$. 
An element $\omega^r$ defines a non-regular character of $\bbT$: 
$$\omega^r(t_1,t_2):=\omega^{r}(t_1)\omega^{r}(t_2)$$ 
for all $(t_1,t_2)\in\bbT=\bbF_q^{\times}\times\bbF_q^{\times}$. Composing with multiplication in $\bbT^{\vee}$, we get an action of 
$(\bbF_q^{\times})^{\vee}$ on  $\bbT^{\vee}$,
which factors on the quotient set $\bbT^{\vee}/W_0$:
$$  \bbT^{\vee}/W_0 \times (\bbF_q^{\times})^{\vee} \longrightarrow  \bbT^{\vee}/W_0, \;(\gamma, \omega^r) \mapsto \gamma\omega^r.$$
If $\gamma\in \bbT^{\vee}/W_0$ is regular (non-regular), then $\gamma\omega^r$
is regular (non-regular). 
\end{Pt*}

\begin{Pt*} \label{orbits}
Restricting characters of $\bbT$ to the subgroup $\bbF_q^{\times}\simeq \{ \diag(a,a) : a\in \bbF_q^{\times}\} $ induces a homomorphism $\bbT^{\vee}\ra   (\bbF_q^{\times})^{\vee}$ which factors into a restriction map  
$$ \bbT^{\vee}/W_0 \ra   (\bbF_q^{\times})^{\vee}, \; \gamma\mapsto \gamma |_{\bbF_q^{\times}}.$$
The relation to the $(\bbF_q^{\times})^{\vee} $-action on the source $\bbT^{\vee}/W_0$ is given by the formula
$$(\gamma\omega^r)|_{ \bbF_q^{\times}} = \gamma|_{ \bbF_q^{\times}} \;\omega^{2r}.$$ 

We describe the fibers of the restriction map $\gamma\mapsto \gamma |_{\bbF_q^{\times}}$.

Let $(\cdot)|_{ \bbF_q^{\times}}^{-1}(\omega^{2r})$ be the fibre at a square element
$\omega^{2r}$. By the above formula, the action of $\omega^{-r}$ on $\bbT^{\vee}/W_0$ induces a bijection with the fibre $(\cdot)|_{\bbF_q^{\times}}^{-1}(1).$
The fibre 
$$(\cdot)|_{\bbF_q^{\times}}^{-1}(1)=\{ 1 \otimes 1 \} \coprod \{ \omega \otimes \omega^{-1}, \omega^2 \otimes \omega^{-2},...,  
\omega^{\frac{q-3}{2}} \otimes \omega^{-{\frac{q-3}{2}}} \} \coprod \{ \omega^{\frac{q-1}{2}} \otimes \omega^{-{\frac{q-1}{2}}} \}$$
has cardinality $\frac{q+1}{2}$ and, in the above list, we have chosen a representative in $\bbT^{\vee}$ for each element in the fibre. 
The $\frac{q-3}{2}$ elements in the middle of this list, i.e. the $W_0$-orbits 
represented by the characters $\omega^r \otimes \omega^{-r}$ for $r=1,...,\frac{q-3}{2}$, are all regular $W_0$-orbits. The two orbits at the two ends of the list
are non-regular orbits (note that $\frac{q-1}{2}\equiv -\frac{q-1}{2} \mod (q-1)$). Since the action of $\omega^{-r}$ preserves regular (non-regular) orbits, any fibre at a square element (there are 
$\frac{q-1}{2}$ such fibres) has the same structure.

On the other hand, let $(\cdot)|_{\bbF_q^{\times}}^{-1}(\omega^{2r-1})$ be the fibre at a non-square element $\omega^{2r-1}$. The action of 
$\omega^{-r}$ induces a bijection 
with the fibre $(\cdot)|_{\bbF_q^{\times}}^{-1}(\omega^{-1})$.  
The fibre 
$$(\cdot)|_{\bbF_q^{\times}}^{-1}(\omega^{-1})= \{ 1 \otimes \omega^{-1}, \omega \otimes \omega^{-2},...,  
\omega^{\frac{q-1}{2}-1} \otimes \omega^{-{\frac{q-1}{2}}} \} $$
has cardinality $\frac{q-1}{2}$ and we have chosen a representative in $\bbT^{\vee}$ for each element in the fibre. 
All elements of the fibre are regular $W_0$-orbits. Since the action of $\omega^{-r}$ preserves regular (non-regular) orbits, any fibre at a non-square element (there are $\frac{q-1}{2}$ such fibres) has the same structure.

Note that $\frac{q-1}{2}( \frac{q+1}{2} +  \frac{q-1}{2} ) = \frac{q^2-q}{2}$ is the cardinality of the set  $\bbT^{\vee}/W_0$.
\end{Pt*}

\begin{Pt*}\label{ZGveeaction}
Recall the commutative $k$-semigroup scheme
$$
V^{(1)}_{\mathbf{\whT},0}=\bbT^{\vee}\times V_{\mathbf{\whT},0}=\bbT^{\vee}\times  \SingDiag_ {2\times 2}\times\bbG_m 
$$ 
together with its $W_0$-action, cf. \ref{SpecS1q}: the natural action of $W_0$ on 
the factors $\bbT^{\vee}$ and  $\SingDiag_ {2\times 2}$ and the trivial one on $\bbG_m$. 
There is a commuting action of the $k$-group scheme
$$\cZ^\vee:= (\bbF_q^{\times})^{\vee}\times \bbG_m$$ 
on $V^{(1)}_{\mathbf{\whT},0}$:  
the (constant finite diagonalizable) group $(\bbF_q^{\times})^{\vee}$ acts only on the factor $\bbT^{\vee}$ and in the way described in \ref{Fq_action}; an element $z_0\in\bbG_m$ acts trivially on $\bbT^{\vee}$, by multiplication with the diagonal matrix $\diag(z_0,z_0)$ on $\SingDiag_ {2\times 2}$ and by multiplication with the square $z_0^2$ on $\bbG_m$. Therefore the quotient 
$V^{(1)}_{\mathbf{\whT},0}/W_0$ inherits a $\cZ^\vee$-action. Now, according to \ref{connected_comp}, one has the decomposition
$$
V^{(1)}_{\mathbf{\whT},0}/W_0=\coprod_{\gamma\in(\bbT^{\vee}/W_0)_{\reg}} V_{\mathbf{\whT},0}\coprod_{\gamma\in(\bbT^{\vee}/W_0)_{\nonreg} } V_{\mathbf{\whT},0}/W_0.
$$
Then the $(\bbF_q^{\times})^{\vee}$-action is by permutations on the index set $\bbT^{\vee}/W_0$, i.e. on the 
set of connected components of $V^{(1)}_{\mathbf{\whT},0}/W_0$; as observed above, it preserves the subsets of regular and non-regular components. The $\bbG_m$-action on $V^{(1)}_{\mathbf{\whT},0}/W_0$ preserves each connected component.

\end{Pt*} 

\begin{Pt*} \label{projections}
The two canonical projections from $V^{(1)}_{\mathbf{\whT},0}$ to $\bbT^{\vee}$ and 
$\bbG_m$ respectively induce two projection morphisms 

$$
\xymatrix{
&V^{(1)}_{\mathbf{\whT},0}/W_0 \ar[dl]_{\pr_{\bbT^{\vee}/W_0}} \ar[dr]^{\pr_{\bbG_m}} &  \\
\bbT^{\vee}/W_0  & &\bbG_m .
}
$$
Then we may compose the map $\pr_{\bbT^{\vee}/W_0}$ with the restriction map $(\cdot) |_{\bbF^{\times}_q} :\bbT^{\vee}/W_0\ra(\bbF_q^{\times})^{\vee}$, set
$$\theta := \big((\cdot) |_{\bbF^{\times}_q} \circ \pr_{\bbT^{\vee}/W_0}\big) \times {\pr_{\bbG_m}} $$
and view $V^{(1)}_{\mathbf{\whT},0}/W_0$ as fibered over the space $\cZ^\vee$:
$$
\xymatrix{
&V^{(1)}_{\mathbf{\whT},0}/W_0 \ar[d]^\theta &  \\
&\cZ^\vee. &
}
$$
The relation to the $\cZ^\vee$-action on the source $V^{(1)}_{\mathbf{\whT},0}/W_0$ is given by the formula
$$\theta (x.(\omega^r,z_0) )=\theta(x)(\omega^{2r},z^2_0)=\theta(x)(\omega^{r},z_0)^2$$ 
for $x\in V^{(1)}_{\mathbf{\whT},0}/W_0$ and $(\omega^r,z_0)\in\cZ^{\vee}$. This formula 
follows from the formula in \ref{orbits} and the definition of the $\bbG_m$-action in \ref{ZGveeaction}.
\end{Pt*}

\begin{Def*} \label{def_space_Sat_cc}
Let $\zeta\in \cZ^\vee$. The \emph{space of mod $p$ Satake parameters with central character $\zeta$} 
is the $k$-scheme $$(V^{(1)}_{\mathbf{\whT},0}/W_0)_{\zeta}:= \theta^{-1}(\zeta).$$
\end{Def*}

\begin{Pt*}\label{connectedcomp}
Let $\zeta=(\zeta |_ {\bbF^{\times}_q}, z_2)\in \cZ^\vee(k)=(\bbF^{\times}_q)^{\vee}\times k^{\times}$. Denote by $(V^{(1)}_{\mathbf{\whT},0}/W_0)_{z_2}$ the fibre of $\pr_{\bbG_m}$ at $z_2\in k^{\times}$. Then by \ref{connected_comp} we have
$$
(V^{(1)}_{\mathbf{\whT},0}/W_0)_{\zeta}=\coprod_{\gamma\in(\bbT^{\vee}/W_0)_{\reg},  \gamma |_ {\bbF^{\times}_q}= \zeta |_ {\bbF^{\times}_q} } V_{\mathbf{\whT},0,z_2}\coprod_{\gamma\in(\bbT^{\vee}/W_0)_{\nonreg}, \gamma |_ {\bbF^{\times}_q}= \zeta |_ {\bbF^{\times}_q} } V_{\mathbf{\whT},0,z_2}/W_0.
$$
Recall that the choice of standard coordinates $x,y$ identifies
$$
V_{\mathbf{\whT},0,z_2}\simeq \bbA^1\cup_0 \bbA^1
$$ 
with two affine lines over $k$, intersecting at the origin, cf. \ref{regularcase}. 
On the other hand, the choice of the Steinberg coordinate $z_1$ identifies
$$
V_{\mathbf{\whT},0,z_2}/W_0\simeq \bbA^1$$ 
with a single affine line over $k$, cf. \ref{non-regularcase}. 
\end{Pt*}

\begin{Lem*}\label{twistSatake} Let $\zeta, \eta\in \cZ^{\vee}$. The action of $\eta$ on 
$V^{(1)}_{\mathbf{\whT},0}/W_0$ induces an isomorphism of $k$-schemes $(V^{(1)}_{\mathbf{\whT},0}/W_0)_{\zeta} \simeq (V^{(1)}_{\mathbf{\whT},0}/W_0)_{\zeta\eta^2}$.
\end{Lem*} 

\begin{proof} 
Follows from the last formula in \ref{projections}. 
\end{proof}

\begin{Pt*}\label{twistHecke} Recall from \ref{connected_comp}  the spherical map 
 $$
\xymatrix{
\Sph:  (V^{(1)}_{\mathbf{\whT},0}/W_0)(k)\ar[r] & \{\textrm{left $\cH^{(1)}_{\overline{\bbF}_q}$-modules}\}/\sim.
}
$$
The $\cH^{(1)}_{\overline{\bbF}_q}$-modules in the image of this map are of length $1$ or $2$, cf. \ref{Sphreg} and \ref{Sphnonreg}.
We write $\Sph(v)^{\rm ss}$ for the semisimplification of the module $\Sph(v)$, for $v\in (V^{(1)}_{\mathbf{\whT},0}/W_0)(k)$.

Let $(\omega^r,z_0)\in \cZ^\vee(k)$. Recall that the standard or irreducible $\cH^{(1)}_{\overline{\bbF}_q}$-modules 
may be `twisted by the character $(\omega^r,z_0)$' : in the regular case, the actions of $X,Y,U^2$ get multiplied by $z_0,z_0,z_0^{2}$ respectively and the component $\gamma$ gets multiplied by $\omega^r$, cf. \cite[2.4]{V04}; in the non-regular case, the action of $U$ gets multiplied by $z_0$, the action of $S$ remains unchanged and the component $\gamma$ gets multiplied by $\omega^r$, cf. \cite[1.6]{V04}.
This gives an action of the group of $k$-points of $\cZ^\vee$ on the standard or irreducible $\cH^{(1)}_{\overline{\bbF}_q}$-modules.
It extends to an action on semisimple $\cH^{(1)}_{\overline{\bbF}_q}$-modules.

\end{Pt*}
\begin{Prop*} \label{Asphequiv} The map $\Sph(-)^{\rm ss}$ is $\cZ^\vee(k)$-equivariant. 
\end{Prop*}
\begin{proof} Let $(\omega^r,z_0)\in \cZ^\vee(k)$.
Let $v\in  (V^{(1)}_{\mathbf{\whT},0}/W_0)(k)$ and let its connected component be indexed by $\gamma\in\bbT^{\vee} /W_0$. 
Suppose that $\gamma$ is regular, choose an ordering $\gamma=(\chi,\chi^s)$ on the set $\gamma$ and standard coordinates. Then
$\Sph(v)=\Sph^{\gamma}(v)$ is a simple two-dimensional standard $\cH^{\gamma}_{\overline{\bbF}_p}$-module, cf.  \ref{Sphreg}, i.e. of the form
$M(x,y,z_2,\chi)$ \cite[3.2]{V04}. Then $$\Sph(v.(\omega^r,z_0))\simeq M(z_0x,z_0y,z_0^2z_2,\chi. \omega^r)\simeq\Sph(v).(\omega^r,z_0).$$
Suppose that $\gamma=\{\chi\}$ is non-regular and choose Steinberg coordinates. (a) If $v\in D(2)_{\gamma}(k)$, then 
$\Sph(v)=\Sph^{\gamma}(2)(v)$ is 
a simple two-dimensional $\cH^{\gamma}_{\overline{\bbF}_p}$-module, cf.  \ref{Sphnonreg}, i.e. of the form 
$M(z_1,z_2,\chi)$ \cite[3.2]{V04}. 
Then $$\Sph(v.(\omega^r,z_0))\simeq M(z_0z_1,z_0^2z_2,\chi. \omega^r)\simeq \Sph(v).(\omega^r,z_0).$$
(b) If  $v\in D(1)_{\gamma}(k)$, then the semisimplified module $\Sph(v)^{\sss}$ is the direct sum of the two characters in the spherical pair 
$\Sph^{\gamma}(1)(v)=\{(0,z_1),(-1,-z_1)\}$ where $z_2=z^2_1$. Similarly $\Sph(v.(\omega^r,z_0))^{\sss}$ is the direct sum of the 
characters $\{(0,z_0z_1),(-1,-z_0z_1)\}$ in the component $\gamma.\omega^r$, and hence is isomorphic to $\Sph(v)^{\sss}.(\omega^r,z_0)$.
\end{proof}

\begin{Pt*}  We now explain the compatibility with central characters for $G$-representations.
In order to do this, let us consider $W$ to be a subgroup of $G$, by sending $s$ to the matrix 
$ \left(\begin{array}{cc}
0 & 1\\
1 & 0
\end{array} \right)
$
and by identifying the group $\Lambda$ with a subgroup of $T$ via $(1,0)\mapsto\diag(\varpi^{-1},1)$ and 
$(0,1)\mapsto\diag(1,\varpi^{-1})$. We obtain for example (recall that $u=(1,0)s\in W$)

$$
u =
\left (\begin{array}{cc}
0 & \varpi^{-1}\\
1 & 0
\end{array} \right),\hskip15pt
u^{-1} =
\left (\begin{array}{cc}
0 & 1\\
\varpi & 0
\end{array} \right),\hskip15pt
us =
\left (\begin{array}{cc}
\varpi^{-1} & 0\\
0 & 1
\end{array} \right),\hskip15pt
su =
\left (\begin{array}{cc}
1 & 0\\
0 & \varpi^{-1}
\end{array} \right).
$$
Moreover, $u^{2}=\diag(\varpi^{-1},\varpi^{-1})$.\footnote{Note that our element $u$ equals the element $u^{-1}$ in \cite{Be11},\cite{Br07} and \cite{V04}.} Since 
$$
\left (\begin{array}{cc}
0 & \varpi^{-1}\\
1 & 0
\end{array} \right) 
\left (\begin{array}{cc}
a & b\\
c & d
\end{array} \right)
\left (\begin{array}{cc}
0 & 1\\
\varpi & 0
\end{array} \right) = \left (\begin{array}{cc}
d & \varpi^{-1}c\\
\varpi b & a
\end{array} \right)
$$
the element $u\in G$ normalizes 
the group $I^{(1)}$.
\end{Pt*}

\begin{Pt*} \label{I(1)invariants}\label{some_invariants} Let $\Mod^{\rm sm}(k[G])$ be the category of smooth $G$-representations over $k$. 
Taking $I^{(1)}$-invariants yields a functor $\pi\mapsto \pi^{I^{(1)}}$ from $\Mod^{\rm sm} (k[G])$ to the category $\Mod(\cH^{(1)}_{\overline{\bbF}_q})$. If $F=\bbQ_p$, it induces a bijection between the irreducible $G$-representations and the irreducible
$\cH^{(1)}_{\overline{\bbF}_p}$-modules, under which supersingular representations correspond to supersingular Hecke modules \cite{V04}. 

\vskip5pt

For future reference, let us recall the $I^{(1)}$-invariants for some classes of representations. 
 If $\pi=\Ind_B^G(\chi)$ is a principal series representation with $\chi=\chi_1\otimes\chi_2$, then $\pi^{I^{(1)}}$ is a standard module in the component $\gamma:=\{\chi|_{\bbT}, \chi^s|_{\bbT}\}$. 
 
In the regular case, one chooses the ordering $(\chi|_{\bbT}, \chi^s|_{\bbT})$ on the set $\gamma$ and standard coordinates $x,y$. Then 
$$\Ind_B^G(\chi)^{I^{(1)}}=M(0,\chi(su),\chi (u^2),\chi|_{\bbT})=M(0,\chi_2(\varpi^{-1}),\chi_1(\varpi^{-1})\chi_2(\varpi^{-1}),\chi|_{\bbT})$$ 
In the non-regular case, one has 
$$\Ind_B^G(\chi)^{I^{(1)}}=M(\chi(su),\chi(u^2),\chi|_{\bbT})=M(\chi_2(\varpi^{-1}),\chi_1(\varpi^{-1})\chi_2(\varpi^{-1}),\chi|_{\bbT}).$$ 
These standard modules are irreducible if and only if $\chi\neq\chi^s$ \cite[4.2/4.3]{V04}.\footnote{Our formulas differ from \cite[4.2/4.3]{V04} by $\chi(\cdot)\leftrightarrow \chi(\cdot)^{-1}$, since we are working with left modules; also compare with the explicit
calculation with right convolution given in \cite[Appendix A.5]{V04}.}

\vskip5pt

Let $F=\bbQ_p$. If $\pi=\pi(r,0,\eta)$ is a standard supersingular representation with parameter $r=0,...,p-1$ and central character $\eta: \bbQ_p^{\times}\ra k^{\times}$, then $\pi^{I^{(1)}}$ is a supersingular module in the component $\gamma=\{\chi,\chi^s\}$ represented by the character $\chi:=(\omega^r\otimes 1)\cdot(\eta |_{\bbF_p^{\times}})$, cf. \cite[5.1/5.3]{Br07}. If $\pi$ is the trivial representation $\mathbbm{1}$ or the Steinberg representation {\rm St}, then $\gamma=1$ and $\pi^{I^{(1)}}$ is the character $(0,1)$ or $(-1,-1)$ respectively. 
\end{Pt*}

\begin{Pt*}\label{translate_the_action}
Let $\pi\in \Mod^{\rm sm} (k[G])$. Since $u\in G$ normalizes 
the group $I^{(1)}$, one has $I^{(1)}u I^{(1)}= u I^{(1)}$. It follows that the convolution action of the Hecke operator $U$ (resp. $U^2$) on $\pi^{I^{(1)}}$ is therefore induced by the action of $u$ (resp. $u^2$ on $\pi$). 
Similarly, the group $I^{(1)}$ is normalized by the Iwahori subgroup $I$ and $I/I^{(1)}\simeq \bbT$. It follows that
the convolution action of the operators $T_t, t\in \bbT$ on $\pi^{I^{(1)}}$ is the factorization of the $\mathbf{T}(o_F)$-action on $\pi$.

\end{Pt*}

\begin{Pt*}\label{ZGveeZvee}
We identify $F^{\times}$ with the center $Z(G)$ via $a\mapsto\diag(a,a)$. A (smooth) character 
$$
\zeta:  Z(G)=F^{\times}\lra k^\times
$$ 
is determined by its value $\zeta(\varpi^{-1})\in k^\times$ and its restriction $\zeta |_{o_F^{\times}}$. Since the latter is trivial on the subgroup $1+\varpi o_F$, we may view it as a character of $\bbF_q^{\times}$; we will write $\zeta |_{\bbF_q^{\times}}$ for this restriction in the following. Thus the group of characters of $Z(G)$ gets identified with the group of $k$-points of the group scheme $\cZ^{\vee} =(\bbF_q^{\times})^{\vee}\times \bbG_m$: 
$$Z(G)^\vee\iso \cZ^{\vee}(k),\; \zeta\mapsto (\zeta |_{\bbF_q^{\times}},\zeta(\varpi^{-1})).$$
\end{Pt*}

\begin{Prop*}\label{central_car_comp} 
Suppose that $\pi\in \Mod^{\rm sm} (k[G])$ has a central character $\zeta: Z(G)\rightarrow k^{\times}$. Then the Satake parameter $S(\pi^{I^{(1)}})$ of $\pi^{I^{(1)}}\in\Mod(\cH^{(1)}_{\overline{\bbF}_q})$ has central character $\zeta$, i.e. it is supported on the closed subscheme
$$
(V^{(1)}_{\mathbf{\whT},0}/W_0)_{(\zeta |_ {\bbF^{\times}_q} , \zeta(\varpi^{-1}))}\subset V^{(1)}_{\mathbf{\whT},0}/W_0.
$$
\end{Prop*}

\begin{proof} 
If $M$ is any $\cH^{(1)}_{\overline{\bbF}_q}$-module, then
$$
M=\bigoplus_{\gamma\in\bbT^{\vee}/W_0}\varepsilon_{\gamma}M=\bigoplus_{\gamma\in\bbT^{\vee}/W_0}\oplus_{\lambda\in\gamma}\varepsilon_{\lambda}M,
$$
and $\bbT\subset \overline{\bbF}_q[\bbT]\subset \cH^{(1)}_{\overline{\bbF}_q}$ acts on $\varepsilon_{\lambda}M$ through the character
$\lambda:\bbT\ra \bbF_q^{\times}$. Now if $M=\pi^{I^{(1)}}$, then the $\bbT$-action on $M$ is the factorization of the  $\mathbf{T}(o_F)$-action on $\pi$, cf. \ref{translate_the_action}. In particular, the restriction of the $\bbT$-action along the diagonal inclusion $\bbF_q^{\times}\subset\bbT$ is the factorization of the action of the central subgroup $o^{\times}_F\subset Z(G)$ on $\pi$, which is given by $\zeta|_{o^{\times}_F}$ by assumption. Hence
$$
\varepsilon_{\gamma}M\neq 0\quad\Longrightarrow\quad \forall\lambda\in\gamma,\ \lambda|_{\bbF_q^{\times}}=\zeta|_{\bbF_q^{\times}}\ \textrm{i.e.}\ \gamma|_{\bbF_q^{\times}}=\zeta|_{\bbF_q^{\times}}.
$$
Moreover, the element $u^2=\diag(\varpi^{-1},\varpi^{-1})\in Z(G)$ acts on $\pi$ by multiplication by $\zeta (\varpi^{-1})$ by assumption. Therefore, by \ref{translate_the_action}, the Hecke operator $z_2:=U^2\in \cH^{(1)}_{\overline{\bbF}_q}$ acts on $\pi^{I^{(1)}}$ by multiplication by $\zeta (\varpi^{-1})$. Thus we have obtained that $S(\pi^{I^{(1)}})$ is supported on 
$$
\coprod_{\gamma\in(\bbT^{\vee}/W_0)_{\reg},  \gamma |_ {\bbF^{\times}_q}= \zeta |_ {\bbF^{\times}_q} } V_{\mathbf{\whT},0,\zeta (\varpi^{-1})}\coprod_{\gamma\in(\bbT^{\vee}/W_0)_{\nonreg}, \gamma |_ {\bbF^{\times}_q}= \zeta |_ {\bbF^{\times}_q} } V_{\mathbf{\whT},0,\zeta (\varpi^{-1})}/W_0
=(V^{(1)}_{\mathbf{\whT},0}/W_0)_{(\zeta |_ {\bbF^{\times}_q} , \zeta(\varpi^{-1}))}.
$$
\end{proof}

\vskip10pt 

\noindent {\small Cédric Pépin, LAGA, Université Paris 13, 99 avenue Jean-Baptiste Clément, 93 430 Villetaneuse, France \newline {\it E-mail address: \url{cpepin@math.univ-paris13.fr}} }

\vskip10pt

\noindent {\small Tobias Schmidt, IRMAR, Universit\'e de Rennes 1, Campus Beaulieu, 35042 Rennes, France \newline {\it E-mail address: \url{tobias.schmidt@univ-rennes1.fr}} }

\end{document}